\documentclass[10pt,leqno, final]{amsart}
\usepackage[colorlinks,linkcolor=blue,citecolor=blue,urlcolor=blue]{hyperref}
\usepackage{amsfonts,amssymb}
\usepackage{color}
\usepackage{graphicx}
\usepackage{mathtools}
\usepackage{extarrows}
\usepackage{algorithmic}   
\usepackage{caption}
\usepackage{subcaption}
\usepackage{tikz}

\usepackage{todonotes}
\usepackage{marginnote}
\usepackage[normalem]{ulem}

\allowdisplaybreaks

\usepackage[notcite]{showkeys}

\newtheorem{theorem}{Theorem}[section]
\newtheorem{lemma}[theorem]{Lemma}
\newtheorem{corollary}[theorem]{Corollary}
\newtheorem{proposition}[theorem]{Proposition}
\newtheorem{assumption}[theorem]{Assumption}

\theoremstyle{definition}
\newtheorem{definition}[theorem]{Definition}
\newtheorem{example}[theorem]{Example}
\newtheorem*{example*}{Example}
\newtheorem{remark}[theorem]{Remark}

\numberwithin{equation}{section}

 \usepackage[
	backend=biber,
	style=numeric,
	maxbibnames=99,
	giveninits=true,
	sorting=nyt,
	sortcites=true,
	isbn=false,
	citecounter=true,
]{biblatex}
\newcommand{\BB}{\mathcal{B}}

\newcommand{\EE}{\mathcal{E}}
\newcommand{\FF}{\mathcal{F}}

\newcommand{\HH}{\mathcal{H}}
\newcommand{\II}{\mathcal{I}}
\newcommand{\JJ}{\mathcal{J}}

\newcommand{\LL}{\mathcal{L}}

\newcommand{\RR}{\mathcal{R}}

\newcommand{\VV}{\mathcal{V}}
\newcommand{\WW}{\mathcal{W}}
\newcommand{\XX}{\mathcal{X}}

\newcommand{\ZZ}{\mathcal{Z}}

\newcommand{\N}{\mathbb{N}}

\newcommand{\R}{\mathbb{R}}

\newcommand{\frakM}{\mathfrak{M}}
\RequirePackage{stmaryrd} 

\newcommand{\embed}{\hookrightarrow}

\newcommand{\weakly}{\rightharpoonup}
\newcommand{\weak}{\rightharpoonup}

\newcommand{\dist}{\operatorname{dist}}

\renewcommand{\d}{\mathrm{d}}

\newcommand{\dual}[2]{\langle #1 , #2 \rangle}

\def\clap#1{\hbox to 0pt{\hss#1\hss}}

\newcommand{\BV}{\mathrm{BV}}

\newcommand{\Drm}{\mathrm{D}}
\newcommand{\frakP}{\mathfrak{P}}
\newcommand{\frakC}{\mathfrak{C}}
\newcommand{\var}{\mathrm{Var}}

\begin{document}

\title[Optimal Control of Rate-Independent Systems]
{Viscous Approximation of Optimal Control Problems Governed by Rate-Independent Systems with Non-Convex Energies}

\author{Merlin Andreia} \address{Technische Universit\"at Dortmund, Fakult\"at f\"ur
  Mathematik, Lehrstuhl LSX, Vogelpothsweg 87, 44227 Dortmund, Germany}
\email{merlin.andreia@tu-dortmund.de}

\author{Christian Meyer} \address{Technische Universit\"at Dortmund, Fakult\"at f\"ur
  Mathematik, Lehrstuhl LSX, Vogelpothsweg 87, 44227 Dortmund, Germany}
\email{christian2.meyer@tu-dortmund.de}


\date{\today} 


\begin{abstract} 
    We consider an optimal control problem governed by a rate-inde\-pendent system with non-convex energy. 
    The state equation is approximated by means of viscous regularization w.r.t.\ to hierarchy of two different 
    Hilbert spaces. The regularized problem corresponds to an optimal control problem subject to a 
    non-smooth ODE in Hilbert space, which is substantially easier to solve than the original optimal control problem. 
    The convergence properties of the viscous regularization are investigated. It is shown that every sequence 
    of globally optimal solutions of the viscous problems admits a (weakly) converging subsequence whose limit 
    is a globally optimal solution of the original problem, provided that the latter admits at least one optimal solution
    with an optimal state that is continuous in time.
\end{abstract}

\maketitle

\section{Introduction}

This paper is concerned with the optimal control of a rate-independent system of the form 
\begin{equation}\label{eq:ris}
    0 \in \partial \RR(\dot z(t)) + D_z \II(\ell(t), z(t)) \quad \text{f.a.a.\ } t \in (0,T), 
    \quad z(0) = z_0,
\end{equation}
where $\RR: \XX \to \R$ is convex and positive 1-homogeneous on a Banach space $\XX$.
In addition to $\XX$, three are two Hilbert spaces $\VV$ and $\ZZ$ such that $\ZZ \embed^{c,d} \VV \embed \XX$.
Moreover, $\II$ denotes the energy given by 
\begin{equation*}
    \II(\ell(t), z) := \frac{1}{2}\,\dual{A z}{z} + \FF(z) - \dual{\ell(t)}{z},
\end{equation*}
where $A \in L(\ZZ,\ZZ^*)$ is a symmetric and coercive operator and $\FF: \ZZ \to \R$ is smooth, but potentially non-convex. 
The external load $\ell: [0,T] \to \VV^*$ serves as control to optimize the system. 
The precise assumptions on the data are specified in Section~\ref{sec:assumc3} below.
An example in function space satisfying all assumptions is given at the end of Section~\ref{sec:assumc3}.

Systems of the form \eqref{eq:ris} involving a positive homogeneous dissipation functional and a 
non-convex energy are challenging to handle as ``classical'' time-continuous solutions satisfying the system
a.e.\ in time do in general not exist. We refer to \cite[Section~2.3]{Tho22} for a counterexample.
For that reason, several alternative solution concepts have been developed, among them 
global energetic solutions and (parametrized) balanced viscosity (BV) solutions as the probably most 
prominent examples. We refer to \cite{MR15} for an overview over the numerous solution concepts.
Here we focus on the concept of parametrized BV solutions. As the name already indicates, 
the idea behind this concept is to approximate solutions to \eqref{eq:ris} by means of viscous regularization.
The passage to the limit for vanishing viscosity is accompanied by a parametrization of the solution 
trajectory and has been carried out by various authors before, see, e.g., \cite{MRS12, MRS16}. 
Here, we aim to adapt this procedure to an optimal control problem governed by \eqref{eq:ris}. 
The precise definition of a parametrized BV solution will be given in Section~\ref{sec:defparaBV} below.

Formally, the optimal control problem under consideration reads
\begin{equation}\tag{OCP}\label{eq:ocp0}
	\left\{\quad  
    \begin{aligned}
        \min \quad &  J(z,\ell) := j(z(T))+\frac{\beta}{2} \, \|\ell\|_{H^1(0,T;\VV^*)}^2\\
		    \text{s.t.} \quad & \text{$z$ is a (parametrized) BV solution of \eqref{eq:ris} associated with $\ell$},
    \end{aligned}
    \right.
\end{equation}
where $j : \VV \to \R$ is a given objective functional and $\beta>0$ is a fixed Tikhonov regularization parameter.
The optimal control problem will be completed by additional constraints 
that ensure local stability of the optimal solution at initial and end time. 
The motivation for these additional constraints will become clear when introducing the notion
of parametrized BV solutions and we will introduce them rigorously in Section~\ref{sec:optcontrol} below.
Inspired by the regularization employed in \cite{MRS12, MRS16}, the viscous approximation of \eqref{eq:ocp0} reads 
\begin{equation}\tag{$\mathrm{vOCP}_\epsilon$} 
    \left\{ \quad
    \begin{aligned}\label{eq:ocpeps}
        \min \quad &  J(z_\epsilon,\ell):=j(z_\epsilon(T))+\frac{\beta}{2} \, \|\ell\|_{H^1(0,T;\VV^*)}^2\\
        \text{s.t.} \quad & 0\in\partial\RR(z_\epsilon '(t)) + \epsilon\, z_\epsilon'(t) + \mathrm{D}_z\II(\ell(t),z_\epsilon(t)),
        \quad z_\epsilon(0)=z_0,
   \end{aligned}
    \right.
\end{equation}
where $\epsilon > 0$ is the regularization parameter. Compared to \eqref{eq:ocp0},
\eqref{eq:ocpeps} provides several advantageous features, for instance, the state equation is 
uniquely solvable for every control $\ell$ in contrast to \eqref{eq:ris}.
Nonetheless, it is still a non-trivial problem, as, for instance, the existence of solutions to the regularized 
state equation is not immediate to see and is frequently shown by means of time discretization \cite{MRS13, Tho22}
or using a second regularization, which turns the state equation into an ODE in the Hilbert space $\ZZ$ \cite{KRZ11}. 
Here, we pursue the second approach and consider a second ``double'' viscous regularization of \eqref{eq:ocp0} 
given by
\begin{equation}\tag{$\mathrm{vOCP}_{\epsilon, \delta}$} 
    \left\{ \quad
    \begin{aligned}\label{eq:ocpepsdelta}
        \min \quad &  J(z_{\epsilon, \delta},\ell):=j(z_{\epsilon, \delta}(T))+\frac{\beta}{2} \, \|\ell\|_{H^1(0,T;\VV^*)}^2\\
        \text{s.t.} \quad & 0\in\partial\RR(z_{\epsilon, \delta} '(t)) + \epsilon\, z_{\epsilon, \delta}'(t) 
        + \delta\, A  z_{\epsilon, \delta}' + \mathrm{D}_z\II(\ell(t),z_{\epsilon, \delta}(t)), \\
        & z_{\epsilon, \delta}(0)=z_0
    \end{aligned}
    \right.
\end{equation}
with a second regularization parameter $\delta > 0$. Since $A$ induces a norm on $\ZZ$, 
the second regularization is nothing else than a viscous regularization in $\ZZ$, while the first regularization 
in \eqref{eq:ocpeps} is a regularization in $\VV$.
From the viewpoint of optimal control, this second regularization 
has the additional advantage that it leads to an optimal control problem governed by a non-smooth ODE. 
Problems of this type are meanwhile well investigated in literature, we only refer to \cite{SWW17}, where a problem 
of exactly this structure is analyzed in detail by smoothing the non-smooth ODE.
Similarly to \eqref{eq:ocp0}, both regularized optimal control problems will be completed by 
additional constraints at initial and end time, see Section~\ref{sec:optcontrol} for a mathematically
rigorous statement of the optimal control problems.

Let us put our work into perspective. 
Optimal control problems governed by rate-independent systems with uniformly convex energies 
have been intensively analyzed in the literature, especially optimal control of the sweeping process 
is well investigated, see, e.g., \cite{Bro87, Wac12, BK13, AO14, Wac15, Wac16, CHNM16, GW18, AC18, BC23} 
and the references therein.
The situation changes, if one turns to non-convex energies. 
While several works prove the existence of optimal solutions \cite{Ste12, ELS13, EL14}, only little research 
has been done concerning their approximation by means of regularization. 
In \cite{MW21}, the viscous regularization of optimal control of perfect plasticity is investigated. 
In this case the energy is not uniformly convex, but still convex. 
In \cite{Rin08, Rin09, MR09}, approximation of optimal control of \eqref{eq:ris} by means of time discretization is 
investigated. 
The work that is closest to ours is \cite{KMS22}, where the same setting as in our case is discussed 
with the only difference that the underlying spaces are all the same and equal $\R^n$.
This assumption substantially simplifies the convergence analysis and makes the second regularization 
in \eqref{eq:ocpepsdelta} superfluous.
The main challenge of the convergence analysis is the so-called reverse approximation, which is nothing else than
the construction of a recovery sequence that is feasible for the regularized problem 
\eqref{eq:ocpeps} and converges (in a suitable sense) to a solution of the original problem \eqref{eq:ocp0}. 
At this point, our convergence analysis in principle follows the strategy of \cite{KMS22} with the difference that 
the second regularization is required for the construction of the recovery sequence. 
This construction benefits from the fact that, in the context of optimal control, 
we do not only have the state $z$, but also the control $\ell$ at our disposal. 
If the loads are fixed, there is in general no hope that every BV solution can be approximated by means of viscous 
regularization, cf.\ the one-dimensional counterexample in \cite[Example~1.8.3]{MR15}. 
In our context however, the external loads are not fixed and can be used for the construction of a recovery sequence. 
However, this sequence has to meet the regularity requirements of the loads, which are functions from $H^1(0,T;\VV^*)$
due to the Tikhonov term in \eqref{eq:ocp0}, \eqref{eq:ocpeps}, and \eqref{eq:ocpepsdelta}. 
This enforces the control to be continuous in time and the investigations in \cite{AM24} illustrate that this is 
indispensable in the context of optimal control, due to a lack of compactness of the solution set for 
external loads that are not continuous in time.
In order to guarantee this regularity of the loads within our construction of the recovery sequence, 
we unfortunately need the additional assumption that 
there exists at least one optimal solution whose state is continuous in time. 
Such a rather restrictive assumption is typical for the reverse approximation property, see also \cite{MW21, KMS22} in this context.

The paper is organized as follows: 
In Section~\ref{sec:assumc3}, we state our standing assumptions on the quantities in \eqref{eq:ris} and 
give a function space example at the end of the section, which fulfills all of our assumptions.
The notion of paramterized BV solutions in rigorously defined in Section~\ref{sec:defparaBV} and we 
state a known convergence result for the viscous regularization of the state equation alone.
Section~\ref{sec:visreg} is then devoted to the double viscous regularization and to the convergence analysis 
for the state equation for $\delta \searrow 0$ and fixed $\varepsilon > 0$. 
We show that solutions of the double viscously regularized system converge to solutions of the single regularization 
and in this way also prove existence of solutions of the latter. This result however is merely a byproduct and 
was already well known before. What is more important, we derive an estimate in Lemma~\ref{lem:z''} which 
only holds for the double viscous regularization and 
is essential for the reverse approximation argument in Section~\ref{sec:revapp}. 
Before we address this issue, we prove the existence of optimal solutions to \eqref{eq:ocp0}--\eqref{eq:ocpepsdelta}.
While the proofs are more or less standard in case of \eqref{eq:ocpeps} and \eqref{eq:ocpepsdelta}, 
the proof of existence for \eqref{eq:ocp0} relies on a deep compactness result of \cite{KT23}.
The reverse approximation property is then established in Section~\ref{sec:revapp}. 
Using this result, we finally prove our main results in Section~\ref{sec:ocpapprox}, 
namely that optimal solutions to \eqref{eq:ocp0} 
can be approximated by solutions of \eqref{eq:ocpeps} and \eqref{eq:ocpepsdelta}, respectively, in the sense
that sequence of optimal solutions to \eqref{eq:ocpeps} resp.\ \eqref{eq:ocpepsdelta} 
admit (weakly) convergent subsequences and the limits are solutions to \eqref{eq:ocp0}, 
see Theorems~\ref{thm:ocpapproxeps} and \ref{thm:ocpapproxepsdelta}.
The paper ends with an appendix containing auxiliary results.

\section{Standing Assumptions}\label{sec:assumc3}

Let us introduce the notation and our standing assumptions that are valid throughout the entire paper. 

\subsubsection*{Notation}
If $X$ is a normed space, then $B_X(0,r)$ denotes the open ball in $X$ around $0$ with radius $r$. 
For the distance of an element $x \in X$ to a subset $M\subset X$, we write 
$\dist_X(x, M) := \inf_{m\in M} \|x - m\|_X$. If $x \notin X$, then $\dist_X(x, M)$ is set to $\infty$.

Given $T> 0$ and a function $v: [0,T] \to X$, we denote its variation by 
\begin{equation*}
    \var_X(v) := \sup\Big\{ \sum_{i=1}^n 
    \|v(t_i) - v(t_{i-1})\|_X \colon n\in \N,\; 0 \leq t_0 < t_1 < \ldots < t_n \leq T\Big\}. 
\end{equation*}
The set of functions with bounded variation is denoted by $\BV([0,T];X)$.
Equipped with the norm $\|v\|_{\BV([0,T];X)} := \sup_{t\in [0,T]} \|v(t)\|_X + \var_X(v)$, 
it becomes a Banach space, provided that $X$ is complete. 
By $\frakM([0,T];X)$, we denote the space of vector-valued regular Borel measures, i.e., 
the space of countably additive regular set functions from $\BB([0,T])$ to $X$, where $\BB([0,T])$ 
denotes the Borel-$\sigma$ algebra on $[0.T]$.
The space of Bochner-measurable, to the power $p\in [1, \infty]$ integrable functions from $(0,T)$ to $X$ 
is denoted by $L^p(0,T;X)$. If such a function admits a distributional time derivative in $L^p(0,T;X)$, it belongs to
the Bochner-Sobolev space $W^{1,p}(0,T;X)$.
Similarly, $C^m([0,T];X)$, $m\in \N_0$, is the space of $m$-times continuously differentiable functions with values in $X$. 

Given two normed vector spaces $X$ and $Y$, we denote the space of linear and bounded operators 
from $X$ to $Y$ by $\LL(X,Y)$. The topological dual of $X$ is denoted by $X^*$ and, for the dual 
pairing of $x\in X$ and $x^* \in X^*$, we write $\dual{x^*}{x}_{X^*, X}$. If $X$ is a Hilbert space, 
we write $\dual{\cdot }{\cdot}_X$ for the scalar product.

Finally, throughout the paper, $C>0$ denotes a generic non-negative constant. 

\subsubsection*{Spaces} In all what follows, $\ZZ$ and $\VV$ are Hilbert spaces and $\XX$ is a Banach space such that 
$\ZZ\hookrightarrow ^{c,d}\VV\hookrightarrow \XX$, where the notation $\hookrightarrow ^{c,d}$ refers to dense and compact embedding.
Moreover, there is a further intermediate Banach space $\WW$ with $\ZZ\hookrightarrow^c\WW\hookrightarrow \VV$, which 
will be essential for the treatment of the nonlinearity $\FF$, cf.\ \eqref{eq:D2FWcont} below.

\subsubsection*{Energy}\label{sec:energyass}
We assume that the energy $\II:\VV^*\times \ZZ\to\R$ has the following structure:
\begin{equation*}
    \II(\ell,z) := \frac{1}{2}\,\langle Az,z\rangle_{\ZZ^*,\ZZ} +\FF(z)-\langle\ell,z\rangle_{\VV^*,\VV}.
\end{equation*}
Herein, $A\in\LL(\ZZ,\ZZ^*)$ is supposed to be symmetric and coercive, i.e., it exists $\alpha>0$ with
\begin{align}\label{eq:Acoer}
\langle Az,z\rangle_{\ZZ^*,\ZZ}\geq \alpha\|z\|_\ZZ^2 \quad\text{ and } \quad\langle Az,w\rangle_{\ZZ^*,\ZZ}= \langle Aw,z\rangle_{\ZZ^*,\ZZ}
\end{align}
for all $z,w\in\ZZ$. Since $A\in\LL(\ZZ,\ZZ^*)$ is coercive and continuous, $|\cdot|_\ZZ:\ZZ\to \R$, given by 
$|z|_\ZZ:=\big(\langle Az,z\rangle_{\ZZ^*,\ZZ}\big)^{\frac{1}{2}}$, defines a norm on $\ZZ$ that is equivalent to the canonical norm $\|\cdot\|_\ZZ$.
Furthermore, we assume that $\FF:\ZZ\to\R$ satisfies
\begin{align}
	&\quad\FF\in C^2(\ZZ;\R),~\FF\geq 0, \label{eq:F1}\\
	&\Drm_z\FF\in C^1(\ZZ;\VV^*),\quad \|\Drm_z^2\FF(z)v\|_{\VV^*}\leq \gamma(1+\|z\|_\ZZ^q)\|v\|_\ZZ~\forall v,z\in\ZZ \label{eq:F2} 
\end{align}  
with some $q \in [0, 1/2)$ and a constant $\gamma > 0$. 

\begin{remark}
    Admittedly, the growth condition on $\Drm_z^2\FF$ with an exponent $q< 1/2$ is rather restrictive, but we emphasize 
    that on bounded sets large curvature of $\FF$ is well allowed.  Moreover, the condition $\FF\geq 0$ can be relaxed 
    by just requiring a lower bound. Nonetheless, for simplicity, we just set this lower bound to zero.
\end{remark}

From these assumptions, it follows that $\FF$ is completely continuous:

\begin{lemma}\label{lem:convF}
    Let $z_k \weakly z$ in $\ZZ$. Then $\FF(z_k) \to \FF(z)$ as $k \to \infty$.
\end{lemma}

\begin{proof}
    Applying the mean value theorem yields the existence of $\tau_k\in(0,1)$ with
    \begin{align*}
	    &\|\FF(z_k)-\FF(z)\|_\ZZ\\
	    &\; =|\langle \Drm_z\FF(z),z_k-z\rangle_{\VV^*,\VV}|
	    +\frac{1}{2}\big|\langle \Drm_z^2\FF(z+\tau_k(z_k-z))(z_k-z),z_k-z\rangle_{\VV^*,\VV}\big|\\
	    & \; \leq \|\Drm_z\FF(z)\|_{\VV^*}\|z_k-z\|_\VV+\frac{1}{2}C\big(1+\|z+\tau_k(z_k-z)\|_\ZZ^q\big)\|z_k-z\|_\ZZ\|z_k-z\|_\VV
	    \to 0,
    \end{align*}
    where we used \eqref{eq:F2}, the boundedness of $(z_k)_k\subset \ZZ$ by its weak convergence,
    and the strong convergence $z_k\to z$ in $\VV$ by the compact embedding $\ZZ\hookrightarrow^c\VV$.
\end{proof}

Moreover, $\Drm_z\FF:\ZZ\to\ZZ^*$ is supposed to be weak-weak continuous, i.e.,
\begin{align}\label{DFweak}
	z_k\rightharpoonup z \text{ in } \ZZ \quad \implies \quad \Drm_z\FF(z_k)\rightharpoonup \Drm_z\FF(z) \text{ in}~\ZZ^*.
\end{align}
We moreover suppose that $D_z^2\FF$ can be extended to 
the intermediate space $\WW$ with values in $\LL(\ZZ,\VV^*)$ and we assume that this extension 
is (strongly) continuous, i.e., 
\begin{equation}\label{eq:D2FWcont}
		z_k\to z \text{ in } \WW\quad \implies\quad \mathrm{D}_z^2\FF(z_k)\to \mathrm{D}_z^2\FF(z) \text{ in } \LL(\ZZ,\VV^*).
\end{equation}
Finally, we require that $\mathrm{D}_z^2\FF:\ZZ\to \LL(\ZZ,\ZZ^*)$ is Lipschitz continuous on bounded sets, 
i.e., for all $r>0$,  there exists $L_r>0$ such that
\begin{align}\label{eq:D2Flip}
	\|(\mathrm{D}_z^2\FF(z_1)-\mathrm{D}_z^2\FF(z_2))v\|_{\ZZ^*}\leq L_r \|z_1-z_2\|_\ZZ\|v\|_\ZZ 
\end{align}
for all $z_1,z_2\in B_\ZZ(0,r)$ and $v\in\ZZ$.

For later purpose, we moreover introduce the load-independent part of the energy, called reduced energy, by
\begin{equation*}
	\EE(z):= \frac{1}{2}\,\langle Az,z\rangle_{\ZZ^*,\ZZ} +\FF(z),
\end{equation*}
so that the energy can be written as $\II(\ell,z)=\EE(z)-\langle\ell,z\rangle_{\VV^*,\VV}$. 

\subsubsection*{Dissipation}
For the dissipation $\RR:\ZZ\to [0,\infty)$\index{dissipation}, we assume that
\begin{align}
	&\RR \mbox{~is~proper,~convex,~and~lower~semicontinuous}, \label{R1}\\
	&\RR \text{ is positive 1-homogeneous, i.e., }  
	\RR(\lambda z)=\lambda\RR(z) ~\forall \,z\in\ZZ,~\lambda>0, \label{R2}\\
 	&\exists \,\underline{c}, \overline{c} >0, \mbox{~such~that}~\underline{c}\,\|z\|_\XX\leq\RR(z)\leq \overline{c}\,\|z\|_\VV~\forall \,z\in\ZZ.  \label{R3}
\end{align}

\begin{lemma}\label{lem:dR(0)bound}
	Let the dissipation $\RR:\ZZ\to[0,\infty)$ comply with \eqref{R1}--\eqref{R3} and $z\in\ZZ$ be arbitrarily given. 
	Then the subdifferential $\partial \RR(z)$ is a bounded subset of $\VV^*$.
\end{lemma}

\begin{proof}
	Let $\xi\in\partial \RR(z)\subset \ZZ^*$ be given. From the convexity of $\RR$ we infer that for all $w\in\ZZ$ it holds $\RR(w)-\RR(z)\leq\RR(w-z)$ such that the definition of the subdifferential along with \eqref{R3} yields
	\begin{align*}
		\langle \xi,w-z\rangle_{\ZZ^*,\ZZ}\leq \RR(w)-\RR(v)\leq \RR(w-z)\leq \overline{c}\,\|w-z\|_\VV .
	\end{align*}
	 Thus, due to the density of $\ZZ$ in $\VV$, $\xi$ can be extended to an element of $\VV^*$ with $\|\xi\|_{\VV^*}\leq \overline{c}$, 
	 which gives the claim.
\end{proof}

\subsubsection*{Initial data}
For a given initial value $z_0\in\ZZ$ we assume that the external load $\ell\in H^1(0,T;\VV^*)$ complies with 
\begin{align}\label{eq:init}\index{initial data}
	-\mathrm{D}_z\II(\ell(0),z_0)=-Az_0-\mathrm{D}_z\FF(z_0)+\ell(0)\in \partial \RR(0),
\end{align}
i.e., the initial state is locally stable w.r.t.\ the load $\ell$.

\subsubsection*{Objective}
The state dependent part $j$ in the objective in \eqref{eq:ocp0} is supposed to be continuous from $\VV$ to $\R$ and bounded from below.
As an example, we mention $j(z):= \frac{1}{2}\,\|z-z_{\textup{des}}\|_\VV^2$ with a given desired final state $z_{\textup{des}} \in \VV$. 
Moreover, the Tikhonov parameter $\beta > 0$ is a fixed positive number.

\vspace*{2ex}

We emphasize that parts of the results of the paper do not require all of the above assumptions. 
For instance, condition~\eqref{eq:D2FWcont} involving the additional space $\WW$ is only needed for 
the reverse approximation result in Section~\ref{sec:revapp}.

The following two auxiliary results will be useful for the upcoming analysis.

\begin{lemma}\label{lem:lambdamu}
There exists constants $\lambda,\mu>0$ such that 
	\begin{align}\label{eq:assenergy1}
		\|z\|_\ZZ^2\leq \lambda \, \II(\ell,z)+\mu \, \|\ell\|_{\VV^*}^2
	\end{align}
holds for all $z\in \ZZ,\ell\in \VV^*$.
\end{lemma}

\begin{proof}
We exploit the coercivity of $A\in\LL(\ZZ,\ZZ^*)$ and non-negativity of $\FF$ by \eqref{eq:Acoer} and \eqref{eq:F1}, respectively, in order to obtain for arbitrary $z\in\ZZ$, $\ell\in\VV^*$, and $\kappa > 0$ by means of Young's inequality that 
\begin{align*}
	\frac{\alpha}{2}\|z\|_\ZZ^2
	& \leq \II(\ell,z)+\|\ell\|_{\VV^*}\|z\|_\VV \\
	& \leq \II(\ell,z)+\frac{1}{2\kappa}\|\ell\|_{\VV^*}^2+\frac{\kappa}{2}\|z\|_\VV^2
    \leq \II(\ell,z)+\frac{1}{2\kappa}\|\ell\|_{\VV^*}^2+\frac{C\kappa }{2}\|z\|_\ZZ^2,
\end{align*}
where we applied the continuous embedding $\ZZ\hookrightarrow \VV$ with embedding constant $C > 0$.
Setting $\kappa:=\frac{\alpha}{2C}$ yields the assertion with $\lambda=\frac{4}{\alpha}$ and $\mu=\frac{4C}{\alpha^2}$.
\end{proof}

\begin{lemma}\label{lem:nu}
	Let $\ell\in H^1(0,T;\VV^*)$ be a given external load.
	Then there exists a constant $\nu > 0$ such that, for all $z\in\ZZ$ and all $t\in[0,T]$, it holds that
	\begin{align}\label{eq:estenergy}
		\|z\|_\ZZ - \nu (\|\ell\|_{H^1(0,T;\VV^*)}^2 +1) \leq \II(\ell(t),z) .
	\end{align}
\end{lemma}

\begin{proof}
	By exploiting the coercivity of $A$, the non-negativity of $\FF$ and the continuous embedding 
	$\ZZ\hookrightarrow^c \VV$, we obtain for all $z\in\ZZ$ and $t\in[0,T]$ that
	\begin{equation*}
		\II(\ell(t),z)\geq \frac{\alpha}{2}\|z\|_\ZZ^2-\|\ell(t)\|_{\ZZ^*}\,\|z\|_\ZZ
		\geq \|z\|_\ZZ - \frac{(\|\ell(t)\|_{\ZZ^*}+1)^2}{2\alpha} ,
	\end{equation*}
	which, along with the continuous embedding of $H^1(0,T;\VV^*) \embed C([0,T];\ZZ^*)$ 
	and Young's inequality, directly implies the assertion.
\end{proof}

\begin{remark}\label{rem:lambdamunu}
    We underline that the constants $\lambda$, $\mu$, and $\nu$ from the above lemmas only depend
    on the coercivity constant $\alpha$ and embedding constants. This is mainly due to the non-negativity of $\FF$
    (which can be easily generalized by requiring that $\FF$ is bounded from below).
\end{remark}

We end this section with an example which fulfills all assumptions stated above.

\begin{example}\label{ex:energy}
	Let $\Omega\subset\R^3$ be an open bounded domain with Lipschitz boundary. We set 
    $\ZZ=H^1_0(\Omega)$, $\VV=L^2(\Omega)$, $\XX=L^1(\Omega)$, and $\WW = L^\kappa(\Omega)$, $2 \leq \kappa < 6$,
	such that $\ZZ\hookrightarrow^{c,d}\VV\hookrightarrow\XX$ and $\ZZ \hookrightarrow^c\WW\hookrightarrow\VV$ 
	by Sobolev embeddings. The dissipation is given by 
	\begin{align*}
		\RR:L^1(\Omega)\to \R,\quad z\mapsto \int_\Omega |z(x)| \,\d x
	\end{align*}
    so that \eqref{R1}--\eqref{R3} are indeed fulfilled. For the operator $A$ we choose 	$A=-\Delta:H^1_0(\Omega)\to H^{-1}(\Omega)$, 
    i.e., $\langle -\Delta z,v\rangle_{\ZZ^*,\ZZ}= \int_\Omega \dual{\nabla z(x)}{\nabla v(x)}_{\R^3} \,\d x$. Finally, the non-linearity is given by
	\begin{align*}
		\FF:\ZZ\to\R,\quad z\mapsto \int_\Omega f(z(x))\, \d x
	\end{align*}
	with a non-negative function $f\in C^2(\R;\R)$ fulfilling the following conditions:
	\begin{align}
        & \text{For all }	z\in H^1_0(\Omega)  \text{ there holds } f'(z)\in L^2(\Omega), \label{eq:f1}\\
        & \text{for all }	z\in L^\kappa(\Omega)  \text{ there holds } f''(z)\in L^3(\Omega), \\
        & \exists\, q\in [0,1/2),\, \gamma > 0:  \quad |f''(r)| \leq \gamma\big(1+ |r|^q\big) 
        \quad\forall\, r\in \R ,\\
        & z_k\rightharpoonup z~\text{in}~H^1(\Omega) 
        \quad\implies\quad f'(z_k)\rightharpoonup f'(z)~\text{in}~L^{\frac{6}{5}}(\Omega), \\
        & z_k\to z~\text{in}~L^\kappa(\Omega)\quad\implies\quad f''(z_k)\to f''(z) ~\text{in}~L^3(\Omega), \\
        & \exists \, s\in [1,4], L>0: \quad |f''(r_1)-f''(r_2)|\leq L\,|r_1-r_2|^s\quad \forall\, r_1,r_2\in \R, \label{eq:f6}
	\end{align}
	where, with a slight abuse of notation, we have denoted the Nemyzki operators associated with $f'$ and $f''$ by the same symbol.
	By straight forward computation, one shows that, under \eqref{eq:f1}--\eqref{eq:f6}, the non-linearity $\FF$ satisfies 
	\eqref{eq:F1}--\eqref{eq:D2Flip}.
	Concrete examples for functions satisfying \eqref{eq:f1}--\eqref{eq:f6} are the double-well type potential
    \begin{equation}\label{eq:doublewell}
        f(x) := 
        \begin{cases}
            x^4 - x^2 + 1, & x \in [-1,1]\\
            \frac{8}{a(a-2)}\,|x|^a + \frac{a - 6}{a - 2}\,x^2 + \frac{4a - 8}{a(a-2)}, & \text{else},
        \end{cases}
    \end{equation}    	
    with $a\in (2,2.5)$
	or a globally Lipschitz continuous function like $f(x) := \sin(x) + 1$. 
	In the former example, one chooses $\kappa \in [2, 6[$ and $q=a-2$ and $s=2$, 
	while, in the latter, the assumptions are fulfilled with 
	$\kappa \in [2, 6[$ and $q = 0$ and $s = 1$. 
\end{example}

%

\section{Parametrized Balanced Viscosity Solutions}\label{sec:defparaBV}

As already mentioned in the introduction,
it is well known that, due to the non-convexity of the energy and the positive homogeneity of the dissipation, 
\eqref{eq:ris} does in general not admit a ``classical'' solution, i.e., a solution $z \in W^{1,1}(0,T;\VV) \cap L^1(0,T;\ZZ)$ that fulfills 
\eqref{eq:ris} f.a.a.\ $t\in (0,T)$, cf., e.g., the counterexample in \cite[Section~2.3]{Tho22}. 
This solution concept is known as \emph{differential solution} and we will come back to this 
notion of solutions in Section~\ref{sec:revapp}.
Because of the lack of existence of a differential solution, 
several alternative solution concepts have been developed, among them \emph{global energetic solutions} and \emph{local solutions}. 
For a comprehensive overview we refer to \cite{MR15}. Since we aim to approximate the optimal control problem \eqref{eq:ocp0} 
by viscous approximation, we focus on the solution concept that arises through exactly this approximation procedure, 
which is known as \emph{parametrized balanced viscosity solution}.

\begin{definition}[Parametrized balanced viscosity (BV) solution]\label{def:paramsol}\index{parametrized balanced viscosity (BV) solution}
    We introduce the set of $\RR$-Lipschitz continuous functions by 
    \begin{equation*}
    \begin{aligned}
         \mathrm{AC}^\infty([0,T];\RR) 
         := \Big\{ z: [0,T] \to \ZZ\; | \; \exists\, L \geq 0, \; & \forall\, 0 \leq s < t \leq T: \\[-1.5ex]
         & \RR(z(t) - z(s)) \leq L(t-s)  \Big\}.    
    \end{aligned}
    \end{equation*}
    On $\mathrm{AC}^\infty([0,T];\RR)$ the limit 
    \begin{equation*}
        \RR[z'](t) := \lim_{h \searrow 0} \RR\Big( \frac{z(t+h) - z(t)}{h} \Big)
    \end{equation*}
    exists f.a.a.\ $t\in (0,T)$ and is called generalized metric derivative, cf. \ \cite[Prop.~2.2]{RMS08}.
    Then, for given data $z_0\in\ZZ$ and $\ell\in H^1(0,T;\VV^*)$, 
	we call a triple $(S,\hat t,\hat z)\in [T,\infty)\times W^{1,\infty}(0,T)\times \mathrm{AC}^\infty([0,S];\RR)\cap L^\infty (0,S;\ZZ)$ a 
	\emph{normalized, $\mathfrak{p}$-parametrized balanced viscosity solution} (in short: paramterized BV solution) 
	of the rate-independent system \eqref{eq:ris}, if the set 
	\begin{equation}\label{eq:defG.}
		G:= \big\{s\in[0,S]: \dist_{\VV^*}(-\mathrm{D}_z\II(\hat \ell(s),\hat z(s)),\partial \RR(0))>0\big\}
	\end{equation} 
	is a relatively open subset of $[0,S]$ and 
	\begin{gather*}
		\hat z \in W^{1,1}_{\textup{loc}}(G;\VV), \quad
		\mathrm{D}_z\II(\hat\ell(\cdot),\hat z(\cdot)) \in L^\infty_{\textup{loc}}(G;\VV^*),\\
		\|\hat z'(\cdot)\|_\VV \dist_{\VV^*}(-\mathrm{D}_z \II(\hat\ell(\cdot),\hat z(\cdot)),\partial \RR(0)) \in L^\infty(G;\R),
	\end{gather*}
	where $\hat \ell:=\ell\circ\hat t:[0,S]\to \VV^*$.
	Moreover, the following conditions shall be satisfied:\\
	\emph{Initial and end time condition:}
	\begin{align}
		\hat z(0)=z_0,\quad \hat t(0)=0,\quad \hat t (S)=T, \label{eq:initial}
	\end{align}
	\emph{Complementarity condition:}
	\begin{alignat}{3}
	\hat t'(s)&\geq 0 & \quad & \text{f.a.a. }s\in(0,S), \label{eq:compl2}\\
	\hat t'(s)\dist_{\VV^*}(-\mathrm{D}_z \II(\hat\ell(s),\hat z(s)),\partial \RR(0))&=0 &\quad& \text{f.a.a. }s\in(0,S)\label{eq:compl}
	\end{alignat}
	\emph{Normalization condition:} For almost all $s\in(0,S)$ it holds
	\begin{align}
	1=\begin{cases} \hat t'(s)+ \RR[\hat z'](s)+\|\hat z'(s)\|_\VV \dist_{\VV^*}(-\mathrm{D}_z \II(\hat\ell(s),\hat z(s)),\partial \RR(0)), &\text{ if } s\in G\\
					 \hat t'(s)+ \RR[\hat z'](s), &\text{ if } s\in [0,S]\setminus G \end{cases}\label{eq:normalization}
	\end{align}
	\emph{Energy identity:} For all $s\in [0,S]$ it holds
	\begin{equation}
	\begin{aligned}\index{energy identity}
		\II(\hat\ell(s),\hat z(s))
		& +\int_{0}^s\RR[\hat z'](r)\,\d r \\
		& +\int_{(0,s)\cap G}\|\hat z'(r)\|_\VV\dist_{\VV^*}(-\mathrm{D}_z \II(\hat\ell(r),\hat z(r)),\partial \RR(0)) \,\d r\\
	    & =\II(\hat\ell(0),z_0)-\int_0^s \langle \hat\ell'(r),\hat z(r)\rangle_{\VV^*,\VV}\,\d r. \label{eq:ener}
	\end{aligned}
	\end{equation}
\end{definition}

As its name indicates the existence of a paramterized BV solution is shown by means of viscous regularization. 
Let us formally sketch the procedure. First one approximates \eqref{eq:ris} by 
\begin{equation}\tag{$\mathrm{RIS}_\epsilon$}
	0 \in \partial \RR_{\epsilon}( z_\epsilon'(t)) + \mathrm{D}_z \II(\ell(t), z_\epsilon(t)) \quad \text{f.a.a.\ } t \in (0,T), 
    \quad z_\epsilon(0) = z_0, \label{eq:RISeps}
\end{equation}
where $\RR_\epsilon(z)=\RR(z)+\frac{\epsilon}{2}\|z\|_\VV^2$.
The existence and uniqueness of a solution to \eqref{eq:RISeps} will be shown below as a byproduct 
of the double viscous regularization in Section~\ref{sec:visreg}.
Secondly, the solution $z_\epsilon$ of \eqref{eq:RISeps} is reparametrized in time by means of the so-called 
\emph{vanishing viscosity contact potential} given by 
\begin{equation}\label{eq:frakp}
    \mathfrak{p}: \ZZ \times \VV^* \to \R,\quad 
    \mathfrak{p}(v, \xi) := \RR(v) + \|v\|_{\VV} \dist_{\VV^*}(w, \partial\RR(0)) .
\end{equation}
Given $\mathfrak{p}$, the parametrization reads 
\begin{equation}\label{eq:defs}
	s_{\epsilon}(t):= t+\int_{0}^{t}\mathfrak{p}( z_{\epsilon}'(r),-\mathrm{D}_{z}\II(\ell(r),z_{\epsilon}(r)))\,\d r,\quad S_{\epsilon}:=s_{\epsilon}(T).
\end{equation}
Since $s_\epsilon$ is strictly monotone increasing, its inverse exists and is denoted by 
\begin{equation}\label{eq:defhatt}
	\hat t_{\epsilon}:=(s_{\epsilon})^{-1}:[0,S_{\epsilon}]\to[0,T].
\end{equation}
Moreover, the parametrized solution is defined by 
\begin{equation}\label{eq:defhatz}
    \hat z_{\epsilon}:=z_{\epsilon} \circ \hat t_{\epsilon}:[0,S_{\epsilon}]\to\ZZ .
\end{equation}
A parametrized BV solution $(S, \hat t, \hat z)$ then arises as limit of $(S_\epsilon, \hat t_\epsilon, \hat z_\epsilon)$ as $\epsilon$ tends to zero. 
The associated limit analysis is rather involved, but has been carried out by various authors before; 
as an example, we refer to \cite{MRS16, Tho22}.
For the viscous approximation of the optimal control problems \eqref{eq:ocp0}, we need a slightly more general result, where 
the external loads are not fixed but vary with $\epsilon$ and converge weakly as $\epsilon$ tends to zero. 
Since, in our setting, the loads enter the system just linearly, the adaptation of the convergence analysis to (weakly) converging loads 
is straight forward. A comprehensive proof can be found in \cite[Section~3.3]{And25}. 

\begin{theorem}\label{thm:exparasol}
	Let $(z_\epsilon)_{\epsilon>0}$ be a sequence of solutions of \eqref{eq:RISeps} associated with $\ell_\epsilon\in H^1(0,T;\VV^*)$ 
	and a fixed initial value $z_0\in\ZZ$. 
	Moreover, let $S_\epsilon,\hat t_\epsilon,\hat z_\epsilon$ be defined as in \eqref{eq:defs}--\eqref{eq:defhatz} 
	and the loads are supposed to satisfy $\ell_\epsilon\rightharpoonup \ell$ in $H^1(0,T;\VV^*)$ as $\epsilon\searrow 0$. 
	Then there exists a subsequence $(\epsilon_n)_{n\in\N}$ and a limit 
	$(S,\hat t,\hat z)$ with $S\in [T,\infty)$, $\hat t \in W^{1,\infty}(0,S)$ and 
	$ \hat z\in \mathrm{AC}^\infty([0,S];\RR)\cap L^\infty(0,S;\ZZ)\cap C([0,S];\VV)$ such that 
    \begin{align}
        & S_{\epsilon_n}\to S, \hspace*{18.5ex}  \hat t_{\epsilon_n}\rightharpoonup^*\hat t \text{ in } W^{1,\infty}(0,S), 
        \label{eq:convparat}\\
        & \hat z_{\epsilon_n} \rightharpoonup^*\hat z \text{ in }  L^{\infty}(0,S;\ZZ),
        \quad  \hat z_{\epsilon_n}\to \hat z \text{ in } C([0,S];\VV),\label{eq:convparaz}\\
        & \hat z_{\epsilon_n}(s_n)\rightharpoonup \hat z (s) \text{ in } \ZZ \text{ for all converging sequences } s_n\to s.
        \label{eq:convparazs}
    \end{align} 
	Here the functions $\hat t_{\epsilon_n}$, $\hat z_{\epsilon_n}$ are constantly extended if $S_{\epsilon_n}< S$.
 	Moreover, the limit $(S,\hat t,\hat z)$ is a normalized, $\mathfrak{p}$-parametrized BV solution associated with $\ell$ 
 	in the sense of Definition \ref{def:paramsol}.
 \end{theorem}

\begin{remark}
    Let us mention that the statement of Theorem~\ref{thm:exparasol} can be slightly sharpened. 
    To be more precise, the limit vanishing viscosity analysis allows to prove additional smoothness results of the 
    limit $\hat z$, see \cite[Theorem~3.3.6]{And25}. Since these results are not needed for our analysis, 
    we do not go into detail here.
\end{remark}

In order to translate the convergent result from Theorem~\ref{thm:exparasol} back into physical time, 
we follow the lines of \cite[Section~3.8.2]{MR15} and define the set
\begin{equation*}
    \frakP(\hat t, \hat z) 
    := \{ z: [0,T] \to \VV \colon \,
    \forall\, t \in [0,T] \; \exists\,s\in [0,S] \colon (t, z(t)) = (\hat t(s), \hat z(s)) \}. 
\end{equation*}
Furthermore, we define the time points where $\hat t$ provides a progress in physical time as
\begin{equation*}
\begin{aligned}
    \frakC(\hat t) & := 
    \begin{aligned}[t]
        \bigcup \big\{ \hat t(s) \colon \, \exists\, s\in [0,S] \colon \; & \forall \, a \in [0,s) \colon\, \hat t(a) < \hat t(s),\\
        &  \forall \,b \in (s,S] \colon\, \hat t(s) < \hat t(b) \big\} 
    \end{aligned} \\
    &= \{ t \in [0,T] \colon \, \text{$\hat t^{-1}(t)$ is a singleton}\},
\end{aligned}
\end{equation*}
where the second equation follows from the monotony of $\hat t$.

\begin{lemma}\label{lem:contset}
    The set $[0,T]\setminus \frakC(\hat t)$ has Lebesgue measure zero.
\end{lemma}

\begin{proof}
    According to the area formula, which is applicable, since $\hat t$ is Lipschitz continuous, there holds that
    \begin{equation*}
        \int_\R \HH^0([0,S] \cap \hat t^{-1}(t))\,\d t = \int_0^S |\hat t'(r)|\,\d r 
        \leq \|\hat t\|_{W^{1,\infty}(0,S)}\, S < \infty.
    \end{equation*}
    Hence, $\HH^0([0,S] \cap \hat t^{-1}(t) < \infty$ for almost all $t\in (0,T)$. Thus, for almost all $t\in (0,T)$, 
    the set $\hat t^{-1}(t)$ is finite and, due to the monotony of $\hat t$, this implies that it is a singleton and 
    therefore belongs to $\frakC(\hat t)$.
\end{proof}

If $t\notin \frakC(\hat t)$, then, by definition of $\frakC(\hat t)$ and monotony of $\hat t$, the interval $[a^t, b^t]$ defined by
\begin{equation*}
    a^t := \min\{s\in [0,S] \colon \hat t(s) = t\}, \quad 
    b^t := \max\{s\in [0,S] \colon \hat t(s) = t\}
\end{equation*}
is non-empty. Note that $a^t$ and $b^t$ exist by continuity of $\hat t$.

\begin{lemma}\label{lem:sprunglimits}
    For every $t\notin \frakC(\hat t)$, there holds that $\hat z(a^t) \neq \hat z(b^t)$.
\end{lemma}

\begin{proof}
    Assume the contrary, i.e., there exists $t\notin \frakC(\hat t)$ with $\hat z(a^t) = \hat z(b^t)$. 
    Due to $\hat t(s) = t$ for all $s\in [a^t, b^t]$, there holds that 
    $\hat \ell(s) = \ell(\hat t(s)) = \ell(t)$ for all $s\in [a^t, b^t]$.
    Thus the energy equality in \eqref{eq:ener} implies 
    \begin{equation*}
    \begin{aligned}
	     \int_{a^t}^{b^t}\RR[\hat z'](r)\,\d r 
		& +\int_{(a^t,b^t)\cap G}\|\hat z'(r)\|_\VV\dist_{\VV^*}(-\mathrm{D}_z \II(\hat\ell(r),\hat z(r)),\partial \RR(0)) \,\d r\\
	    & =\II(\ell(t),\hat z(a^t)) - \II(\ell(t),\hat z(b^t)) 
	    + \int_{a^t}^{b^t} \langle \hat\ell'(r),\hat z(r)\rangle_{\VV^*,\VV}\,\d r = 0 
    \end{aligned}
    \end{equation*}
    Owing to the non-negativity of the integrands, this gives 
    $\RR[\hat z'](s) = 0$ a.e.\ in $(a^t, b^t)\setminus G$ and 
    $\RR[\hat z'](s) + \|\hat z'(s)\|_\VV\dist_{\VV^*}(-\mathrm{D}_z \II(\hat\ell(s),\hat z(s)),\partial \RR(0)) = 0$
    a.e.\ in $(a^t,b^t)\cap G$. The normalization condition \eqref{eq:normalization} therefore implies 
    that $\hat t'(s) = 1$ a.e.\ in $(a^t, b^t)$ contradicting that $\hat t$ is constant on $[a^t, b^t]$.
\end{proof}

\begin{lemma}\label{lem:continuitypts}
    For every $z \in \frakP(\hat t, \hat z)$,
    the set $\frakC(\hat t)$ coincides with the points of continuity of $z$. To be more precise, 
    if $t\in \frakC(\hat t)$, then $z$ is continuous with values in $\VV$ there and, on the other hand, 
    if $z$ is continuous in $t\in (0,T)$ with values in $\VV$, then $t\in \frakC(\hat t)$.
    As a consequence, $z$ is continuous a.e.\ in $[0,T]$.
\end{lemma}

\begin{proof}
    Let $t\in \frakC(\hat t)$ be arbitrary and consider an arbitrary sequence $(t_n)_{n\in \N} \subset [0,T]$ converging to $t$. 
    By definition of $\frakP(\hat t, \hat z)$, there exists $s_n$ with $t_n = \hat t(s_n)$ and $z(t_n) = \hat z(s_n)$ for all $n\in \N$. 
    Due to $(s_n)_{n\in \N} \subset [0,S]$, there exists a converging subsequence. Consider an arbitrary of these denoted by 
    $(s_{n_k})_{k\in \N}$, i.e., $s_{n_k} \to s \in [0,S]$. Then the continuity of $\hat t$ implies 
    \begin{equation*}
        \hat t (s) = \lim_{k\to \infty} \hat t(s_{n_k}) = \lim_{k\to\infty} t_{n_k} = t.
    \end{equation*}
    Since $t\in \frakC(\hat t)$, the limit $s$ is uniquely characterized by $\hat t(s) = t$ and thus the whole sequence 
    $(s_n)_{n\in \N}$ converges to $s$. The continuity of $\hat z$, then gives
    \begin{equation*}
        z(t_n) = \hat z(s_n) \to \hat z(s)= z(t) \quad \text{in }\VV,
    \end{equation*}
    which implies desired continuity of $z$ on $\frakC(\hat t)$.
    
    On the other hand, let $t\in (0,T)$ be a point of continuity of $z$. Assume by contrary that $t\notin \frakC(\hat t)$
    such that there exist $0 \leq a^t < b^t \leq S$ with $\hat t (a^t) = t = \hat t(b^t)$.
    Consider sequences $t_n \nearrow t$ and $\tau_n \searrow t$. By definition of $\frakP(\hat t, \hat z)$, there exist 
    sequences $(s_n)_{n\in\N}$ and $(\sigma_n)_{n\in\N}$ such that $t_n = \hat t(s_n)$, $z(t_n) = \hat z(s_n)$
    and $\tau_n = \hat t(\sigma_n)$, $z(\tau_n) = \hat z(\sigma_n)$. 
    Again we can select converging subsequences $(s_{n_k})_{k\in \N}$ and $(\sigma_{n_k})_{k\in \N}$. 
    Since $t_{n_k} < t$, the monotony of $\hat t$ implies $s_{n_k} < a^t$. Suppose that 
    $s := \lim_{k\to \infty} s_{n_k} < a^t$. Then the continuity of $\hat t$ leads to a contradiction, since
    \begin{equation*}
        t = \lim_{k\to\infty }t_{n_k} = \lim_{k\to \infty} \hat t(s_{n_k}) = \hat t(s) < t,
    \end{equation*}
    where the last inequality follows from the definition of $a^t$ and $s< a^t$.
    Therefore, every convergent subsequence converges to $a^t$ and thus the whole sequence $(s_n)$ 
    converges to $a^t$. By exactly the same argument, it follows that  $\sigma_n \to b^t$.
    The continuity of $\hat z$ along with Lemma~\ref{lem:sprunglimits} then yield the desired contradiction:
    \begin{equation*}
        z(t) = \lim_{n\to\infty} z(t_n)
        = \lim_{n\to\infty} \hat z(s_n) 
        = \hat z(a^t)
        \neq \hat z(b^t) = \lim_{n\to\infty} \hat z(\sigma_n)
        = \lim_{n\to\infty} z(\tau_n) = z(t).
    \end{equation*}
    Note that at initial and end time, i.e., when $t=0$ or $t=T$, then $z$ may well be (left- resp.\ right-)continuous, 
    even if $0 \notin \frakC(\hat t)$ and $T \notin \frakC(\hat t)$, respectively, since one may well choose 
    $b^0$ and $a^T$ for the curve parameter associated with $t = 0$ and $t=T$, 
    when picking the concrete element from $\frakP(\hat t, \hat z)$.
    
    Since $[0,T]\setminus \frakC(\hat t)$ has Lebesgue measure zero, we see that $z$ is continuous a.e.\ in $(0,T)$ as claimed.
\end{proof}

Given the above results, we can now state a convergence result in physical time:

\begin{corollary}\label{cor:convinphystime}
    In the situation of Theorem~\ref{thm:exparasol}, there holds for every $z \in \frakP(\hat t, \hat z)$ that 
    \begin{equation*}
        z_{\epsilon_n}(t) \weak z(t) \quad \text{in } \ZZ \quad \forall \, t \in \frakC(\hat t). 
    \end{equation*}
\end{corollary}

\begin{proof}
    Let $t\in \frakC(\hat t)$ be arbitrary. By definition of $\frakC(\hat t)$, there exists a unique $s\in (0,S)$ such that 
   $t = \hat t(s)$ and, by definition of $\frakP(\hat t, \hat z)$, we thus have $z(t) = \hat z(s)$ for all $z\in \frakP(\hat t, \hat z)$.
   Furthermore, due to $t\in [0,T]$, the construction of $\hat t_\epsilon$ implies that 
   for every $n\in \N$, there exists $s_n \in [0,S_{\epsilon_n}]$ such that 
    \begin{equation}\label{eq:sneq}
        \hat t_{\epsilon_n}(s_n) = t = \hat t(s).
    \end{equation}       
    Due to $(s_n)_{n\in \N} \subset [0,S]$, there exists a converging subsequence. 
    Take an arbitrary of these denoted by $(s_{n_k})_{k\in \N}$ with limit $\tilde s \in [0,S]$. 
    Then the uniform convergence of $\hat t_{\epsilon_n}$ to $\hat t$ by \eqref{eq:convparat}
    and \eqref{eq:sneq} imply that
    \begin{equation*}
       \hat t(\tilde s) = \lim_{k\to\infty} \hat t_{\epsilon_{n_k}}(s_{n_k}) = \hat t(s).
    \end{equation*}    
    Since $s$ is uniquely defined by $\hat t(s) = t$ as seen above, we obtain $\tilde s = s$ and, 
    because the converging subsequence was arbitrary, this implies the convergence of the whole 
    sequence $(s_n)_{n\in \N}$ to $s$. The pointwise weak convergence from \eqref{eq:convparazs} 
    thus yields
    \begin{equation*}
        z_{\epsilon_n}(t) =
        \hat z_{\epsilon_n}(s_n) 
        \weak \hat z(s) = z(t),
    \end{equation*}
    which establishes the claim.
\end{proof}

\begin{remark}
    The above result can be sharpened by introducing the notion of \emph{balanced-viscosity solution}, 
    which live in the physical time and are not parametrized. 
    We refer to \cite[Theorem~3.8.12]{MR15} for our semilinear setting and to \cite{MRS16} for a more general setting.
    Since the above result suffices for our purpose, we do not go into more details.
\end{remark}

\section{Double Viscous Regularization}\label{sec:visreg}

The existence of solutions to \eqref{eq:RISeps} is typically shown by means of time discretization, cf., e.g., \cite{MRS13}. 
Alternatively, one can apply a second viscous regularization involving a $\ZZ$-(semi-)norm instead of just the $\VV$-norm, 
see, e.g., \cite{KRZ11}.
From the view point of optimal control, the second approach has two important advantages. First, it leads to an
optimal control problem that can be treated with standard techniques from (nonsmooth) optimization, 
see Section~\ref{sec:optcontrol} below.
Secondly, and maybe more importantly, it allows to establish an estimate on the second derivative of the reduced energy $\EE$, 
see \eqref{eq:estE} in Lemma~\ref{lem:z''}. This estimate is of major importance for the reverse approximation property 
in Section~\ref{sec:revapp}, cf.\ the a priori estimate in \eqref{eq:H1Zbounddel} and the proof of Lemma~\ref{lem:strongconvH1}.
For these reasons, we study the additional viscous regularization involving the $\ZZ$-(semi-)norm in detail in the following.
This double viscous regularization reads as
\begin{align}\tag{$\mathrm{RIS}_{\epsilon,\delta}$}
	0 \in \partial \RR_{\epsilon,\delta}( z'(t)) + \Drm_z \II(\ell(t), z(t)) \quad \text{f.a.a.\ } t \in (0,T), 
    \quad z(0) = z_0, \label{eq:RISepsdel}
\end{align}
where 
\begin{equation*}
    \RR_{\epsilon,\delta}(z):= \RR(z)+\frac{\epsilon}{2} \, \|z\|_\VV^2+\frac{\delta}{2} \, |z|_\ZZ^2,
    \quad \epsilon,\delta>0.
\end{equation*}
In all what follows, we call \eqref{eq:RISeps} $\epsilon$-viscous regularization, whereas \eqref{eq:RISepsdel} is called 
$(\delta, \epsilon)$-viscous regularization.

\begin{lemma}\label{lem:exRIS}
	Let $\epsilon,\delta>0$, $z_0\in\ZZ$, and $\ell\in H^1(0,T;\VV^*)$ be given. 
	Then there exists a unique solution $z_{\epsilon,\delta} \in H^{2}(0,T;\ZZ)$ of \eqref{eq:RISepsdel}. 
	Moreover, $z_{\epsilon,\delta}$ satisfies the a priori estimate
	\begin{equation}\label{eq:boundz}
    \begin{aligned}
		\sup_{t\in[0,T]}\|z_{\epsilon,\delta}(t)\|_\ZZ + \frac{\epsilon}{2} \int_0^T \|z_{\epsilon,\delta}'(r)\|_\VV^2\, \d r
        	+ \frac{\delta}{2} \int_0^T |z_{\epsilon,\delta}'(r)|_\ZZ^2\, \d r \qquad & \\
	    \leq C(1+\|\ell\|_{H^1(0,T;\VV^*)}^2), &
    \end{aligned}
	\end{equation}
	where $C>0$ is independent of $\epsilon$, $\delta$ and $\ell$.
\end{lemma}

\begin{proof}
We start with the a priori estimate. To that end, assume that $z\in H^1(0,T;\ZZ)$ is a solution of \eqref{eq:RISepsdel}.
Then we apply the Fenchel-Young equality and the chain rule for Sobolev functions
to \eqref{eq:RISepsdel} to obtain
\begin{equation*}
	\RR_{\epsilon,\delta}(z'(t))+\RR_{\epsilon,\delta}^*(-\Drm_z\II(\ell(t),z(t)))
	= -\frac{\d}{\d t}\II(\ell(t),z(t))-\langle\ell'(t),z(t)\rangle_{\VV^*,\VV}
\end{equation*} 
for almost all $t\in (0,T)$. Integration over $(0,t)$ leads to
\begin{equation}
\begin{aligned}\label{eq:energyeq}
	\II(\ell(t),z(t))+\int_0^t\RR_{\epsilon,\delta}(z'(r))+\RR_{\epsilon,\delta}^*(-\Drm_z\II(\ell(r),z(r))\,\d r \\
	=\II(\ell(0),z(0))-\int_0^t \langle\ell'(r),z(r)\rangle_{\VV^*,\VV}\, \d r.
\end{aligned}
\end{equation}
With Young's inequality and \eqref{eq:assenergy1} at hand, the integral on the right hand side can be estimated by
\begin{equation}
\begin{aligned}\label{eq:estintright}
	\int_0^t \big|\langle\ell'(r),z(r)\rangle_{\VV^*,\VV}\big| \,\d r
	& \leq \frac{1}{2}\bigg(\|\ell\|_{H^1(0,t;\VV^*)}^2+\int_0^t\|z(r)\|_\ZZ^2\,\d r\bigg)\\
	& \leq \frac{1}{2}\bigg( (1+\mu)\|\ell\|_{H^1(0,t;\VV^*)}^2+\int_0^t\lambda \II(\ell(r),z(r))\,\d r\bigg).
\end{aligned}
\end{equation}
By exploiting the non-negativity of $\RR_{\epsilon,\delta}$ and $\RR_{\epsilon,\delta}^*$, equation \eqref{eq:energyeq} gives
\begin{equation*}
	\II(\ell(t),z(t))
	\leq \II(\ell(0),z_0)+\frac{1}{2}\bigg( (1+\mu)\|\ell\|_{H^1(0,t;\VV^*)}^2+\int_0^t\lambda \II(\ell(r),z(r))\,\d r\bigg)
\end{equation*}
such that applying Gronwall's inequality results in
\begin{align}\label{eq:estener}
	\II(\ell(t),z(t))\leq e^{\frac{1}{2}\lambda T}\bigg(\II(\ell(0),z_0)+\frac{1+\mu}{2}\|\ell\|_{H^1(0,T;\VV^*)}^2\bigg).
\end{align}
Eventually, we apply estimate \eqref{eq:estenergy} and end up with
\begin{align}\label{eq:boundpiz}
	\|z(t)\|_\ZZ\leq e^{\frac{1}{2}\lambda T}\bigg(\II(\ell(0),z_0)+\frac{1+\mu}{2}\|\ell\|_{H^1(0,T;\VV^*)}^2\bigg)
	+ \nu(\|\ell\|_{H^1(0,T;\VV^*)}^2 +1).
\end{align}
After inserting \eqref{eq:estener} in \eqref{eq:estintright} and exploiting 
\eqref{eq:estenergy}, the energy equality \eqref{eq:energyeq} results in
\begin{equation}\label{eq:Rest}
\begin{aligned}
    \int_0^T \RR_{\epsilon,\delta} (z'(r)) & + \RR_{\epsilon,\delta}^*(-\mathrm{D}_z\II(\ell(r),z(r)))\,\d r \\
    & \leq \nu(\|\ell\|_{H^1(0,T;\VV^*)}^2 +1)+C\Big(\II(\ell(0),z_0)+\|\ell\|_{H^1(0,T,\VV^*)}^2\Big),
\end{aligned}
\end{equation}
which in view of the non-negativity of $\RR$ and $\RR_{\epsilon,\delta}^*$ in turn leads to
\begin{multline*}
	\frac{\epsilon}{2} \int_0^T \|z'(r)\|_\VV^2\, \d r
	+ \frac{\delta}{2} \int_0^T |z'(r)|_\ZZ^2\, \d r \\
	 \leq  \nu(\|\ell\|_{H^1(0,T;\VV^*)}^2 +1) +C\Big(\II(\ell(0),z_0)+\|\ell\|_{H^1(0,T,\VV^*)}^2\Big).
\end{multline*}
Together with \eqref{eq:boundpiz}, this implies the desired a priori estimate.

To show existence of a solution, first note that, by employing Fenchel duality, \eqref{eq:RISepsdel} can equivalently be written as
\begin{equation}\tag{ODE$_{\epsilon,\delta}$}\label{eq:ODE}
	 z'(t) =\partial\RR_{\epsilon,\delta}^*\big(-\Drm_z\II(\ell(t),z(t))\big)\quad \text{f.a.a.\ } t \in (0,T),  \quad z(0) = z_0,
\end{equation} 
Note that, since $\partial\RR_{\epsilon,\delta}^*$ is single valued due to the uniform convexity of $\RR_{\epsilon,\delta}$ in $\ZZ$, 
this is an ODE with values in $\ZZ$.
The uniform convexity moreover implies that $\partial\RR_{\epsilon,\delta}^*: \ZZ^* \to \ZZ$ is globally Lipschitz continuous. 
Thus, the a priori estimates allow us to apply a standard truncation argument that guarantees 
the global Lipschitz continuity of the whole right hand side in \eqref{eq:ODE}, 
see \cite[Theorem~3.1.1]{And25} or \cite{KMS22} for details.
The Picard-Lindelöf theorem then yields the existence of a unique solution $z\in W^{1,\infty}(0,T;\ZZ)$ of \eqref{eq:ODE}.

Finally, we verify the improved regularity of  the solution of \eqref{eq:ODE}, i.e., $z\in H^{2}(0,T;\ZZ)$. 
Due to our assumptions on $\II$, the fact that $z\in W^{1,\infty}(0,T;\ZZ)$, and 
$\ell\in H^1(0,T;\VV^*)$, the inner function of the right hand side in \eqref{eq:ODE} is absolutely continuous. 
Therefore, the composition with the Lipschitz continuous function $\partial\RR_{\epsilon,\delta}$
is absolutely continuous, too, and thus $z'\in W^{1,1}(0,T;\ZZ)$. Hence, the second derivative exists 
a.e.\ and we obtain f.a.a.\ $t\in(0,T)$ that
\begin{align*}
    \|z''(t)\|_\ZZ&=\lim_{h\searrow 0}\frac{\|z'(t+h)-z'(t)\|_\ZZ}{h}\\
    &=\lim_{h\searrow 0}\frac{\big\|\partial\RR_{\epsilon,\delta}^*\big(-\Drm_z\II(\ell(t+h),z(t+h))\big)
    -\partial\RR_{\epsilon,\delta}^*\big(-\Drm_z\II(\ell(t),z(t))\big)\big\|_\ZZ}{h}\\
    &\leq L\Big\|-\frac{\d}{\d t}\Drm_z\II(\ell(t),z(t))\Big\|_{\ZZ^*},
\end{align*}
where $L>0$ denotes the Lipschitz constant of $\partial\RR_{\epsilon,\delta}$.
Consequently, $z''$ is bounded by an $L^2$-integrable function a.e., which gives $z\in H^{2}(0,T;\ZZ)$ as claimed.
\end{proof}

\begin{remark}\label{rem:Rest}
    In view of the structure of $\RR_{\epsilon, \delta}$ and the non-negativity of $\RR_{\epsilon, \delta}^*$, 
    the estimate in \eqref{eq:Rest} implies 
    \begin{equation*}
        \int_0^T \RR(z_{\epsilon, \delta}')(t)\,\d t
        \leq C(1+\|\ell\|_{H^1(0,T;\VV^*)}^2)
	\end{equation*}
	with a $C>0$  independent of $\epsilon$, $\delta$ and $\ell$. 
	This inequality will become important in the proof of Theorem~\ref{thm:W11bound} at the end of this section.
\end{remark}

\begin{lemma}\label{lem:z'(0)}
	Under our assumptions on the initial state, the solution $z_{\epsilon,\delta}$ of \eqref{eq:RISepsdel} fulfills $z_{\epsilon,\delta}'(0)=0$.
\end{lemma}

\begin{proof}
	Assumption \eqref{eq:init} implies $\RR(v)\geq \langle -\mathrm{D}_z\II(\ell(0),z_0),v\rangle_{\ZZ^*,\ZZ}$ 
	such that the definition of the conjugate functional gives
	\begin{align*}
		\RR_{\epsilon,\delta}^*(-\mathrm{D}_z\II(\ell(0),z_0))
		&=\sup_{v\in \ZZ}\Big(\langle -\mathrm{D}_z\II(\ell(0),z_0),v\rangle_{\ZZ^*,\ZZ}-\RR_{\epsilon,\delta}(v)\Big)\\
         &\leq\sup_{v\in \ZZ}\Big(\langle -\mathrm{D}_z\II(\ell(0),z_0),v\rangle_{\ZZ^*,\ZZ}-\RR(v)\Big)\leq 0.
	\end{align*}
    Since $\RR_{\epsilon,\delta}^*$ is non-negative, we observe that $\RR_{\epsilon,\delta}^*$ is minimal at $-\mathrm{D}_z\II(\ell(0),z_0)$ 
    and therefore, $0\in \partial \RR_{\epsilon,\delta}^*(-\mathrm{D}_z\II(\ell(0),z_0))$.
    Now, since $z_{\epsilon,\delta}'$ is continuous by the embedding $H^1(0,T;\ZZ)\hookrightarrow C([0,T];\ZZ)$,
	 \eqref{eq:ODE} holds for all $t\in [0,T]$. 
	 In particular, $z_{\epsilon,\delta}'(0)=\partial \RR_{\epsilon,\delta}^*(-\mathrm{D}_z\II(\ell(0),z_0))=0$ as claimed.
\end{proof}

Now, we have everything at hand to prove the essential estimate on the reduced energy $\EE$ mentioned above:

\begin{lemma}\label{lem:z''}
	The solution $z_{\epsilon,\delta}$ of \eqref{eq:RISepsdel} satisfies
	\begin{equation}
	\begin{aligned}\label{eq:z''}
		\epsilon\dual{z_{\epsilon,\delta}''(t)}{z_{\epsilon,\delta}'(t)}_{\VV^*,\VV}&+\delta\dual{Az_{\epsilon,\delta}''(t)}{z_{\epsilon,\delta}'(t)}_{\ZZ^*,\ZZ}\\
		&\qquad+\Big\langle\frac{\d}{\d t}\mathrm{D}_z\II(\ell(t),z_{\epsilon,\delta}(t)),z_{\epsilon,\delta}'(t)\Big\rangle_{\ZZ^*,\ZZ}=0.
	\end{aligned} 
	\end{equation}
	for almost all $t\in(0,T)$.
	Furthermore, we have for all $t\in[0,T]$ that
	\begin{equation}\label{eq:estE}
	\begin{aligned}
		\frac{\epsilon}{2}\,\|z_{\epsilon,\delta}'(t)\|_\VV^2 
		+ \frac{\delta}{2}\, |z_{\epsilon, \delta}'(t)|_{\ZZ}^2
		+ \int_0^t\mathrm{D}_z^2\EE(z_{\epsilon,\delta}(r))[z_{\epsilon,\delta}'(r),z_{\epsilon,\delta}'(r)]\,\d r \qquad\qquad & \\
		= \int_0^t\langle \ell'(r),z_{\epsilon,\delta}'(r)\rangle_{\VV^*,\VV}\,\d r . &
	\end{aligned}
	\end{equation}
\end{lemma}

\begin{proof}
	To prove \eqref{eq:z''} define 
	\begin{equation*}
        h(t):= \epsilon z_{\epsilon,\delta}'(t)+\delta A z_{\epsilon,\delta}'(t)+\mathrm{D}_z\II(\ell(t),z_{\epsilon,\delta}(t))\in H^{1}(0,T;\ZZ^*).
	\end{equation*}
	Then, the definition of $\RR_{\epsilon,\delta}$ along with \eqref{eq:RISepsdel} 
    implies that $0 \in h(t) + \partial \RR(z_{\epsilon,\delta}'(t))$ f.a.a.\ $t\in[0,T]$ 
    and, again, the continuity of $z_{\epsilon, \delta}'$ yields that this equality even holds for all $t\in [0,T]$.
	This in turn yields
	\begin{equation*}
		\langle -h(t),z_{\epsilon,\delta}'(t)\rangle_{\ZZ^*,\ZZ} =\RR(z_{\epsilon,\delta}'(t)), \quad 
		\langle -h(s), z_{\epsilon,\delta}'(t)\rangle_{\ZZ^*,\ZZ} \leq \RR(z_{\epsilon,\delta}'(t))
	\end{equation*}
	for all $t,s\in[0,T]$. Thus we obtain for every $t\in (0,T)$ and every $\tau > 0$ small enough (so that $s := t \pm \tau \in [0,T]$) that
	$\frac{1}{\tau}\langle h(t\pm \tau)-h(t),z_{\epsilon,\delta}'(t)\rangle_{\ZZ^*,\ZZ}\geq 0$.
    Now, due to $h\in H^{1}(0,T;\ZZ^*)$, almost every $t\in(0,T)$ is a point of differentiability of $h$ and therefore, 
    we can pass to the limit $\tau\searrow 0$ a.e.\ in $(0,T)$. 
    This gives $\langle h'(t),z_{\epsilon,\delta}'(t)\rangle_{\ZZ^*,\ZZ}= 0$ f.a.a.\ $t\in (0,T)$, 
    which, in view of the definition of $h$, is nothing else than \eqref{eq:z''}.

    To verify \eqref{eq:estE}, we integrate \eqref{eq:z''} over $(0,t)$ to obtain 
    \begin{equation}\label{eq:int46}
	\begin{aligned}
		0 &= \int_0^t
        \begin{aligned}[t]
            & \frac{\epsilon}{2}\frac{\d}{\d t}\|z_{\epsilon,\delta}'(r)\|_\VV^2
		    + \frac{\delta}{2}\frac{\d}{\d t}|z_{\epsilon,\delta}'(r)|_\ZZ^2 \\
        		& + \langle \frac{\d}{\d t}\mathrm{D}_z\II(\ell(r),z_{\epsilon,\delta}(r)),z_{\epsilon,\delta}'(r)\rangle_{\ZZ^*,\ZZ}\,\d r
        \end{aligned} \\
		&= \frac{\epsilon}{2}\|z_{\epsilon,\delta}'(t)\|_\VV^2-\frac{\epsilon}{2}\|z_{\epsilon,\delta}'(0)\|_\VV^2
		+\frac{\delta}{2}|z_{\epsilon,\delta}'(t)|_\ZZ^2 -\frac{\delta}{2}|z_{\epsilon,\delta}'(0)|_\ZZ^2\\
		&\qquad + \int_0^t \mathrm{D}_z^2\EE(z_{\epsilon,\delta}(r))[z_{\epsilon,\delta}'(r),z_{\epsilon,\delta}'(r)]
		-\langle\ell'(r),z_{\epsilon,\delta}'(r)\rangle_{\VV^*,\VV}\,\d r.
	\end{aligned}	
    \end{equation}
	Eventually, thanks to Lemma~\ref{lem:z'(0)}, we end up with \eqref{eq:estE}.
\end{proof}

With this lemma at hand, we can improve the a priori estimates from Lemma~\ref{lem:exRIS}.

\begin{lemma}\label{lem:zdelH1Z}
	Let $\delta>0$, $\epsilon \in \,]0,1]$, and $\ell\in H^1(0,T;\VV^*)$ be given. 
    Then the solution $z_{\epsilon,\delta}$ of \eqref{eq:RISepsdel} satisfes
    \begin{equation}\label{eq:zdelH1Z}
            \epsilon \, \|z_{\epsilon,\delta}'\|_{L^2(0,T;\ZZ)}^2 
        + \epsilon^2 \, \|z_{\epsilon,\delta}'\|_{L^\infty(0,T;\VV)}^2  \leq C (1 + \|\ell\|_{H^1(0,T;\VV^*)}^2)^{2q+1}
    \end{equation}        
    where $q$ is the exponent from assumption~\eqref{eq:F2} and $C>0$ is independent of $\epsilon$, $\delta$, and $\ell$. 
\end{lemma}

\begin{proof}
Let us abbreviate $\|\ell\| := \|\ell\|_{H^1(0,T;\VV^*)}$.
By exploiting the coercivity of $A$, we infer from \eqref{eq:estE} that
\begin{align*}
	\alpha\int_0^T\|z_{\epsilon,\delta}'(r)\|_\ZZ^2\,\d r&\leq \int_0^T \langle Az_{\epsilon,\delta}'(r),z_{\epsilon,\delta}'(r)\rangle_{\ZZ^*,\ZZ}\,\d r\\
	&\leq -\int_0^T \mathrm{D}_z^2\FF(z_{\epsilon,\delta}(r))[z_{\epsilon,\delta}'(r),z_{\epsilon,\delta}'(r)]\,\d r+\int_0^T \langle \ell'(r),z_{\epsilon,\delta}'(r)\rangle_{\VV^*,\VV}\,\d r
\end{align*}
Next, along with \eqref{eq:F2} and the $L^\infty(0,T;\ZZ)$-estimate from \eqref{eq:boundz}, 
Young's inequality with arbitrary $\rho>0$ yields for the first integral on the right hand side that
\begin{align*}
	-\int_0^T \mathrm{D}_z^2\FF(z_{\epsilon,\delta}(r))[z_{\epsilon,\delta}'(r),z_{\epsilon,\delta}'(r)]\,\d r
	& \leq \int_0^T \gamma(1+\|z_{\epsilon,\delta}(r)\|_\ZZ^q)\|z_{\epsilon,\delta}'(r)\|_\ZZ \|z_{\epsilon,\delta}'(r)\|_\VV \,\d r \\
	& \leq \int_0^T C(1+ \|\ell\|^{2})^q \Big(\rho\|z_{\epsilon,\delta}'(r)\|_\ZZ^2+\frac{1}{4\rho}\|z_{\epsilon,\delta}'(r)\|_\VV^2\Big)\,\d r .
\end{align*}
Thus, if we choose $\rho = \frac{\alpha}{2 C(1+ \|\ell\|^{2})^q}$, we obtain
\begin{equation*}
    \frac{\alpha}{2}\int_0^T\|z_{\epsilon,\delta}'(r)\|_\ZZ^2\,\d r
    \leq \Big(\frac{C^2(1+ \|\ell\|^{2})^{2q}}{2\alpha} + \frac{1}{2}\Big) \|z_{\epsilon,\delta}'\|_{L^2(0,T;\VV)}^2 + \frac{1}{2}\,\|\ell\|^2.
\end{equation*}
Thus, the $H^1(0,T;\VV)$-bound from \eqref{eq:boundz} implies 
\begin{equation*}
    \epsilon\|z_{\epsilon,\delta}'\|_{L^2(0,T;\ZZ)}^2 \leq C(1 + \|\ell\|^2)^{2q+1} +  \epsilon\, \|\ell\|^2.
\end{equation*}
Since $\epsilon \leq 1$ by assumption, the desired $H^1(0,T;\ZZ)$-estimate follows.

To prove the $W^{1,\infty}(0,T;\VV)$-bound, let $t\in[0,T]$ be arbitrary. 
As in the proof of Lemma \ref{lem:z''}, we integrate \eqref{eq:z''} over $(0,t)$ to obtain \eqref{eq:int46}.
Then, from  the coercivity of $A$ and Lemma~\ref{lem:z'(0)}, we deduce
\begin{align*}
	&\frac{\epsilon}{2}\|z_{\epsilon,\delta}'(t)\|_\VV^2+\alpha \int_0^t \|z_{\epsilon,\delta}'(r)\|_\ZZ^2\,\d r\\
	&\quad \leq \int_0^t\Big|\mathrm{D}_z^2\FF(z_{\epsilon,\delta}(r))[z_{\epsilon,\delta}'(r),z_{\epsilon,\delta}'(r)]\Big|\,\d r
	+\int_0^t \langle \ell'(r),z_{\epsilon,\delta}'(r)\rangle_{\VV^*,\VV}\,\d r\\
	& \quad \leq 
	\begin{aligned}[t]
    	    \int_0^t C(1+ \|\ell\|^{2})^q \Big(\rho\|z_{\epsilon,\delta}'(r)\|_\ZZ^2
	    & +\frac{1}{4\rho}\|z_{\epsilon,\delta}'(r)\|_\VV^2\Big)\,\d r \\
	    & +\frac{1}{2}\, \|\ell\|^2
	    +\frac{1}{2}\, \|z_{\epsilon,\delta}'\|_{L^2(0,T;\VV)}^2
	\end{aligned}
\end{align*} 
where we again used \eqref{eq:F2}, the $L^\infty(0,T;\ZZ)$-bound from \eqref{eq:boundz} and 
Young's inequality with $\rho>0$ arbitrary in the last step. 
By choosing $\rho=\frac{\alpha}{C(1+ \|\ell\|^{2})^q}$, we again infer from the $H^1(0,T;\VV)$-estimate in \eqref{eq:boundz} that
\begin{align*}
	\epsilon^2 \|z_{\epsilon,\delta}'(t)\|_\VV^2
    & \leq \Big(\frac{C^2(1+ \|\ell\|^{2})^{2q}}{2\alpha} + 1\Big) \epsilon\,\|z_{\epsilon,\delta}'\|_{L^2(0,T;\VV)}^2 +\epsilon\, \|\ell\|^2\\
    & \leq C(1 + \|\ell\|^2)^{2q+1} +  \epsilon\, \|\ell\|^2,
\end{align*}
which gives the desired estimate for the $L^\infty(0,T;\VV)$-norm of $z'_{\epsilon, \delta}$.
\end{proof}

\begin{lemma}\label{lem:epsunique}
	Let $z_1,z_2\in H^1(0,T;\ZZ)$ be solutions of \eqref{eq:RISeps}. Then it holds $z_1=z_2$.
\end{lemma}

\begin{proof}
	Since $z_1,z_2$ are both solutions of \eqref{eq:RISeps}, they fulfill
	\begin{equation*}
	    \RR(v) \geq \RR(z_i'(t)) - \epsilon \dual{z_i'(t)}{v - z_i'(t)}_{\VV^*, \VV} 
	    - \dual{\Drm_z\II(\ell(t), z_i'(t))}{v - z_i'(t)}_{\ZZ^*, \ZZ}, \;\, i = 1,2    
	\end{equation*}
	for all $v\in \ZZ$ and almost all $t\in(0,T)$.
	Testing the variational inequality for $z_1$ with $z_2$ and vice versa and adding the arising inequalities gives
	\begin{align*}
		0\geq \epsilon \|z_1'(t)-z_2'(t)\|_\VV^2
		+ \langle \mathrm{D}_z \II(\ell(t), z_1(t))-\mathrm{D}_z \II(\ell(t), z_2(t)),z_1'(t)-z_2'(t)\rangle_{\ZZ^*,\ZZ}.
	\end{align*}
	By exploiting the concrete structure of the energy and the assumed symmetry of $A$ from \eqref{eq:Acoer}, we deduce
    \begin{equation}\label{eq:gronwallprep}
    \begin{aligned}
		& \frac{1}{2} \frac{\d}{\d t} \langle Az_1(t)-Az_2(t),z_1(t)-z_2(t)\rangle_{\ZZ^*,\ZZ}
		+\epsilon \|z_1'(t)-z_2'(t)\|_\VV^2\\
		&\qquad \leq -\langle \mathrm{D}_z\FF(z_1(t))-\mathrm{D}_z\FF(z_2(t)),z_1'(t)-z_2'(t)\rangle_{\VV^*,\VV}\\
		& \qquad \leq C\|z_1'(t)-z_2'(t)\|_\VV \|z_1(t)-z_2(t)\|_\ZZ\\
		& \qquad \leq \epsilon\|z_1'(t)-z_2'(t)\|_\VV^2+ C_\epsilon \|z_1(t)-z_2(t)\|_\ZZ^2    
    \end{aligned}        
    \end{equation}
	with a constant $C_\epsilon>0$ depending on $\epsilon$. 
	In the second last step we used the mean value theorem in combination 
	with the boundedness of the second derivative of $\FF$ by \eqref{eq:F2} and the fact that 
	$z_1,z_2\in H^1(0,T;\ZZ)\embed L^\infty(0,T;\ZZ)$.
	Eventually, integrating over $(0,s)$ and using the coercivity of $A$ leads to 
	\begin{align*}
		\frac{\alpha}{2}\|z_1(s)-z_2(s)\|_\ZZ^2\leq \int_0^s C_\epsilon \|z_1(t)-z_2(t)\|_\ZZ^2\, \d t,
	\end{align*}
	where we exploited the initial condition $z_1(0)=z_2(0)=z_0$.
	Finally, from Gronwall's inequality we infer that $z_1=z_2$, since $s$ was arbitrary.
\end{proof}

By means of the previous results,
we can now prove the existence of a solution of \eqref{eq:RISeps} 
as the limit of solutions of \eqref{eq:RISepsdel} for $\delta \searrow 0$.
As a byproduct, we obtain the existence of solutions to \eqref{eq:RISeps}, which however is a known result
as already mentioned above. 
What is more, instead of considering a fixed external load, 
we let the loads converge weakly in $H^1(0,T;\VV^*)$ as $\delta \searrow 0$, 
similarly to Theorem~\ref{thm:exparasol}. 
This is essential for the approximability of solutions to \eqref{eq:ocpeps} by solutions of \eqref{eq:ocpepsdelta}, 
see Section~\ref{sec:ocpapprox} below.

\begin{theorem}[Convergence of the double viscous regularization]\label{thm:exRISeps}
	Let $\epsilon>0$ be fixed and $(\ell_{\epsilon,\delta})_{\delta>0}\subset H^1(0,T;\VV^*)$ be a weakly converging sequence, 
	i.e., $\ell_{\epsilon,\delta}\rightharpoonup \ell_{\epsilon}$ as $\delta\searrow 0$ with some limit $\ell_\epsilon \in H^1(0,T;\VV^*)$. 
    Denote by $(z_{\epsilon,\delta})_{\delta>0}$ the sequence of solutions to the $(\epsilon, \delta)$-viscous regularization
    \eqref{eq:RISepsdel} associated with $\ell_{\epsilon,\delta}$. Then there exists a limit $z_\epsilon\in H^1(0,T;\ZZ)\cap W^{1,\infty}(0,T;\VV)$  
	such that 
	\begin{equation}
		z_{\epsilon,\delta}\rightharpoonup z_\epsilon \; \text{ in }  H^1(0,T;\ZZ), \quad 
		z_{\epsilon,\delta}\rightharpoonup^* z_\epsilon \; \text{ in }  W^{1,\infty} (0,T;\VV), \label{eq:convergencedelta}
    \end{equation}
    as $\delta \searrow 0$. Moreover, $z_\epsilon$ is the solution of \eqref{eq:RISeps} associated with $\ell_\epsilon$.
\end{theorem}

\begin{proof}
    The proof is based on standard energetic arguments. For the sake of completeness, we present it in detail.
    First note that, thanks to its weak convergence, the sequence $(\ell_{\epsilon, \delta})_{\delta> 0}$ is bounded in 
    $H^1(0,T;\VV^*)$. Therefore, due to \eqref{eq:boundz} and \eqref{eq:zdelH1Z}, 
	the sequence $(z_{\epsilon,\delta})_{\delta>0}$ is bounded in $H^1(0,T;\ZZ)\cap W^{1,\infty}(0,T;\VV)$ by a constant 
	not depending on $\delta$. Therefore, by passing $\delta$ to zero, we can extract a weakly converging subsequence such that
	\begin{equation}\label{eq:zweakH1}
		z_{\epsilon,\delta_n}\rightharpoonup z_\epsilon \text{ in }  H^1(0,T;\ZZ), 
		\quad z_{\epsilon,\delta_n}\rightharpoonup^* z_\epsilon \text{ in }  W^{1,\infty} (0,T;\VV).
	\end{equation}
	Since $H^1(0,T;\ZZ) \embed C([0,T];\ZZ)$, the point evaluation in time is linear and continuous thus weakly continuous.
	Moreover, due to the Aubin-Lions lemma, the embedding $H^1(0,T;\ZZ)\hookrightarrow C([0,T];\VV)$ is compact
	so that point evaluation in time with values in $\VV$ is even a compact operator. 
	Thus we obtain
	\begin{equation}\label{eq:pktweak2}
		z_{\epsilon,\delta_n}(t)\rightharpoonup z_\epsilon(t) \text{ in } \ZZ, 
		\quad z_{\epsilon,\delta_n}(t)\to z_\epsilon(t) \text{ in } \VV \text{ for all } t\in [0,T].
	\end{equation}
	Consequently, there holds that $z_0=z_{\epsilon,\delta_n}(0)\to z_\epsilon(0)$ in $\VV$ as $\delta_n \searrow 0$,
	such that the limit satisfies the initial condition $z_\epsilon(0)=z_0$.
	
	By definition of the subdifferential, it is readily seen that the solution $z_{\epsilon,\delta_n}$ of \eqref{eq:RISepsdel} fulfills 
    \begin{equation}\label{eq:energyineqepsdelta}
	\begin{aligned}
		\int_0^t\RR(v(r))\,\d r
		& \geq \int_0^t\RR(z_{\epsilon,\delta_n}'(r))\,\d r + \int_0^t\epsilon\langle z_{\epsilon,\delta_n}'(r),z_{\epsilon,\delta_n}'(r)-v(r)\rangle_{\VV^*,\VV}\,\d r\\
		& \quad + \int_0^t\delta_n\langle Az_{\epsilon,\delta_n}'(r),z_{\epsilon,\delta_n}'(r)-v(r)\rangle_{\ZZ^*,\ZZ}\,\d r\\
		& \quad + \int_0^t\langle \mathrm{D}_z\II(\ell_{\epsilon,\delta_n}(r),z_{\epsilon,\delta_n}(r)),z_{\epsilon,\delta_n}'(r)-v(r)\rangle_{\ZZ^*,\ZZ}\,\d r
	\end{aligned}
    \end{equation}
	for all $t\in[0,T]$ and all $v\in L^2(0,T;\ZZ)$. Now we discuss each term on the right hand side separately. 
	First, due to the weak convergence 
	$z_{\epsilon,\delta_n}'\rightharpoonup z_\epsilon'$ in $L^2(0,T;\ZZ)$ and the convexity and lower semicontinuity of $\RR$ from \eqref{R1}, 
	one deduces that
	\begin{equation}\label{eq:Repsconv}
		\liminf_{n\to\infty}\int_0^t\RR(z_{\epsilon,\delta_n}'(r))\,\d r\geq \int_0^t\RR(z_{\epsilon}'(r))\,\d r, 
	\end{equation}
	cf. \cite[Lemma~A.1]{And25} or \cite[Lemma~A.3.5]{Sie20}.
	Regarding the second term, we exploit the weak lower semicontinuity of the norm and \eqref{eq:zweakH1} in order to obtain
	\begin{equation}\label{eq:epsconv}
	\begin{aligned}
		 \liminf_{n\to\infty} \int_0^t \epsilon
		 \langle z_{\epsilon,\delta_n}'(r),z_{\epsilon,\delta_n}'(r)-v(r)\rangle_{\VV^*,\VV}\,\d r \qquad\qquad\qquad &\\
		\geq \int_0^t\epsilon\langle z_{\epsilon}'(r),z_{\epsilon}'(r)-v(r)\rangle_{\VV^*,\VV}\,\d r . &
	\end{aligned}
	\end{equation}
	Next, because of \eqref{eq:zdelH1Z} and the boundedness of $(\ell_{\epsilon,\delta_n})_{n\in\N}$, the third integral satisfies
	\begin{equation*}
		\lim_{n\to \infty} \int_0^t\delta_n\langle Az_{\epsilon,\delta_n}'(r),z_{\epsilon,\delta_n}'(r)-v(r)\rangle_{\ZZ^*,\ZZ}\,\d r = 0 . 
	\end{equation*}
	Concerning the fourth integral, we first consider the part of the energy involving the external load. 
	By applying integration by parts, this term can be written as
	\begin{align*}
		& \int_0^T\langle \ell_{\epsilon,\delta_n}(r),z_{\epsilon,\delta_n}'(r)-v(r)\rangle_{\VV^*,\VV}\,\d r \\
	    & \quad = \langle \ell_{\epsilon,\delta_n}(T),z_{\epsilon,\delta_n}(T)\rangle_{\VV^*,\VV}
	    - \langle \ell_{\epsilon,\delta_n}(0),z_{\epsilon,\delta_n}(0)\rangle_{\VV^*,\VV} \\
	    & \qquad -\int_0^T\langle \ell_{\epsilon,\delta_n}'(r),z_{\epsilon,\delta_n}(r)\rangle_{\VV^*,\VV}\,\d r
	    -\int_0^T\langle \ell_{\epsilon,\delta_n}(r),v(r)\rangle_{\VV^*,\VV}\,\d r.
	\end{align*}
	The weak continuity of the point evaluation in combination with the weak convergence of $(\ell_{\epsilon,\delta_n})_{n\in\N}$ in $H^1(0,T;\VV^*)$ 
	implies pointwise weak convergence of the loads in $\VV^*$ and thus, thanks to \eqref{eq:pktweak2}, we deduce
	\begin{multline*}
		\langle \ell_{\epsilon,\delta_n}(T),z_{\epsilon,\delta_n}(T)\rangle_{\VV^*,\VV}
		-\langle \ell_{\epsilon,\delta_n}(0),z_{\epsilon,\delta_n}(0)\rangle_{\VV^*,\VV}\\
		\to \langle \ell_\epsilon(T),z_\epsilon(T)\rangle_{\VV^*,\VV}-\langle \ell_\epsilon(0),z_\epsilon(0)\rangle_{\VV^*,\VV}.
	\end{multline*}
	Furthermore, exploiting the strong convergence $z_{\epsilon,\delta_n}\to z_\epsilon$ in $L^2(0,T;\VV)$ by compact embedding, leads to
	\begin{align*}
		-\int_0^T\langle \ell_{\epsilon,\delta_n}'(r)&,z_{\epsilon,\delta_n}(r)\rangle_{\VV^*,\VV}\,\d r -\int_0^T\langle \ell_{\epsilon,\delta_n}(r),v(r)\rangle_{\VV^*,\VV}\,\d r\\
		&\to -\int_0^T\langle {\ell_\epsilon}'(r),z_\epsilon(r)\rangle_{\VV^*,\VV}\,\d r -\int_0^T\langle \ell_\epsilon(r),v(r)\rangle_{\VV^*,\VV}\,\d r
	\end{align*}
	such that we obtain after applying integrations by parts again
	\begin{align*}
		\int_0^T\langle \ell_{\epsilon,\delta_n}(r),z_{\epsilon,\delta_n}'(r)-v(r)\rangle_{\VV^*,\VV}\,\d r\to \int_0^T\langle \ell_\epsilon(r),z_\epsilon'(r)-v(r)\rangle_{\VV^*,\VV}\,\d r,~ n\to\infty. 
	\end{align*}
	Together with Lemma~\ref{lem:limAF} in Appendix~\ref{sec:limAF},
	which is applicable due to  \eqref{eq:zweakH1}, we eventually infer
    \begin{equation}\label{eq:DzIepsconv}
	\begin{aligned}
		& \liminf_{n\to\infty}\int_0^t \langle \mathrm{D}_z\II(\ell_{\epsilon,\delta_n}(r),
		z_{\epsilon,\delta_n}(r)),z_{\epsilon,\delta_n}'(r)-v(r)\rangle_{\ZZ^*,\ZZ}\,\d r\\
		&\quad = \liminf_{n\to\infty}\int_0^t \langle Az_{\epsilon,\delta_n}(r),z_{\epsilon,\delta_n}'(r)-v(r)\rangle_{\ZZ^*,\ZZ}\,\d r\\
		&\qquad\qquad+\lim_{n\to\infty}\int_0^t\langle \mathrm{D}_z\FF(z_{\epsilon,\delta_n}(r)),z_{\epsilon,\delta_n}'(r)-v(r)\rangle_{\ZZ^*,\ZZ}\,\d r \\
		&\qquad\qquad\qquad-\lim_{n\to\infty}\int_0^t\langle \ell_{\epsilon,\delta_n}(r),z_{\epsilon,\delta_n}'(r)-v(r)\rangle_{\ZZ^*,\ZZ}\,\d r\\
		&\quad \geq \int_0^t \big( \langle Az_{\epsilon}(r),z_{\epsilon}'(r)-v(r)\rangle_{\ZZ^*,\ZZ} \\
		&\quad\quad  \qquad + \langle \mathrm{D}_z\FF(z_{\epsilon}(r)),z_{\epsilon}'(r)-v(r)\rangle_{\ZZ^*,\ZZ}
		-  \langle \ell_\epsilon(r),z_{\epsilon}'(r)-v(r)\rangle_{\ZZ^*,\ZZ}\big) \d r\\
		&\quad =\int_0^t \langle \mathrm{D}_z\II(\ell_\epsilon(r),z_{\epsilon}(r)),z_{\epsilon}'(r)-v(r)\rangle_{\ZZ^*,\ZZ}\,\d r.
	\end{aligned}    
    \end{equation}
	All in all, we have thus shown that the limit $z_\epsilon$ satisfies 
	\begin{equation}\label{eq:VIeps}
	\begin{aligned}
        	\int_0^t\RR(v(r))\,\d r \geq 
        \int_0^t\RR(z_{\epsilon}'(r)) &+ \epsilon\langle z_{\epsilon}'(r),z_{\epsilon}'(r)-v(r)\rangle_{\VV^*,\VV}\\
		&+ \langle \mathrm{D}_z\II(\ell_\epsilon(r),z_{\epsilon}(r)),z_{\epsilon}'(r)-v(r)\rangle_{\ZZ^*,\ZZ}\,\d r.
	\end{aligned}
	\end{equation}
    Standard arguments based on the fundamental lemma of the calculus of variations show that this   
    is an equivalent formulation of \eqref{eq:RISeps}. 
    Therefore, since $z_\epsilon$ also satisfies the initial condition, it is a solution of  \eqref{eq:RISeps} and, 
    since the latter one is unique by Lemma~\ref{lem:epsunique}
    the entire sequence converges as claimed in \eqref{eq:convergencedelta}.
\end{proof}

Owing to the convergence results in \eqref{eq:convergencedelta},
the ($\delta$-independent) a priori bounds from \eqref{eq:boundz} and \eqref{eq:zdelH1Z} 
readily transfer to the limit $z_\epsilon$ yielding the following results that will turn out to be useful 
in the upcoming analysis:

\begin{corollary}\label{cor:boundszeps}
	The solution $z_\epsilon\in H^1(0,T;\ZZ)$ of \eqref{eq:RISeps} fulfills the a priori estimate
	\begin{align}
	    \sup_{t\in[0,T]} \|z_\epsilon(t)\|_\ZZ &\leq C(1 + \|\ell\|_{H^1(0,T;\VV^*)}^2), \label{eq:boundzeps} \\
	    \epsilon \|z_{\epsilon}\|_{W^{1,\infty}(0,T;\VV)}
	 	+ \epsilon \|z_{\epsilon}\|_{H^1(0,T;\ZZ)}^2 &\leq C(1 + \|\ell\|_{H^1(0,T;\VV^*)}^2)^{2q+1}, 
	 	\label{eq:boundDzeps}
	\end{align}
	  where with a constant $C>0$ is independent of $\epsilon$.
\end{corollary}

\begin{corollary}\label{cor:DzI}
There exists a constant $C>0$ not depending on $\epsilon$ such that
\begin{equation}\label{eq:boundDzIeps}
	\|\mathrm{D}_z\II(\ell(t),z_{\epsilon}(t))\|_{\VV^*}\leq C(1 + \|\ell\|_{H^1(0,T;\VV^*)}^2)^{2q+1}
	\quad \forall\, t\in [0,T].
\end{equation}
\end{corollary}

\begin{proof}
    As $z_\epsilon$ is the solution of \eqref{eq:RISeps}, it holds for almost every $t\in (0,T)$ that
    \begin{equation*}
        - \mathrm{D}_z\II(\ell(t),z_\epsilon(t))
        \in\partial\RR(z_\epsilon '(t)) + \epsilon\, z_\epsilon'(t) 
        \subset \partial\RR(0) + \epsilon\, z_\epsilon'(t) 
    \end{equation*}
    Now the boundedness of $\partial \RR(0)\subset \VV^*$ by Lemma \ref{lem:dR(0)bound} and the 
    a priori estimate from \eqref{eq:boundDzeps} yields the claim for almost all $t\in [0,T]$.
    To show the assertion for an arbitrary $t\in [0,T]$, pick a sequence $(t_n)_n \subset [0,T]$ such that 
    the estimate in \eqref{eq:boundDzIeps} holds in $t_n$ for all $n\in \N$ and $t_n \to t$ as $n\to\infty$.
    Then the sequence $\mathrm{D}_z\II(\ell(t_n),z_{\epsilon}(t_n))$ is bounded in $\VV^*$ 
    and thus converges weakly to some $w \in \VV^*$ with $\|w\|_{\VV^*} \leq C(1 + \|\ell\|_{H^1(0,T;\VV^*)}^2)^{2q+1}$.
    Moreover, due to $z_\epsilon \in H^1(0,T;\ZZ)$ and 
    $\ell\in H^1(0,T;\VV^*)$, both $z_\epsilon$ and $\ell$ are continuous in time (with values in $\ZZ$ and 
    $\VV^*$, respectively). Therefore, $\mathrm{D}_z\II(\ell(t_n),z_\epsilon(t_n))$ converges in $\ZZ^*$ 
    to $\mathrm{D}_z\II(\ell(t),z_\epsilon(t))$ and the uniqueness of the limit implies 
    $w = \mathrm{D}_z\II(\ell(t),z_\epsilon(t))$, which finally yields the claim.
\end{proof}

Next, we verify additional convergence properties as $\delta$ tends to zero. 
Especially, the second one in Proposition~\ref{prop:deltaest} will be useful for the 
passage to the limit in the optimal control context in Section~\ref{sec:ocpapprox} below.

\begin{proposition}\label{prop:deltastrong}
    In the situation of Theorem~\ref{thm:exRISeps}, there holds that 
    $z_{\epsilon, \delta} \to z_\epsilon$ in $H^1(0,T;\VV)$.
\end{proposition}

\begin{proof}
    We return to \eqref{eq:energyineqepsdelta}, which we test with $v = 0$ to obtain
    \begin{equation*}
    \begin{aligned}
        & \int_0^T \epsilon\, \|z_{\epsilon, \delta}'(t)\|_\VV^2 + \RR(z_{\epsilon, \delta}'(t))\, \d t \\
        & \qquad \leq \int_0^T \dual{-\Drm_z \II(\ell_{\epsilon, \delta}(t), z_{\epsilon, \delta}(t))}{z_{\epsilon, \delta}'(t)}_{\ZZ^*, \ZZ}\,\d t \\
        & \qquad = - \II(\ell_{\epsilon, \delta}(T), z_{\epsilon, \delta}(T)) +  \II(\ell_{\epsilon, \delta}(0), z_0)
        - \int_0^T \dual{\ell_{\epsilon, \delta}'(t)}{z_{\epsilon, \delta}(t)}_{\VV^*, \VV} \, \d t.
    \end{aligned}
    \end{equation*}
    Because of the coercivity of $A$ and Lemma~\ref{lem:convF}, the reduced energy $\EE$ is 
    lower semicontinuous w.r.t.\ weak convergence in $\ZZ$.  
    Moreover, because of the weak continuity of the point evaluation in time, $\ell_{\epsilon, \delta}(T)$ converges weakly in $\VV$ 
    to $\ell_\epsilon(T)$. Thus, the weak convergence of $z_{\epsilon, \delta}$ in $H^1(0,T;\ZZ)$ from \eqref{eq:convergencedelta} 
    and the pointwise weak convergence  $z_{\epsilon, \delta}(T) \weak z_\epsilon(T)$ in $\ZZ$
    as consequence thereof as well as 
    the weak lower semicontinuity of $\RR$, cf.\ \eqref{eq:Repsconv}, and of the squared norm yield
    \begin{equation*}
    \begin{aligned}
        & \int_0^T \epsilon\, \|z_{\epsilon}'(t)\|_\VV^2 + \RR(z_{\epsilon}'(t))\, \d t \\
        & \quad \leq \liminf_{\delta \searrow 0}  \int_0^T \epsilon\, \|z_{\epsilon, \delta}'(t)\|_\VV^2 + \RR(z_{\epsilon, \delta}'(t))\, \d t \\
        & \quad \leq \limsup_{\delta \searrow 0}  \int_0^T \epsilon\, \|z_{\epsilon, \delta}'(t)\|_\VV^2 + \RR(z_{\epsilon, \delta}'(t))\, \d t \\
        & \quad \leq \limsup_{\delta \searrow 0} - \II(\ell_{\epsilon, \delta}(T), z_{\epsilon, \delta}(T)) +  \II(\ell_{\epsilon, \delta}(0), z_0)
        - \int_0^T \dual{\ell_{\epsilon, \delta}'(t)}{z_{\epsilon, \delta}(t)}_{\VV^*, \VV} \, \d t \\
        &\quad \leq - \II(\ell_{\epsilon}(T), z_{\epsilon}(T)) +  \II(\ell_{\epsilon}(0), z_0)
        - \int_0^T \dual{\ell_{\epsilon}'(t)}{z_{\epsilon}(t)}_{\VV^*, \VV} \, \d t \\
        & \quad = \int_0^T \epsilon\, \|z_{\epsilon}'(t)\|_\VV^2 + \RR(z_{\epsilon}'(t))\, \d t,
    \end{aligned}
    \end{equation*}
    where the last equality follows from \eqref{eq:VIeps} tested with $v = 0$ and $v = 2\, z_\epsilon$.
    Hence, we have that 
    \begin{equation*}
        \epsilon\, \| z_{\epsilon}'\|_\VV^2 + \int_0^T \RR(z_{\epsilon}'(t))\, \d t
        = \lim_{\delta\searrow 0}\epsilon\, \| z_{\epsilon,\delta}'\|_\VV^2 + \int_0^T \RR(z_{\epsilon, \delta}'(t))\, \d t
    \end{equation*}
    and, as both addends are weakly lower semicontinuous, this implies the norm convergence 
    $\| z_{\epsilon,\delta}'\|_\VV \to \| z_{\epsilon}'\|_\VV$. Together with the weak convergence by Theorem~\ref{thm:exRISeps}, 
    this in turn yields strong convergence of the time derivative, i.e., $z_{\epsilon, \delta}' \to z_\epsilon'$ in $L^2(0,T;\VV)$. 
    Since the function itself converge anyway by Theorem~\ref{thm:exRISeps} and the compact embedding of 
    $H^1(0,T;\ZZ) \embed L^2(0,T;\VV)$, this finishes the proof.
\end{proof}

\begin{proposition}\label{prop:deltaest}
    Let $\ell \in H^1(0,T;\VV^*)$ and $\epsilon, \delta > 0$ be given and 
    denote by $z_{\epsilon, \delta}$ the solution of \eqref{eq:RISepsdel} associated with $\ell$, 
    while $z_\epsilon$ is the solution of \eqref{eq:RISeps} for $\ell$. 
    Then there exists a constant $c_\epsilon>0$, depending on $\epsilon$, such that
    \begin{equation*}
        \sup_{t\in [0,T]} \|z_{\epsilon, \delta}(t) - z_\epsilon(t)\|_\ZZ \leq c_\epsilon \, \sqrt{\delta}.
    \end{equation*}
\end{proposition}

\begin{proof}
    The proof follows the lines of Lemma~\ref{lem:epsunique}.
    By testing the variational inequality for $z_\epsilon$ with $z_{\epsilon, \delta}$ and vice versa and adding 
    the arising inequalities, we obtain
    \begin{equation*}
    \begin{aligned}
        0 & \geq \epsilon \|z_{\epsilon,\delta}'(t)-z_\epsilon'(t)\|_\VV^2
        + \delta \dual{Az_{\epsilon, \delta}'}{z_{\epsilon, \delta}' - z_\epsilon'}_{\ZZ^*,\ZZ}\\
        & \quad +\langle \mathrm{D}_z \II(\ell(t), z_{\epsilon, \delta}(t))
        -\mathrm{D}_z \II(\ell(t), z_\epsilon(t)),z_{\epsilon, \delta}'(t)-z_\epsilon'(t)\rangle_{\ZZ^*,\ZZ}
        \quad \text{f.a.a.\ }t\in (0,T).
    \end{aligned}
    \end{equation*}
    Analogously to \eqref{eq:gronwallprep}, this implies     
    \begin{align*}
		& \frac{1}{2}\, \frac{\d}{\d t} \langle Az_{\epsilon, \delta}(t)-Az_\epsilon(t),z_{\epsilon, \delta}(t)-z_\epsilon(t)\rangle_{\ZZ^*,\ZZ}
		+ \epsilon\, \|z_{\epsilon, \delta}'(t)-z_\epsilon'(t)\|_\VV^2\\
		&\qquad \leq -\langle \mathrm{D}_z\FF(z_{\epsilon, \delta}(t))-\mathrm{D}_z\FF(z_\epsilon(t)),
		z_{\epsilon, \delta}'(t)-z_\epsilon'(t)\rangle_{\VV^*,\VV} 
		- \delta  \dual{A z_{\epsilon, \delta}'}{z_{\epsilon, \delta}' - z_\epsilon'}_{\ZZ^*,\ZZ}\\
		& \qquad \leq \epsilon\, \|z_{\epsilon, \delta}'(t)-z_\epsilon'(t)\|_\VV^2+ C_\epsilon\, \|z_{\epsilon, \delta}(t)-z_\epsilon(t)\|_\ZZ^2 
		+ \delta \, |z_{\epsilon, \delta}'(t)|_{\ZZ}\,|z_{\epsilon, \delta}' - z_\epsilon'|_\ZZ.
	\end{align*} 
	Now let $s\in [0,T]$ be arbitrary. Integrating the above inequality from $0$ to $s$ and using that 
	$z_\epsilon = z_{\epsilon, \delta} = z_0$ yield together with the coercivity of $A$ that
	\begin{equation*}
	\begin{aligned}
	    \frac{\alpha}{2}\,\|z_1(s)-z_2(s)\|_\ZZ^2 
	    & \leq C_\epsilon \int_0^s \|z_{\epsilon, \delta}(t)-z_\epsilon(t)\|_\ZZ^2 \,\d t \\
	    & \quad + \delta \|z_{\epsilon, \delta}\|_{H^1(0,T;\ZZ)} \,\|z_{\epsilon, \delta} - z_\epsilon\|_{H^1(0,T;\ZZ)}.
	\end{aligned}
	\end{equation*}
	Thanks to Lemma~\ref{lem:zdelH1Z} and Corollary~\ref{cor:boundszeps}, 
	$\|z_{\epsilon, \delta}\|_{H^1(0,T;\ZZ)}$ and $\|z_{\epsilon}\|_{H^1(0,T;\ZZ)}$ are bounded by a constant 
	independent of $\delta$ (but well depending on $\epsilon$). Gronwall's inequality then gives the result.
\end{proof}

We end this section by verifying a $W^{1,1}(0,T;\ZZ)$-bound, which is independent of $\epsilon$ and $\delta$
and will be used at the very end of this paper in the proof of Theorem~\ref{thm:ocpapproxepsdelta}.
For this purpose, let us define a tailored scalar product given by
\begin{equation}\label{eq:scaprod}
	\langle z_1,z_2\rangle_{\epsilon,\delta}
	:= \langle z_1,z_2\rangle_{\VV^*,\VV}+\frac{\delta}{\epsilon}\,\langle Az_1,z_2\rangle_{\ZZ^*,\ZZ}, \quad z_1,z_2\in\ZZ
\end{equation}
with its the associated norm $\|z\|_{\epsilon,\delta}:=\sqrt{\langle z,z\rangle_{\epsilon,\delta}}$, which is an equivalent norm 
on $\ZZ$ satisfying
\begin{equation}\label{eq:normequiv}
     \alpha\, \frac{\delta}{\epsilon}\, \|z\|_\ZZ^2 \leq \|z\|_{\epsilon, \delta}^2
     \leq \Big( C_\VV + \|A\|_{\LL(\ZZ, \ZZ^*)} \,\frac{\delta}{\epsilon}\Big) \|z\|_\ZZ^2 
     \quad \forall\, z\in \ZZ,
\end{equation}
where $C_\VV > 0$ is the embedding constant of $\ZZ \embed \VV$.

\begin{theorem}\label{thm:W11bound}
    Let $M > 0$ be given. Then there exists a constant $C_M>0$ depending on $M$, but not on $\epsilon,\delta>0$ 
	such that, for every $\ell \in H^1(0,T;\VV^*)$ with $\|\ell\|_{H^1(0,T;\VV^*)} \leq M$,
	the solution $z_{\epsilon,\delta}\in H^{2}(0,T;\ZZ)$ of \eqref{eq:RISepsdel} fulfills 
	\begin{equation}\label{eq:W11bound}
        \| z_{\epsilon,\delta}\|_{W^{1,1}(0,T;\ZZ)}\leq C_M.
	\end{equation}
\end{theorem}

\begin{proof}
    By inserting the form of the energy $\II$ in  \eqref{eq:z''} and rearranging terms, we obtain
    f.a.a.\ $t\in (0,T)$ that
	\begin{align*}
		& \langle A z_{\epsilon,\delta}'(t),z_{\epsilon,\delta}'(t)\rangle_{\ZZ^*,\ZZ} \\
		\quad & = - \epsilon \langle z_{\epsilon,\delta}''(t),z_{\epsilon,\delta}'(t)\rangle_{\epsilon,\delta}
		-\langle \mathrm{D}_z^2\FF(z_{\epsilon,\delta}(t))z_{\epsilon,\delta}'(t),z_{\epsilon,\delta}'(t)\rangle_{\VV^*,\VV} 
		+\langle \ell'(t),z_{\epsilon,\delta}'(t)\rangle_{\VV^*,\VV}.
	\end{align*}
    For the term involving $\Drm_z^2 \FF$, we employ an argument based on Ehrling's lemma and the estimates 
    from  \eqref{eq:F2} and \eqref{R3} to obtain for every $r > 0$ that 
    	\begin{align}\label{eq:estF''z'z'}
		\big|\langle \mathrm{D}_z^2\FF(z)v,v\rangle_{\VV^*,\VV}\big|
        	\leq \frac{\alpha}{2}\, \|v\|_\ZZ^2 + C_r\,\RR(v)\,\|v\|_\VV
        	\quad \forall\, v \in \ZZ, \, z \in B_\ZZ(0,r)
	\end{align}
	with a constant $C_r > 0$ depending only on $\alpha$ and $r>0$, 
    see, e.g., \cite[Lemma~1.1]{Kne19} or \cite[Lemma~3.1.6]{And25}.	
    Thanks to the a priori estimate in \eqref{eq:boundz}, $\sup_{t\in [0,T]}\|z_{\epsilon,\delta}(t)\|_\ZZ$ 
    is bounded by a constant $r_M > 0$ depending on $M$, but not on $\epsilon$ and $\delta$.
	Thus \eqref{eq:estF''z'z'} along with the coercivity of $A$ yields
    \begin{equation}\label{eq:zdotest}
    	\begin{aligned}
		& \frac{\alpha}{2}\|z_{\epsilon,\delta}'(t)\|_\ZZ^2 \\
		& \quad \leq -\epsilon \langle z_{\epsilon,\delta}''(t),z_{\epsilon,\delta}'(t)\rangle_{\epsilon,\delta}
		+C_{r_M}\RR(z_{\epsilon,\delta}'(t))\|z_{\epsilon,\delta}'(t)\|_\VV+\|\ell'(t)\|_{\VV^*}\|z_{\epsilon,\delta}'(t)\|_\VV.
	\end{aligned}     
    \end{equation}
    	Let us now pick an arbitrary $t \in \omega := \{ t\in [0,T] \colon \|z_{\epsilon, \delta}'(t)\|_{\ZZ} > 0\}$, 
    	where \eqref{eq:zdotest} is fulfilled (which is the case a.e.\ in $\omega$). 
    	Then dividing \eqref{eq:zdotest} by $\|z_{\epsilon,\delta}'(t)\|_{\epsilon,\delta}\neq 0$ 
    and using \eqref{eq:normequiv} results in
	\begin{align*}
	    \frac{\alpha}{2}\|z_{\epsilon,\delta}'(t)\|_\ZZ
        &\leq - \epsilon \Big( C_\VV + \|A\|_{\LL(\ZZ, \ZZ^*)} \,\frac{\delta}{\epsilon}\Big)
        \bigg\langle z_{\epsilon,\delta}''(t),\frac{z_{\epsilon,\delta}'(t)}
        {\|z_{\epsilon,\delta}'(t)\|_{\epsilon,\delta}}\bigg\rangle_{\epsilon,\delta} \\
        &\qquad  + C_\VV\, \big( C_{r_M} \RR(z_{\epsilon,\delta}'(t))+\|\ell'(t)\|_{\VV^*}\big),
	\end{align*}
	where $C_\VV$ again denoted the embedding constant of $\ZZ\embed \VV$. Now integrating over $[0,T]$
    and employing Lemma~\ref{lem:diffnorm} from the appendix and Lemma~\ref{lem:z'(0)} gives
	\begin{align*}
        \frac{\alpha}{2}\int_0^T\|z_{\epsilon,\delta}'(r)\|_\ZZ \,\d r 
        &= \frac{\alpha}{2}\int_\omega\|z_{\epsilon,\delta}'(r)\|_\ZZ\, \d r \\
		& \leq - \epsilon \Big( C_\VV + \|A\|_{\LL(\ZZ, \ZZ^*)} \,\frac{\delta}{\epsilon}\Big) 
		\int_\omega  \bigg\langle z_{\epsilon,\delta}''(r),\frac{z_{\epsilon,\delta}'(r)}{\|z_{\epsilon,\delta}'(r)\|_{\epsilon,\delta}}\bigg\rangle_{\epsilon,\delta}\d r \\
        & \quad + C_\VV \,C_{r_M} \int_\omega\RR(z_{\epsilon,\delta}'(r))\,\d r + C_\VV \int_\omega\|\ell'(r)\|_{\VV^*}\,\d r\\
        & \leq  - \epsilon \Big( C_\VV + \|A\|_{\LL(\ZZ, \ZZ^*)} \,\frac{\delta}{\epsilon}\Big)  
        \Big(\|z_{\epsilon,\delta}'(T)\|_{\epsilon,\delta}-\|z_{\epsilon,\delta}'(0)\|_{\epsilon,\delta}\Big) \\
        & \quad + C_\VV \,C_{r_M} \int_0^T\RR(z_{\epsilon,\delta}'(r))\,\d r + C_\VV \|\ell'\|_{L^1(0,T;\VV^*)}\\
        & \leq C_\VV \,C_{r_M} \, C(M^2 + 1)  + C_\VV \, \sqrt{T} \, M,
	\end{align*}
	where we used the estimate from Remark~\ref{rem:Rest} for the last inequality.
    In combination with \eqref{eq:boundz}, this yields the assertion.
\end{proof}

\section{Existence of Optimal Controls}\label{sec:optcontrol}

We are now in the position to state a mathematically rigorous formulation of the optimal control problem \eqref{eq:ocp0}:
\begin{equation}\tag{OCP}\label{eq:optcontrol}
    \left\{ \quad
    \begin{aligned}
	    \min \quad &  J(S,\hat{z},\ell) :=j(\hat{z}(S))+\frac{\beta}{2} \, \|\ell\|_{H^1(0,T;\VV^*)}^2\\
        \text{s.t.} \quad & \ell\in H^1(0,T;\VV^*),\quad(S,\hat{t},\hat{z})\in\LL(z_0,\ell),\\
		&-\mathrm{D}_z\II(\ell(0),z_0)\in\partial \RR(0),~-\mathrm{D}_z\II(\ell(T),\hat z(S))\in\partial \RR(0),
	\end{aligned}
	\right.
\end{equation}
where $\LL(z_0,\ell)$ denotes the set of parametrized BV solutions according to Definition~\ref{def:paramsol} 
associated with $z_0$ and $\ell$, i.e., 
\begin{align*}
	\LL(z_0,\ell) := 
	\big\{ (S,\hat{t},\hat{z})\in\, & [T,\infty) \times W^{1,\infty}(0,S) \times \mathrm{AC}^\infty([0,S];\RR)\cap L^\infty(0,S;\ZZ): \\
	& (S,\hat{t},\hat{z}) \text{ is a parametrized BV solution associated with $\ell$}\big\}.
\end{align*}
Compared to the formal definition of \eqref{eq:ocp0} in the introduction, we observe that we have included two additional constraints 
to the problem. The first one concerns the initial state and enforces the solution to be locally stable at the beginning. 
It ensures that the condition in \eqref{eq:init} is fulfilled for all feasible controls $\ell$ such that the results of the previous section apply, 
in particular Lemma~\ref{lem:z'(0)}.

The second additional constraint is motivated by application, as it guarantees that the state is locally stable at end time, too. 
This means that the trajectory does not end in a viscous jump and consequently, the final state can be seen in physical time, 
which certainly makes sense from an application point of view.

The existence of a globally optimal solution of \eqref{eq:optcontrol} is an immediate consequence 
of the following compactness result concerning $\mathfrak{p}$-parametrized BV solutions from \cite{KT23}, 
whose involved proof is anything but obvious, see also \cite{Tho22} for details. 
 
 \begin{theorem}[{Compactness of the feasible set, \cite[Thm.~3.12]{KT23}}]\label{thm:compact}
 	Let $z_0\in\ZZ$ and $\rho > \|z_0\|_\ZZ$ be given and define the set 
    \begin{equation}
    \begin{aligned}
 	    M_\rho:=\Big\{(S,\hat t,\hat z,\ell): &~(S,\hat t,\hat z)\in\LL(z_0,\ell),  \\[-1ex]
 	    &-\mathrm{D}_z\II(\ell(0),z_0)\in\VV^* ,
 	    \; \|z_0\|_\ZZ+\|\ell\|_{H^1(0,T;\VV^*)}\leq \rho\Big\}.
    \end{aligned}
    \end{equation}
 	Then $M_\rho$ is compact in the following sense:
 	
    For every sequence $(S_n,\hat t_n,\hat z_n,\ell_n)_{n\in\N}\subset M_\rho$ 
    there exists a (not relabeled) subsequence and a limit $(S,\hat t,\hat z,\ell)\in M_\rho$ with
	\begin{gather}
		S_n\to S \text{ in }\R,\quad \hat t_n\rightharpoonup^*\hat t \text{ in } W^{1,\infty}(0,S),\quad \hat t (S)=T,\\
		\ell_n\rightharpoonup \ell \text{ in } H^1(0,T;\VV^*), \quad \hat z_n\to \hat z \text{ in } C([0,S];\VV), \label{eq:convell} \\
        \hat z_n(S_n) \to \hat z(S) \text{ in } \VV,
        \quad 
		\mathrm{D}_z\EE(\hat z_n(S_n))\rightharpoonup\mathrm{D}_z\EE(\hat z(S))~\text{in}~\VV^*. \label{eq:convatS}
	\end{gather}
	Herein, the functions $\hat z_n$ and $\hat t_n$ are constantly extended by their value at $S_n$, if $S > S_n$.
\end{theorem}
 
\begin{remark}
    We point out that the above theorem does not provide all results of \cite[Thm.~3.12]{KT23}, but only those 
    needed for the existence of optimal solutions proven below. 
    We moreover emphasize that the second convergence in \eqref{eq:convatS} is not part of \cite[Thm.~3.12]{KT23}, 
    but it can be deduced from $\hat z_n(S_n) \to \hat z(S)$ in $\VV$ and the boundedness of $\{\hat z_n(S_n)\}$ in 
    $\ZZ$, which implies $\hat z_n(S_n) \weakly \hat z(S)$ in $\ZZ$,  giving in turn  the weak convergence of 
    $\mathrm{D}_z\EE(\hat z_n(S_n))$. We refer to \cite[Remark~4.1.2]{And25} for details.
\end{remark}
 
 \begin{theorem}\label{thm:existocp0}
 	Let $z_0\in\ZZ$ be given and assume that there exists $\ell_0\in\VV^*$ such that 
    \begin{equation}\label{eq:assuell0}
        -\mathrm{D}_z\II(\ell_0,z_0)\in \partial\RR(0). 
    \end{equation}
     Then \eqref{eq:optcontrol} has a globally optimal solution.
 \end{theorem}
 
 \begin{proof}
 	First of all, the feasible set of \eqref{eq:optcontrol} is non-empty, 
 	because, thanks to \eqref{eq:assuell0}, 
 	the constant function $z\equiv z_0$ together with $\hat t=\text{id}$ (identity) and $S=T$ 
 	is a parametrized solution associated with $\ell\equiv \ell_0$, 
	i.e., $(T,\text{id},z_0)\in\LL(z_0, \ell_0)$. Note that the end time constraint in \eqref{eq:optcontrol} is also fulfilled
	due to \eqref{eq:assuell0}.

    The rest of the proof is a direct consequence of Theorem~\ref{thm:compact} by applying the 
    direct method of the calculus of variations. For convenience of the reader, we explain the details.
    Let $(S_n,\hat t_n,\hat z_n,\ell_n)_{n\in\N}$ be an infimal sequence, i.e., 
    \begin{equation*}
        \lim_{n\to\infty} J(S_n,\hat{z}_n,\ell_n)=\inf J(S,\hat{z},\ell)=: I,
    \end{equation*}
    where the infimum is considered over all feasible points of \eqref{eq:optcontrol} 
    and is finite since $j:\VV\to\R$ is assumed to be bounded from below. 
    Therefore, thanks to the Tikhonov term in the objective,
    we have that  $\sup_{n\in\N}\|\ell_n\|_{H^1(0,T;\VV^*)}\leq C$. 
    Moreover, $-\mathrm{D}_z\II(\ell_n(0),z_0)\in \partial\RR(0) \subset \VV^*$ for all $n\in\N$ by Lemma \ref{lem:dR(0)bound}
    and thus, we are allowed to apply Theorem~\ref{thm:compact}. 
    Hence, there is a (not relabeled) subsequence of the infimal sequence converging 
    to a limit satisfying $(S,\hat{t},\hat{z})\in\LL(z_0,\ell)$. 
    From \eqref{eq:convell} and the weak continuity of the point evaluation, 
    we infer that $\ell_n(0)\rightharpoonup\ell(0)$ in $\VV^*$ such that 
    $-\mathrm{D}_z\II(\ell(0),z_0)\in\partial\RR(0)$ by exploiting the convexity and closedness of the subdifferential. 
    In the same way we 	deduce from $\ell_n(T)\rightharpoonup\ell(T)$ in $\VV^*$ and 
    \eqref{eq:convatS} that $-\mathrm{D}_z\II(\ell(T),\hat z(S))\in\partial\RR(0)$, too. Thus, the limit is feasible 
	for \eqref{eq:optcontrol}. Finally, thanks to \eqref{eq:convatS}, the continuity of $j$ by assumption
	and weak lower semicontinuity of the norm in combination with \eqref{eq:convell}, we infer
    \begin{equation}\label{eq:optimality}
	\begin{aligned}
		J(S,\hat{z},\ell) &=j(\hat{z}(S))+\frac{\beta}{2} \, \|\ell\|_{H^1(0,T;\VV^*)}^2 \\
		&\leq \liminf_{n\to\infty}j(\hat{z}_n(S_n))+\frac{\beta}{2} \, \|\ell_n\|_{H^1(0,T;\VV^*)}^2=I
	\end{aligned}
    \end{equation}
	such that $(S,\hat t,\hat z,\ell)$ is indeed a minimizer of \eqref{eq:optcontrol}.
 \end{proof}

Next we turn to the $\epsilon$-viscous regularization of \eqref{eq:optcontrol}. 
The mathematically rigorous formulation of \eqref{eq:ocpeps} reads as follows:
\begin{equation}\tag{$\mathrm{vOCP}_\epsilon$}\label{eq:optcontroleps}
    \left\{ \quad 
    \begin{aligned}
        \min \quad &  J_\epsilon(z_\epsilon,\ell):=j(z_\epsilon(T))+\frac{\beta}{2} \, \|\ell\|_{H^1(0,T;\VV^*)}^2\\
        \text{s.t.} \quad & \ell\in H^1(0,T;\VV^*),\quad z_\epsilon\in H^1(0,T;\ZZ),\\
        & 0\in\partial\RR_\epsilon(z_\epsilon '(t))+\mathrm{D}_z\II(\ell(t),z_\epsilon(t)),\quad z_\epsilon(0)=z_0,\\
        & -\mathrm{D}_z\II(\ell(0),z_0)\in\partial\RR(0),\\
        & \dist_{\ZZ^*}(-\mathrm{D}_z\II(\ell(T),z_\epsilon(T)),\partial \RR(0))\leq \epsilon^{\frac{1}{4}},
    \end{aligned}    
    \right.
\end{equation}
where $j:\VV\to \R$ and $\beta>0$ are the same as in \eqref{eq:optcontrol}. 
Note that, while we again require local stability at the initial time, the 
end time constraint is relaxed. This will be of importance for the approximation results in Section~\ref{sec:ocpapprox}.

\begin{theorem}\label{thm:ocpepsexist}
    Let $z_0\in \ZZ$ be given and assume that there exists an $\ell_0\in \VV^*$ such that \eqref{eq:assuell0} is satisfied.
    Then, for every $\epsilon>0$, there exists a globally optimal solution 
    $(z_\epsilon^*,\ell_\epsilon^*)\in H^1(0,T;\ZZ)\times H^1(0,T;\VV^*)$ of \eqref{eq:optcontroleps}.
\end{theorem}

\begin{proof}
    The proof again follows the classical direct method of the calculus of variations. 
    For the sake of completeness, we present the details. To ease notation, we suppress the subscript $\epsilon$ 
    in the following, as $\epsilon$ is fixed here. 
    First, similarly to the proof of Theorem~\ref{thm:existocp0}, the constant functions $z \equiv z_0$ and 
    $\ell\equiv \ell_0$ are a feasible point of \eqref{eq:optcontroleps} due to \eqref{eq:assuell0}. 
    Thus there exists an inifimizing sequence denoted by $(z_n,\ell_n)_{n\in\N}$. 
    As in the proof of Theorem~\ref{thm:existocp0}, one shows that the sequence of loads is bounded in 
    $H^1(0,T;\VV^*)$  due to the Tikhonov term.
    Together with the a priori bound \eqref{eq:boundDzeps} ($\epsilon$ is fixed here), 
    this ensures the existence of a weakly converging subsequence (denoted by the same index) such that
	\begin{equation}\label{eq:weakconvzn}
		(z_n,\ell_n)\rightharpoonup (z^*,\ell^*) \text{ in } H^1(0,T;\ZZ)\times H^1(0,T;\VV^*).
	\end{equation}
	
	To verify the feasibility of $(z^*,\ell^*)$, first note that the weak continuity of the point evaluation 
	in time from $H^1(0,T;\ZZ)\times H^1(0,T;\VV^*)$ to $\ZZ\times \VV^*$ implies 
	\begin{equation}\label{eq:pwkonvzn}
		z_n(t)\rightharpoonup z^*(t) \text{ in } \ZZ,\quad \ell_n(t)\to \ell^*(t) \text{ in } \VV^*
		\quad \forall\, t\in [0,T]
	\end{equation}	
	and thus, the initial condition is satisfied in the limit. Moreover, \eqref{eq:pwkonvzn} implies that 
	$\mathrm{D}_z\II(\ell_n(t),z_n(t))$ converges weakly in $\ZZ^*$ for all $t\in [0,T]$ and, from
	the a priori bound in \eqref{eq:boundDzIeps}, we deduce that this pointwise convergence also holds in $\VV^*$.
	Therefore, the weak closedness of $\partial\RR(0)$ 
	gives $-\mathrm{D}_z\II(\ell^*(0),z_0)\in\partial\RR(0)$ and the weak lower semicontinuity of the distance 
	yield 
    \begin{equation*}
        \dist_{\ZZ^*}(-\mathrm{D}_z\II(\ell^*(T),z^*(T)),\partial\RR(0))\leq \epsilon^{\frac{1}{4}}
    \end{equation*}    	
	so that the additional constraints in \eqref{eq:optcontroleps} are fulfilled, too.
	
    It remains to verify that $(z^*, \ell^*)$ satisfies \eqref{eq:RISeps}, which can be done along the lines 
    of the proof of Theorem \ref{thm:exRISeps}. 
    We start with the equivalent reformulation of \eqref{eq:RISeps} in terms of an energy inequality
    that is fulfilled by $(z_n, \ell_n)$ for every $t\in [0,T]$ and every $v\in L^2(0,T;\ZZ)$:
    	\begin{equation}\label{eq:energyineqzn}
	\begin{aligned}
        	\int_0^t\RR(v(r))\,\d r \geq 
        & \int_0^t\RR(z_n'(r)) \,\d r + \int_0^t \epsilon\langle z_n'(r),z_n'(r)-v(r)\rangle_{\VV^*,\VV} \,\d r\\
		&+ \int_0^t \langle \mathrm{D}_z\II(\ell_n(r),z_n(r)),z_n'(r)-v(r)\rangle_{\ZZ^*,\ZZ}\,\d r.
	\end{aligned}
	\end{equation}
    Using the weak convergence in \eqref{eq:weakconvzn}, the first two terms on the right hand side can 
    be treated as in \eqref{eq:Repsconv} and \eqref{eq:epsconv}.
    For the third term, we can again apply Lemma~\ref{lem:limAF} to obtain a $\liminf$-inequality analogously to 
    \eqref{eq:DzIepsconv}. All in all this shows that \eqref{eq:energyineqzn} transfers to the limit, 
    which shows the feasibility of $(z^*, \ell^*)$.

    Finally, its optimality follows analogously to \eqref{eq:optimality}, by employing the continuity of $j$
    and the weak lower semicontinuity of the squared norm.
\end{proof}

\begin{remark}\label{rem:ocpexistVstar}
    By exactly the same arguments, one proves that \eqref{eq:optcontroleps} also admits a solution, if 
    the end time constraint is replaced by 
    \begin{equation*}
         \dist_{\VV^*}(-\mathrm{D}_z\II(\ell(T),z_\epsilon(T)),\partial \RR(0))\leq \epsilon^{\frac{1}{4}}.
    \end{equation*}
    As the above proof shows, the a priori bound in \eqref{eq:boundDzIeps} implies that 
    $-\mathrm{D}_z\II(\ell(T),z_\epsilon(T))$ converges weakly in $\VV^*$, too, and thus the 
    weak lower semicontinuity (this time w.r.t.\ weak convergence in $\VV^*$) implies the feasibility of the limit. 
    However, as it will turn out in Section~\ref{sec:ocpapprox}, taking the distance w.r.t.\ the $\ZZ^*$-norm 
    is favorable in order to show the convergence of the second regularization, cf.~Proposition~\ref{prop:ocpconvdelta} below.
\end{remark}

In addition to \eqref{eq:optcontroleps}, we consider yet another regularized optimal control problem, 
where the only difference is that \eqref{eq:RISeps} is replaced by the double viscous regularization \eqref{eq:RISepsdel}, 
which can be equivalently reformulated as \eqref{eq:ODE} as seen in the proof of Lemma~\ref{lem:exRIS}:
\begin{equation}\tag{$\mathrm{vOCP}_{\epsilon, \delta}$}\label{eq:optcontrolepsdelta}
    \left\{ \quad 
    \begin{aligned}
        \min \quad &  J_\epsilon(z_{\epsilon, \delta},\ell):=j(z_{\epsilon, \delta}(T))+\frac{\beta}{2} \, \|\ell\|_{H^1(0,T;\VV^*)}^2\\
        \text{s.t.} \quad & \ell\in H^1(0,T;\VV^*),\quad z_{\epsilon, \delta} \in H^2(0,T;\ZZ),\\
        & z_{\epsilon,\delta}'(t) =\partial\RR_{\epsilon,\delta}^*\big(-\Drm_z\II(\ell(t),z_{\epsilon,\delta}(t))\big),
        \quad z_{\epsilon, \delta}(0)=z_0,\\
        & -\mathrm{D}_z\II(\ell(0),z_0)\in\partial\RR(0),\\
        & \dist_{\ZZ^*}(-\mathrm{D}_z\II(\ell(T),z_{\epsilon, \delta}(T)),\partial \RR(0))\leq \epsilon^{\frac{1}{4}} + \delta^{\frac{1}{4}}.
    \end{aligned}    
    \right.
\end{equation}

\begin{remark}\label{rem:ode}
    As indicated in the introduction, the advantage of \eqref{eq:optcontrolepsdelta} is that it is an 
    optimal control problem governed by a non-smooth ODE in Hilbert space. 
    Therefore standard smoothing techniques such as for instance a Moreau-Yosida regularization of 
    $\RR_{\epsilon,\delta}^*$ can be applied to obtain a smooth optimal control problem that is amenable to 
    the classical adjoint calculus for smooth optimal control problems. We refer to \cite{Wac15, SWW17, GW18} 
    for details on smoothing procedures for problems of this type.
\end{remark}

\begin{theorem}
    Let $z_0\in \ZZ$ be given and assume that there exists an $\ell_0\in \VV^*$ such that \eqref{eq:assuell0} 
    is fulfilled.
    Then, for every $\epsilon>0$ and every $\delta > 0$, \eqref{eq:optcontrolepsdelta} admits a globally optimal solution.
\end{theorem}

\begin{proof}
    The proof is exactly along the lines of the proof of Theorem~\ref{thm:ocpepsexist}, 
    with the only difference that the energy inequality associated with \eqref{eq:RISepsdel} additionally contains 
    the term $\int_0^t\delta\langle Az'(r),z'(r)-v(r)\rangle_{\ZZ^*,\ZZ}\,\d r$, cf.\ \eqref{eq:energyineqepsdelta}.
    Due to the coercivity of $A$, this term is clearly lower semicontinuous 
    w.r.t.\ weak convergence in $H^1(0,T;\ZZ)$ and therefore one can pass to the limit in the energy inequality 
    completely analogously to the proof of Theorem~\ref{thm:ocpepsexist}.
\end{proof}

\section{Reverse Approximation}\label{sec:revapp}

This section is devoted to the major challenge in proving the approximability of 
\eqref{eq:optcontrol} by its viscous regularization \eqref{eq:optcontroleps}, namely
the construction of a recovery sequence for optimal solutions of \eqref{eq:optcontrol}.
As indicated in the introduction, there is in general no hope that every paramterized BV solution 
can be approximated via viscous regularization. In the optimal control context however, we have more 
flexibility, since we can vary the control in forms of the external loads, too.
Unfortunately, the objective involves the $H^1(0,T;\VV^*)$-norm of the control
(which is essential for the existence of optimal solutions as explained in the introduction) 
and therefore, we need strong convergence of the load-part of the recovery sequence in $H^1(0,T;\VV^*)$. 
In order to achieve this, we follow an idea of \cite{KMS22}, 
where the finite dimensional case with $\XX=\VV=\ZZ=\R^n$ is investigated.
The essential idea is to enrich the energy by a penalization term, which on the one hand vanishes in the limit 
and on the other hand guarantees the uniform convexity of the energy. 
However, in order to prove that the penalization tends to zero in the limit, which ultimately implies the 
converge of the loads in $H^1(0,T;\VV^*)$, see the proof of Theorem~\ref{thm:revapp}, 
we need the following additional assumption:

\begin{assumption}\label{ass:diffsol}
	There exists at least one optimal solution 
	$(S,\hat t,\hat z,\ell)\in [T,\infty)\times W^{1,\infty}(0,S;\R)\times H^1(0,S;\ZZ)\times H^1(0,T;\VV^*)$ 
	of \eqref{eq:optcontrol} that satisfies 
	\begin{align}\label{eq:t'>del}
		\hat t'(s)\geq \rho \text{ for almost all } s\in(0,S)
	\end{align}
	with $0<\rho\leq 1$.
\end{assumption}

\begin{remark}
   We have to admit that  Assumption~\ref{ass:diffsol} is very restrictive, as it implies that \eqref{eq:optcontrol} 
   has at least one optimal solution with an optimal state $\hat z$ that is continuous in time. 
   Due to the non-convexity of the energy, we know however that this cannot be guaranteed in general, 
   which is, in essence, the motivation for the variety of alternative solutions concepts such at the parametrized  BV solution
   from Definition~\ref{def:paramsol}.
    Yet, we underline that we do not require that every global minimizer admits an optimal state that is continuous in time;
    we only need this property for at least one global minimizer.
\end{remark}

The condition \eqref{eq:t'>del} implies that $\hat t $ is strictly monotone increasing. 
Hence, its inverse function $\hat t^{-1}:[0,T]\to[0,S]$ exists and 
satisfies $1\leq \frac{\d}{\d t}\hat t^{-1}(t)<\frac{1}{\rho}$ for almost all $t\in(0,T)$,
since additionally $\hat t'(s)\leq 1$ almost everywhere by \eqref{eq:normalization}. 
Thus we have $\hat t^{-1}\in W^{1,\infty}(0,T)$.
With the inverse of $\hat t$ at hand, we define 
\begin{equation}\label{eq:tildezdef}
    \tilde z:=\hat z \circ \hat t^{-1},
\end{equation}
which is just the optimal state in physical time.
This state is even a differential solution, as proven in Appendix~\ref{sec:diffsol}.

\begin{lemma}\label{lem:diffsol}
    The transformed state $\tilde z$ is an element of $H^1(0,T;\ZZ)$ and 
    a differential solution of \eqref{eq:ris}, i.e., it satisfies 
    \begin{equation}\label{eq:RIStilde}
	    0\in\partial\RR(\tilde z'(t))+\mathrm{D}_z\II(\ell(t),\tilde z(t)) \; \text{ a.e.\ in }(0,T), \quad \tilde z(0)=z_0.
    \end{equation}
\end{lemma}

Given $\tilde z$ from \eqref{eq:tildezdef}, we construct the penalization by adding a quadratic penalty term in 
$\VV$ to the energy:
\begin{equation*}
    \begin{aligned}
        & \JJ_\eta : [0,T] \times \VV^* \times \ZZ \to \R,\\
    	    & \JJ_\eta(t,\ell,z):=\frac{1}{2}\langle Az,z\rangle_{\ZZ^*,\ZZ}+\FF(z)
	    -\langle \ell,z\rangle_{\VV^*,\VV}+\frac{\eta}{2}\|z-\tilde z(t)\|_\VV^2, 
    \end{aligned}
\end{equation*}
where $\eta > 0$ will be chosen so large such that $\JJ_\eta$ becomes uniformly convex along the solution trajectory, 
see Lemma~\ref{lem:Eetaconvex} below. 
For the rest of this section, it will be convenient 
to define the following energy functionals associated with $\JJ_\eta$:
\begin{equation}\label{eq:etaenergies}
\begin{aligned}
	\EE_\eta : \ZZ\ni z & \mapsto \EE(z)+\frac{\eta}{2}\|z\|_\VV^2 \in \R, \\
	\II_\eta : \VV^* \times \ZZ \ni (\ell,z) & \mapsto \II(\ell,z)+\frac{\eta}{2}\|z\|_\VV^2 \in \R.
\end{aligned}
\end{equation}
Given the penalized energy, we construct the elements of our recovery sequence as solutions to 
the following viscous problem:
\begin{equation}\label{eq:RISepseta}
	0\in \partial\RR_\epsilon(z_\epsilon'(t))+\mathrm{D}_z\JJ_\eta(t,\ell(t),z_\epsilon(t)), \quad z_\epsilon(0)=z_0.
\end{equation}
Straight forward computation shows that, expressed in terms of the energy $\II_\eta$, the above viscous system reads
\begin{equation}\label{eq:RISepseta2}
	0\in \partial\RR_\epsilon(z_\epsilon'(t))+\mathrm{D}_z\II_\eta(\ell(t)+\eta\tilde z(t) ,z_\epsilon(t)), \quad z_\epsilon(0)=z_0.
\end{equation}
In order to apply the results from Section~\ref{sec:visreg} to this equation, first note 
that, by Lemma~\ref{lem:diffsol}, there holds that $\tilde z\in H^1(0,T;\ZZ)$ such that $\ell +\eta \tilde z \in H^1(0,T;\VV^*)$.
Moreover, it is easily seen that the nonlinearity 
\begin{equation*}
    \FF_\eta(z) := \FF(z)+\frac{\eta}{2}\|z\|_\VV^2
\end{equation*}
complies with our standing assumptions from Section \ref{sec:energyass}, 
i.e., \eqref{eq:F1}--\eqref{eq:D2Flip}. Therefore, the results of Section~\ref{sec:visreg} are indeed applicable. 
This in particular implies that the corresponding $(\delta, \epsilon)$-regularized system given by
\begin{equation}\label{eq:RISepsdeltaeta}
	0\in \partial\RR_{\epsilon,\delta}(z_{\epsilon,\delta}'(t))+\mathrm{D}_z\II_\eta(\ell(t) + \eta \,\tilde z(t), z_{\epsilon,\delta}(t)), 
	\quad z_{\epsilon,\delta}(0)=z_0.
\end{equation}
admits a unique solution, cf.~Lemma~\ref{lem:exRIS}. Furthermore, 
by passing to the limit $\delta\searrow 0$, we deduce from Theorem \ref{thm:exRISeps} that 
\begin{equation}\label{eq:H1Vconveta}
	z_{\epsilon,\delta}\rightharpoonup z_{\epsilon} \text{ in } H^1(0,T;\ZZ)
	\quad\text{and}\quad
	z_{\epsilon,\delta}(t)\rightharpoonup z_\epsilon(t) \text{ in } \ZZ \quad \forall\, t\in [0,T],
\end{equation}
where $z_\epsilon$ is a solution of \eqref{eq:RISepseta2}.
Moreover, $z_{\epsilon, \delta}$ and $z_\epsilon$ satisfy the a priori bounds from 
Section~\ref{sec:visreg}, especially the pointwise bounds from Lemma~\ref{lem:exRIS} and Corollary~\ref{cor:boundszeps}.
Let us investigate how these bounds depend on the penalization parameter $\eta$. 

\begin{lemma}\label{lem:uniboundeta}
    There exists a constant $C>0$ independent of $\epsilon$, $\delta$, and $\eta$ such that 
    \begin{equation}\label{eq:uniboundeta}
        \|z_{\epsilon, \delta}(t)\|_{\ZZ} + \|z_{\epsilon}(t)\|_{\ZZ} \leq 
        C ( \eta + 1 ) \quad \forall\, t\in [0,T].
    \end{equation}
\end{lemma}

\begin{proof}
    In principle, the proof follows the lines of the one of Lemma~\ref{lem:exRIS}. 
    First we note that Lemma~\ref{lem:lambdamu} and Lemma~\ref{lem:nu} also hold with 
    $\II_\eta$ instead of $\II$ with the same constants $\lambda$, $\mu$, and $\nu$, since 
    $\II_\eta$ and $\II$ differ only in the nonlinearity and the modified nonlinearity $\FF_\eta$ is also non-negative,
    which is the only property of $\FF$ used in the proofs of Lemmas~\ref{lem:lambdamu} and \ref{lem:nu}, 
    cf.\ Remark~\ref{rem:lambdamunu}.
    
    To ease notation, we abbreviate $z_{\epsilon, \delta}$ simply by $z$.
    Our starting point is again to rewrite \eqref{eq:RISepsdeltaeta} in terms of 
    the Fenchel-Young equality and apply the chain rule to obtain 
    \begin{equation*}
    \begin{aligned}
        	& \RR_{\epsilon,\delta}(z'(t))+\RR_{\epsilon,\delta}^*(-\Drm_z\II_\eta(\ell(t) + \eta\, \tilde z(t),z(t))) \\
	    & \qquad = -\frac{\d}{\d t}\II_\eta(\ell(t) + \eta \tilde z(t),z(t)) - \eta \dual{\tilde z'(t)}{z(t)}_{\VV^*, \VV} 
	    - \langle\ell'(t),z(t)\rangle_{\VV^*,\VV}.
    \end{aligned}
    \end{equation*}
    Thanks to the non-negativity of $\RR_{\epsilon, \delta}$ and its Fenechel conjugate,
    integration with $t\in [0,T]$ arbitrary yields 
    \begin{equation*}
    \begin{aligned}
        	\II_\eta(\ell(t) + \eta \,\tilde z(t),z(t))
        & \leq  \II_\eta(\ell(0) + \eta\, z_0,z_0) \\
        &\quad - \eta  \int_0^t \dual{\tilde z'(r)}{z(r)}_{\VV^*, \VV} \,\d r 
        -  \int_0^t  \langle\ell'(r),z(r)\rangle_{\VV^*,\VV}\,\d r.
    \end{aligned}
    \end{equation*}
    For the last term we again apply Young's inequality in combination with Lemma~\ref{lem:lambdamu} 
    analogously to \eqref{eq:estintright}. Together with $\II_\eta(\ell + \eta \,\tilde z,z) = \II_\eta(\ell, z) - \eta \dual{\tilde z}{z}_{\VV^*, \VV}$,
    this yields    
    \begin{equation*}
    \begin{aligned}
        	\II_\eta(\ell(t),z(t))
	    & \leq \eta \dual{\tilde z(t)}{z(t)}_{\VV^*, \VV} + \II_\eta(\ell(0),z_0) - \eta \,\|z_0\|_\VV^2 +
	    \frac{1+\mu}{2}\, \|\ell\|_{H^1(0,t;\VV^*)}^2 \\[-1ex]
        & \quad  + \eta \, \| z \|_{L^2(0,T;\VV)} \, \|\tilde z\|_{H^1(0,T;\VV)}   
	    + \frac{\lambda}{2} \int_0^t \II_\eta(\ell(r),z(r))\,\d r \\
	    & \leq  \II(\ell(0),z_0) + \frac{1+\mu}{2}\, \|\ell\|_{H^1(0,t;\VV^*)}^2 +  \frac{\lambda}{2} \int_0^t \II_\eta(\ell(r),z(r))\,\d r \\
	    & \quad + (C + \sqrt{T})\,\eta \sup_{t\in [0,T]} \|z(t)\|_\VV\,\|\tilde z\|_{H^1(0,T;\VV)} ,
    \end{aligned}
    \end{equation*}
    where we used that $\II_\eta(\ell(0),z_0) = \II(\ell(0),z_0) + \frac{\eta}{2} \|z_0\|_\VV^2$ and
    employed the continuous embedding $H^1(0,T;\VV) \embed L^\infty(0,T;\VV)$ with embedding constant $C>0$.
    Gronwall's inequality then implies
    \begin{align*}
	    \II_\eta(\ell(t),z(t))
	    \leq 
        e^{\frac{1}{2}\lambda T} \Big( & \II(\ell(0),z_0) + \frac{1+\mu}{2}\|\ell\|_{H^1(0,T;\VV^*)}^2 \\
    	    &+ (C + \sqrt{T}) \,\eta \sup_{t\in [0,T]} \|z(t)\|_\VV \, \|\tilde z\|_{H^1(0,T;\VV)} \Big).
    \end{align*}
    Using again Lemma~\ref{lem:lambdamu} and the continuous embedding $\ZZ\embed \VV$, this in turn gives
    \begin{equation*}
    \begin{aligned}
        \|z(t)\|_\ZZ^2 
        & \leq  \lambda \, \II_\eta(\ell(t),z(t)) + \mu \, \|\ell(t)\|_{\VV^*}^2\\
        & \leq 
        \begin{aligned}[t]
            \lambda \, e^{\frac{1}{2}\lambda T} \Big( & \II(\ell(0),z_0) + \frac{1+\mu}{2}\|\ell\|_{H^1(0,T;\VV^*)}^2 \\[-1ex]
        	    & + (C + \sqrt{T}) \,\eta \sup_{t\in [0,T]} \|z(t)\|_\VV\, \|\tilde z\|_{H^1(0,T;\VV)} \Big) + \mu \, \|\ell\|_{L^\infty(0,T;\VV^*)}^2 
        \end{aligned} \\
        &=: C_1 \,\eta\sup_{t\in [0,T]} \|z(t)\|_\VV + C_2,
    \end{aligned}
    \end{equation*}        
    with constants $C_1, C_2> 0$ that neither depend on $\epsilon$ and $\delta$ nor on $\eta$, but only on the coercivity constant $\alpha$, 
    embedding constants, and the data $\ell$ and $\tilde z$.    
    Since this holds for every $t\in [0,T]$, applying Young's inequality implies the assertion for $z_{\epsilon, \delta}$.
    Due to the weak lower semicontinuity of the norm, the inequality readily transfers to the limit by means of \eqref{eq:H1Vconveta}.
\end{proof}

With Lemma~\ref{lem:uniboundeta} at hand, we can now address the uniform convexity of the energy along the solution trajectory, if 
$\eta$ is chosen sufficiently large.

\begin{lemma}\label{lem:Eetaconvex}
	There exists an $\bar \eta < \infty$ such that, for all $\eta \geq \bar\eta$, the penalized energy $\EE_\eta$ satisfies 
	\begin{align}\label{eq:Econvex}
		\mathrm{D}_z^2\EE_\eta(z_{\epsilon, \delta}(t))[v,v]\geq \frac{\alpha}{2}\, \|v\|_\ZZ^2 \quad\text{for all }v\in \ZZ
		\text{, all } t\in [0,T]
        \text{, and all } \epsilon, \delta > 0,
	\end{align}
	where $\alpha > 0$ is the coercivity constant from \eqref{eq:Acoer}.
\end{lemma}

\begin{proof}
    Let $\epsilon, \delta > 0$ and $t\in [0,T]$ be fixed, but arbitrary and abbreviate $r:= \sup_{t\in[0,T]} \|z_{\epsilon, \delta}(t)\|_\ZZ$.
    According to \eqref{eq:F2}, there holds 
    $|\mathrm{D}_z^2\FF(z_{\epsilon, \delta}(t))[v,v]| \leq \gamma(1+r^q) \|v\|_\ZZ \|v\|_\VV$ with $q\in [0,1/2)$
    such that Young's inequality gives for all $v\in \ZZ$ and all $\kappa > 0$ that  
    \begin{equation*}
        |\mathrm{D}_z^2\FF(z_{\epsilon, \delta}(t))[v,v]|\leq \frac{\gamma(1+r^q)\kappa}{2}\|v\|_\ZZ^2+\frac{\gamma(1+r^q)}{2\kappa}\|v\|_\VV^2
    \end{equation*}
    Hence, by choosing $\kappa = \frac{\alpha}{\gamma(1+r^q)}$, it follows from the coercivity of $A$ 
    and Lemma~\ref{lem:uniboundeta} that
    \begin{align*}
	    \mathrm{D}_z^2\EE_\eta(z_{\epsilon, \delta}(t))[v,v]&= \langle A v, v\rangle_{\ZZ^*,\ZZ}+ \mathrm{D}_z^2\FF(z)[v,v]+\eta \|v\|_\VV^2\\
	    &\geq  \frac{\alpha}{2}\|v\|_\ZZ^2+\Big(\eta -\frac{\gamma^2}{2\alpha} (1+r^q)^2\Big)\|v\|_\VV^2 \\
	    & \geq \frac{\alpha}{2}\|v\|_\ZZ^2+\Big(\eta -\frac{\gamma^2}{2\alpha}\, \big[1+(C (\eta + 1))^q\big]^2 \Big)\|v\|_\VV^2.
    \end{align*}
    Because of $q < 1/2$ by assumption, we obtain $\eta -\frac{\gamma^2}{2\alpha}\,[1+(C (\eta + 1))^q]^2 \to \infty$ as $\eta \to \infty$, 
    which implies the existence of $\bar\eta$. Note that $\bar\eta$ does neither depend on $t$ nor $\epsilon$ and $\delta$,
    owing to  Lemma~\ref{lem:uniboundeta}.
\end{proof}

\begin{remark}\label{rem:Eetaconvex}
    The above proof shows that all trajectories $z: [0,T]\to \ZZ$ satisfying \eqref{eq:uniboundeta} with $\eta = \bar\eta$ 
    also satisfy the inequality in \eqref{eq:Econvex}, i.e., 
    $\mathrm{D}_z^2\EE_{\bar\eta}(z(t))[v,v]\geq \frac{\alpha}{2}\, \|v\|_\ZZ^2$ for all $v\in \ZZ$. 
    This observation will be useful in the proof of Lemma~\ref{lem:zepsweak} below.
\end{remark}

For the rest of this section, we fix the penalization parameter to $\eta = \bar\eta$ and assume that $\bar\eta$ is chosen so large
such that \eqref{eq:uniboundeta} also holds for $\tilde z$, i.e., 
\begin{equation}
    \|\tilde z(t)\|_\ZZ \leq C(\bar\eta + 1) \quad \forall\, t\in [0,T],
\end{equation}
which is possible since $\tilde z \in H^1(0,T;\ZZ) \embed C([0,T];\ZZ)$ by assumption.
Lemma~\ref{lem:Eetaconvex} now allows us to derive an $\epsilon$-independent bound on the time derivative. 
For this purpose, we return to the essential estimate \eqref{eq:estE} on the reduced energy 
from Lemma~\ref{lem:z''}, which now reads as follows: 
\begin{equation}\label{eq:D2Eetaest}
\begin{aligned}
	\frac{\epsilon}{2}\, \|z_{\epsilon,\delta}'(t)\|_\VV^2 + \frac{\delta}{2}\, |z_{\epsilon, \delta}'(t)|_{\ZZ}^2
	+ \int_0^t\mathrm{D}_z^2\EE_{\bar\eta}(z_{\epsilon,\delta}(r))[z_{\epsilon,\delta}'(r),z_{\epsilon,\delta}'(r)]\,\d r 
	\qquad\qquad & \\[-1ex]
	\leq \int_0^t \langle \ell'(r)+\bar\eta\, \tilde z'(r),z_{\epsilon,\delta}'(r)\rangle_{\VV^*,\VV}\,\d r. &
\end{aligned}
\end{equation}
Together with Lemma~\ref{lem:Eetaconvex}, this implies
\begin{equation}\label{eq:H1Zbounddel}
	 \|z_{\epsilon,\delta}'\|_{L^2(0,T;\ZZ)} \leq \frac{2}{\alpha}\,\| \ell'+ \bar \eta\, \tilde z'\|_{L^2(0,T;\VV^*)},
\end{equation} 
and, by means of \eqref{eq:H1Vconveta}, this bound also transfers to the limit $z_\epsilon\in H^1(0,T;\ZZ)$. 
Therefore, in combination with Lemma~\ref{lem:uniboundeta}, we obtain that
\begin{equation}\label{eq:H1Zboundeps}
	 \|z_{\epsilon,\delta}\|_{H^1(0,T;\ZZ)} + \|z_{\epsilon}\|_{H^1(0,T;\ZZ)} \leq C(\bar\eta + 1)
	 \quad \forall\, \epsilon, \delta > 0,
\end{equation}
where $C>0$ does neither depend on $\epsilon$ nor on $\delta$. 

\begin{lemma}\label{lem:zepsweak}
	The solution $z_{\epsilon}$ of \eqref{eq:RISepseta} satisfies
	\begin{align}\label{eq:H1Vconvtilde}
		z_\epsilon\rightharpoonup \tilde z \quad \text{in } H^1(0,T;\ZZ) \text{ as } \epsilon \searrow 0.
	\end{align}
\end{lemma}

\begin{proof}
	Due to the boundedness from \eqref{eq:H1Zboundeps},
	 there exists a weakly converging subsequence (also denoted by $z_\epsilon$) and a weak limit $z^*\in H^1(0,T;\ZZ)$ with 
	\begin{equation}\label{eq:H1Vweakconv}
	z_\epsilon\rightharpoonup z^* \text{ in } H^1(0,T;\ZZ),~ \epsilon\searrow 0.
	\end{equation}
    We show that $z^*$ is a differential solution of the penalized system, i.e., it solves
	\begin{align}\label{eq:etaris}
		0\in \partial\RR({z^*}'(t))+\mathrm{D}_z\II_{\bar\eta}(\ell(t) + \bar\eta\,\tilde z(t) ,z^*(t)),\quad z^*(0)=z_0.
	\end{align}
	The arguments are similar to the proof of Theorem \ref{thm:exRISeps}, so we only sketch them briefly.
	By definition of the subdifferential, the solution $z_\epsilon$ also fulfills the energy inequality associated with \eqref{eq:RISepseta2}, 
	which reads
	\begin{align*}
	    \int_0^T \RR(v(r))\,\d r
	    & \geq \int_0^T\RR(z_\epsilon'(r))\,\d r + \int_0^T \epsilon\langle z_\epsilon'(r),z_\epsilon'(r)-v(r)\rangle_{\VV^*,\VV}\,\d r\\
		& \quad + \int_0^T \langle \mathrm{D}_z\II_{\bar\eta}(\ell(r)+\bar\eta\,\tilde z(r),z_\epsilon(r)),z_\epsilon'(r)-v(r)\rangle_{\ZZ^*,\ZZ}\,\d r
	\end{align*}
	for all $v\in L^2(0,T;\ZZ)$.
    In view of \eqref{eq:H1Vweakconv}, the first term on the right hand side can be treated as in \eqref{eq:Repsconv}.
	Concerning the second term, we exploit the boundedness from \eqref{eq:H1Zboundeps} in order to obtain
	\begin{equation*}
		\int_0^T\epsilon\langle z_\epsilon'(r),z_\epsilon'(r)-v(r)\rangle_{\VV^*,\VV}\,\d r\to 0,\quad \epsilon\searrow 0.
	\end{equation*}
	By exploiting the explicit structure of the penalized energy, the last term can be written as 
	\begin{align*}
		&\int_0^T\langle \mathrm{D}_z\II_{\bar\eta}(\ell(r)+\eta\,\tilde z(r),z_\epsilon(r)),z_\epsilon'(r)-v(r)\rangle_{\ZZ^*,\ZZ}\,\d r\\
		& =\int_0^T \langle Az_\epsilon(r),z_\epsilon'(r)-v(r)\rangle_{\ZZ^*,\ZZ}\,\d r+\int_0^T\langle \mathrm{D}_z\FF(z_\epsilon(r)),z_\epsilon'(r)-v(r)\rangle_{\ZZ^*,\ZZ}\,\d r\\
		&\quad -\int_0^T\langle \ell(r)+\bar\eta\, \tilde z(r),z_\epsilon'(r)-v(r)\rangle_{\VV^*,\VV}
		- \bar\eta\, \langle z_\epsilon(r),z_\epsilon'(r)-v(r)\rangle_{\VV^*,\VV}\,\d r.
	\end{align*}
	Analogously to \eqref{eq:DzIepsconv}, Lemma~\ref{lem:limAF} from Appendix~\ref{sec:limAF} applies to the first two terms.
	Further, the weak convergence from \eqref{eq:H1Vweakconv} along with $\ZZ\embed^c \VV$ implies 
	the strong convergence $z_n\to z$ in $L^2(0,T;\VV)$ and thus
	\begin{align*}
		& \lim_{\epsilon\searrow 0} \int_0^T\langle \ell(r) + \bar\eta\, \tilde z(r),z_\epsilon'(r)-v(r)\rangle_{\VV^*,\VV}
		- \bar\eta \langle z_\epsilon(r),z_\epsilon'(r)-v(r)\rangle_{\VV^*,\VV}\,\d r\\
		&\quad =\int_0^T\langle \ell(r)+\bar\eta\, \tilde z(r),{z^*}'(r)-v(r)\rangle_{\VV^*,\VV}
		- \bar\eta \langle z^*(r),{z^*}'(r)-v(r)\rangle_{\VV^*,\VV}\,\d r.
	\end{align*}
	Altogether we have shown that
	\begin{align*}
        0 & \geq 
       \begin{aligned}[t]
	        \liminf_{\epsilon\searrow0} \int_0^T\RR(z_\epsilon'(r)) & - \RR(v(r)) + \epsilon\langle z_\epsilon'(r),z_\epsilon'(r)-v(r)\rangle_{\VV^*,\VV}\\
	        	&+\langle \mathrm{D}_z\II_{\bar\eta}(\ell(r)+\bar\eta\,\tilde z(r),z_\epsilon(r)),z_\epsilon'(r)-v(r)\rangle_{\ZZ^*,\ZZ}\,\d r
	    \end{aligned}\\
        &\geq \int_0^T\RR({z^*}'(r))-\RR(v(r)) 
        + \langle \mathrm{D}_z\II_{\bar\eta}(\ell(r) + \bar\eta\, \tilde z(r),z^*(r)),{z^*}'(r)-v(r)\rangle_{\ZZ^*,\ZZ}\,\d r
	\end{align*}
	for all $v\in L^2(0,T;\ZZ)$. By standard arguments using the fundamental lemma of the calculus of variations, 
	this energy inequality is equivalent to \eqref{eq:etaris} so that 
     $z^*$ is indeed a differential solution of \eqref{eq:etaris}.
          
	Finally, we show that $z^*=\tilde z$. Thanks to construction of the penalized energy, we have
	$\mathrm{D}_z\II_{\bar\eta}(\ell(t)+\bar\eta\,\tilde z(t),\tilde z(t))=\mathrm{D}_z\II(\ell(t),\tilde z(t))$ and thus, 
	in light of \eqref{eq:RIStilde}, $\tilde z$ is a differential solution of \eqref{eq:etaris}, too. 
	Using the uniform convexity of the penalized energy, one can then show that $z^*=\tilde z$, which concludes the proof.
	Since the arguments are classical, see, e.g.,  \cite[Section~3.4.4]{MR15}, we postpone the corresponding proof to Appendix~\ref{sec:diffuni}. 
\end{proof}

Before we are able to show that even $z_\epsilon\to \tilde z$ strongly in $H^1(0,T;\ZZ)$, 
we need a further property of $\tilde z$, which is as follows:

\begin{lemma}
	The differential solution $\tilde z$ of \eqref{eq:RIStilde} satisfies
	\begin{align}\label{eq:D2Etilde}
		\mathrm{D}_z^2\EE_{\bar\eta}(\tilde z(t))[\tilde z'(t),\tilde z'(t)]
		-\langle\ell'(t) + \bar\eta\, \tilde z'(t),\tilde z'(t)\rangle_{\VV^*,\VV}=0
		\quad \text{ f.a.a. } t\in(0,T).
	\end{align}
\end{lemma}

\begin{proof}
    As seen at the end of the previous proof, $\tilde z$ is also a differential solution of 
    \eqref{eq:etaris}, i.e., it holds that
    \begin{equation*}
        0 \in \partial\RR(\tilde z'(t)) + \Drm_z \II_{\bar\eta}(\ell(t) + \bar\eta\,\tilde z(t), \tilde z(t)) 
        \quad \text{a.e.\ in } (0,T).
    \end{equation*}
	Since $\partial\RR(v)\subset\partial\RR(0)$ for all $v\in\ZZ$ and $\tilde z$ is continuous, we have 
    $-\Drm_z \II_{\bar\eta}(\ell(t) + \bar\eta\,\tilde z(t), \tilde z(t)) \in\partial\RR(0)$ for all $t\in[0,T]$. 
    The characterization of the subdifferential of one-homogeneous functionals then leads to
	\begin{alignat}{3}
		\RR(\tilde z'(t)) &=\langle -\Drm_z \II_{\bar\eta}(\ell(t) + \bar\eta\,\tilde z(t), \tilde z(t)) ,\tilde z'(t)\rangle_{\ZZ^*,\ZZ} 
		& \quad & \text{f.a.a. } t\in (0,T),\label{eq:subtilde1}\\
		\RR(v) &\geq \langle -\Drm_z \II_{\bar\eta}(\ell(t) + \bar\eta\,\tilde z(t), \tilde z(t)) ,v\rangle_{\ZZ^*,\ZZ} 
		& \quad & \text{for all } v\in\ZZ, \tau\in[0,T].\label{eq:subtilde2}
	\end{alignat}
	Now let $t\in(0,T)$ be a Lebesgue point of $\tilde z'$ such that \eqref{eq:subtilde1} is fulfilled, 
	which is true for almost all $t\in (0,T)$.  If we set $v=\tilde z'(t)$ and $\tau=t\pm h$
	with $h>0$ in \eqref{eq:subtilde2} and subtract \eqref{eq:subtilde1} from the arising inequality, then
	\begin{align*}
		\langle \mathrm{D}_z\II_{\bar\eta}(\ell(t\pm h) + \bar\eta\, \tilde z(t \pm h),\tilde z(t\pm h))
		-\mathrm{D}_z\II_{\bar\eta}(\ell(t) + \bar\eta\, \tilde z(t),\tilde z(t)),
		\tilde z'(t)\rangle_{\ZZ^*,\ZZ}\geq 0.
	\end{align*}
	is obtained. Exploiting the explicit structure of the energy $\II$ and rearranging the terms leads to
	\begin{multline*}
		\langle \mathrm{D}_z\EE_{\bar\eta}(\tilde z(t\pm h))-\mathrm{D}_z\EE_{\bar\eta}(\tilde z(t)),\tilde z'(t)\rangle_{\ZZ^*,\ZZ} \\
		-\langle \ell(t\pm h)-\ell(t) + \bar\eta (\tilde z(t\pm h)- \tilde z(t)),\tilde z'(t)\rangle_{\VV^*,\VV}\geq 0.
	\end{multline*}
	Finally, dividing by $\pm h$, passing $h\searrow 0$, and using Lebesgue's differentiation theorem 
	give \eqref{eq:D2Etilde} in all points $t\in(0,T)$, which are additionally Lebesgue points of $\ell'$. 
	Since this holds for almost all points $t\in(0,T)$, the proof is completed.
\end{proof}

We are now in the position to prove the strong convergence of $z_\epsilon$ to $\tilde z$ in $H^1(0,T;\ZZ)$. 
This is the point where the intermediate space $\WW$ and the additional assumption on the nonlinearity $\FF$ 
from \eqref{eq:D2FWcont} comes into play.

\begin{lemma}\label{lem:strongconvH1}
    It holds that 
	\begin{align}\label{eq:stH1Vconvtilde}
		z_\epsilon\to \tilde z \text{ in } H^1(0,T;\ZZ), \text{ as } \epsilon \searrow 0.
	\end{align}
\end{lemma}

\begin{proof}
Integrating \eqref{eq:D2Etilde} over $(0,T)$ and subtracting the resulting equation from \eqref{eq:D2Eetaest} results in
\begin{align*}
	0 & \geq 
    \begin{aligned}[t]
        	\int_0^T  \Big( & \mathrm{D}_z^2\EE_{\bar\eta}(z_{\epsilon,\delta}(r))[z_{\epsilon,\delta}'(r),z_{\epsilon,\delta}'(r)] 
	    -\mathrm{D}_z^2\EE_{\bar\eta}(\tilde z(r))[\tilde z'(r),\tilde z '(r)]\\
	     & -\langle \ell'(r)+\bar\eta\, \tilde z'(r),z_{\epsilon,\delta}'(r)-\tilde z'(r)\rangle_{\VV^*,\VV}\Big) \d r    
    \end{aligned}\\
    &= \int_0^T \mathrm{D}_z^2\EE_{\bar\eta}(z_{\epsilon,\delta}(r))[z_{\epsilon,\delta}'(r)-\tilde z'(r),z_{\epsilon,\delta}'(r)-\tilde z'(r)]\,\d r\\
    &\quad - \int_0^T\langle \ell'(r)+\bar\eta\, \tilde z'(r),z_{\epsilon,\delta}'(r)-\tilde z'(r)\rangle_{\VV^*,\VV}\,\d r\\
    &\quad + 
    \begin{aligned}[t]
        \int_0^T \Big( & 2\,\mathrm{D}_z^2\EE_{\bar\eta}(z_{\epsilon,\delta}(r))[\tilde z'(r),z_{\epsilon,\delta}'(r)]\\
        & - \mathrm{D}_z^2\EE_{\bar\eta}(\tilde z(r))[\tilde z'(r),\tilde z '(r)]
        -\mathrm{D}_z^2\EE_{\bar\eta}(z_{\epsilon,\delta}(r))[\tilde z'(r),\tilde z '(r)]\Big) \d r
    \end{aligned}\\
    &\geq  \frac{\alpha}{2} \int_0^T \|z_{\epsilon,\delta}'(r)-\tilde z'(r)\|_\ZZ^2\,\d r
    - \int_0^T\langle \ell'(r)+\bar\eta\, \tilde z'(r),z_{\epsilon,\delta}'(r)-\tilde z'(r)\rangle_{\VV^*,\VV}\,\d r\\
    &\quad + 
    \begin{aligned}[t]
        \int_0^T \Big( & 2\,\mathrm{D}_z^2\EE_{\bar\eta}(z_{\epsilon,\delta}(r))[\tilde z'(r),z_{\epsilon,\delta}'(r)]\\
        & - \mathrm{D}_z^2\EE_{\bar\eta}(\tilde z(r))[\tilde z'(r),\tilde z '(r)]
        -\mathrm{D}_z^2\EE_{\bar\eta}(z_{\epsilon,\delta}(r))[\tilde z'(r),\tilde z '(r)]\Big) \d r,
    \end{aligned}        
\end{align*}
where we used the uniform convexity of $\EE_{\bar\eta}$ along the solution trajectory from \eqref{eq:Econvex}.
Next we pass to the limit $\delta\searrow 0$ in the terms on the right hand side.
First of all, the weak lower semicontinuity of the norm together with the weak convergence from \eqref{eq:H1Vconveta} yields
\begin{align*}
	\frac{\alpha}{2} \int_0^T \|z_{\epsilon}'(r)-\tilde z'(r)\|_\ZZ^2\,\d r
	\leq \liminf_{\delta\searrow 0} \frac{\alpha}{2} \int_0^T \|z_{\epsilon,\delta}'(r)-\tilde z'(r)\|_\ZZ^2\,\d r.
\end{align*}
Another use of \eqref{eq:H1Vconveta} gives the convergence of the load term to the limit 
$\int_0^T\langle \ell'(r)+ \bar\eta\, \tilde z'(r),z_{\epsilon}'(r)-\tilde z'(r)\rangle_{\VV^*,\VV}\,\d r$ as $\delta \searrow 0$.
To treat the third term, we show that 
\begin{align}\label{eq:D2EL2conv}
	\mathrm{D}_z^2\EE_{\bar\eta}(z_{\epsilon,\delta})\tilde z'
	\to \mathrm{D}_z^2\EE_{\bar\eta}(z_{\epsilon})\tilde z' \text{ in } L^2(0,T;\VV^*),~ \delta\searrow 0.
\end{align}
Due to 
\begin{equation*}
    \mathrm{D}_z^2\EE_{\bar\eta}(z)[v, w]
    =  \dual{Av}{w}_{\ZZ^*, \ZZ} + \Drm_z^2 \FF(z)[v,w] + \bar\eta\, \dual{v}{w}_{\VV^*, \VV}
    \quad \forall\, z, v, w\in \ZZ,
\end{equation*}
there holds that 
\begin{equation*}
    \mathrm{D}_z^2\EE_{\bar\eta}(z_{\epsilon,\delta})\tilde z'
	- \mathrm{D}_z^2\EE_{\bar\eta}(z_{\epsilon})\tilde z'
	= \mathrm{D}_z^2\FF(z_{\epsilon,\delta})\tilde z'
	- \mathrm{D}_z^2\FF(z_{\epsilon})\tilde z' .
\end{equation*}
Because of $\ZZ\hookrightarrow^c\WW$, the pointwise weak convergence \eqref{eq:H1Vconveta} yields
$z_{\epsilon,\delta}(t)\to z_\epsilon(t)$ in $\WW$ such that the continuity of $\Drm_z^2 \FF$ from \eqref{eq:D2FWcont} implies 
\begin{equation*}
    \mathrm{D}_z^2\EE_{\bar\eta}(z_{\epsilon,\delta}(t))\tilde z'(t)\to \mathrm{D}_z^2\EE_{\bar\eta}(z_{\epsilon}(t))\tilde z' (t)
    \quad \text{in $\VV^*$ f. a.a.\ $t\in (0,T)$.}
\end{equation*}
Moreover, by exploiting \eqref{eq:F2} in combination with the uniform bound from \eqref{eq:uniboundeta} gives 
 $\|\mathrm{D}_z^2\FF(z_{\epsilon,\delta}(t))\tilde z'(t)\|_{\VV^*}\leq \gamma(1 + (C(\bar\eta + 1))^q)\,\|\tilde z'(t)\|_\ZZ$
and so, \eqref{eq:D2EL2conv} follows from Lebesgue's dominated convergence theorem. 
Together with the weak convergence in $H^1(0,T;\ZZ)$ from \eqref{eq:H1Vconveta} this allows to pass to the limit in 
the third term.
All in all, by passing to the limit $\delta\searrow 0$, we thus arrive at
\begin{equation}\label{eq:zepsztildeH1conv}
\begin{aligned}
    & \frac{\alpha}{2} \int_0^T \|z_{\epsilon}'(r)-\tilde z'(r)\|_\ZZ^2\,\d r \\
    & \quad \leq \int_0^T\langle \ell'(r)+\bar\eta\, \tilde z'(r),z_{\epsilon}'(r)-\tilde z'(r)\rangle_{\VV^*,\VV}\,\d r \\
    & \qquad + 
    \begin{aligned}[t]
        \int_0^T \Big( & 2\,\mathrm{D}_z^2\EE_{\bar\eta}(z_{\epsilon}(r))[\tilde z'(r),z_{\epsilon}'(r)] \\[-1ex]
	    & - \big(\mathrm{D}_z^2\EE_{\bar\eta}(\tilde z(r)) + \mathrm{D}_z^2\EE_{\bar\eta}(z_{\epsilon}(r))\big)
	    [\tilde z'(r),\tilde z '(r)]\Big) \d r .
    \end{aligned}
\end{aligned}
\end{equation}
Now we repeat exactly the same arguments for the above inequality using the 
weak convergence of $(z_\epsilon)_{\epsilon}$ to $\tilde z$ in $H^1(0,T;\ZZ)$ by \eqref{eq:H1Vconvtilde}.
This implies that the right hand side in \eqref{eq:zepsztildeH1conv} converges to zero as $\epsilon \searrow 0$, 
which ultimately implies $z_\epsilon'\to\tilde z'$ in $L^2(0,T;\ZZ)$. Since $z_\epsilon(0)=\tilde z(0)=z_0$, 
this finally gives the claim.
\end{proof}

Before we formulate the main result of this section, the reverse approximation property, we prove the following lemma, whose statement deals with an end time property for the solutions of the penalized viscous regularized systems. This property 
ensures feasibility of the recovery sequence for the regularized optimal control problems, which will be part of the next section.

\begin{lemma}\label{lem:dist<eps}
	There exists a constant $C>0$ independent of $\epsilon$ 
	such that the solution $z_\epsilon$ of \eqref{eq:RISepseta} satisfies
	\begin{equation*}
		\dist_{\VV^*}(-\mathrm{D}_z\JJ_{\bar\eta}(T, \ell(T), z_\epsilon(T)),\partial\RR(0))\leq C\,\sqrt{\epsilon}.
	\end{equation*}
\end{lemma}

\begin{proof}
    We return to \eqref{eq:D2Eetaest}, which, thanks to \eqref{eq:Econvex}, implies that
    \begin{equation*}
    \begin{aligned}
        	\frac{\epsilon}{2}\|z_{\epsilon,\delta}'(t)\|_\VV^2 
        	&\leq \int_0^t \langle \ell'(r)+\bar\eta\, \tilde z'(r),z_{\epsilon,\delta}'(r)\rangle_{\VV^*,\VV}\,\d r\\
		&\leq \| \ell'+\bar\eta\, \tilde z'\|_{L^2(0,T;\VV^*)}\,\|z_{\epsilon,\delta}'\|_{L^2(0,T;\VV)}
        \leq \frac{2}{\alpha}\,\| \ell'+\bar\eta\, \tilde z'\|_{L^2(0,T;\VV^*)}^2,    
    \end{aligned}
    \end{equation*}
    where we used \eqref{eq:H1Zbounddel} in the last step. 
    By the weak convergence from \eqref{eq:H1Vconveta}, this bound directly carries over to the weak limit, i.e.,    
    \begin{equation}\label{eq:sqrteps}
        	 \|z_{\epsilon}'\|_{L^\infty(0,T;\VV)}\leq \frac{C}{\sqrt{\epsilon}}
    \end{equation}
    with $C :=  \frac{2}{\sqrt{\alpha}}\, \| \ell'+\bar\eta\, \tilde z'\|_{L^2(0,T;\VV^*)}$.
    As $z_\epsilon$ is the solution of \eqref{eq:RISepseta}, it holds for almost every $t\in (0,T)$ that
    \begin{equation*}
        -\mathrm{D}_z\JJ_{\bar\eta}(t, \ell(t), z_\epsilon(t)) - \epsilon\, z_\epsilon'(t) 
        \in \partial\RR(z_\epsilon '(t)) 
        \subset \partial\RR(0) 
    \end{equation*}
    and therefore, in combination with \eqref{eq:sqrteps},
    \begin{equation}\label{eq:distae}
        \dist_{\VV^*} (-\mathrm{D}_z\JJ_{\bar\eta}(t, \ell(t), z_\epsilon(t)), \partial\RR(0)) 
        \leq \epsilon \, \|z_\epsilon'(t)\|_{\VV}
        \leq \sqrt{\epsilon}\,C
        \quad \text{a.e.\ in }(0,T).
    \end{equation}        
	To prove the estimate for $t = T$, we argue analogously to the proof of Corollary~\ref{cor:DzI}. 
	Since $z_\epsilon$ and $\ell$ are continuous in time with values in $\ZZ$ and $\VV^*$, respectively, 
	$\mathrm{D}_z\JJ_{\bar\eta}(\,\cdot\,, \ell(\cdot), z_\epsilon(\cdot))$ is continuous in time with values in 
	$\ZZ^*$. Moreover, we can approximate $T$ with time points $t_n$ where 
	\eqref{eq:distae} is valid and thus,  
	$\mathrm{D}_z\JJ_{\bar\eta}(t_n, \ell(t_n), z_\epsilon(t_n))$ is bounded in $\VV^*$ (by the boundedness of $\partial\RR(0)$)
	and therefore converges weakly in $\VV^*$. 
	By continuity with values in $\ZZ^*$, this weak limit must be 
	$\mathrm{D}_z\JJ_{\bar\eta}(T, \ell(T), z_\epsilon(T))$.
    Thus, the weak lower semicontinuity of the distance implies that  \eqref{eq:distae} also holds in $T$.
\end{proof}

With this at hand, we are now able to prove the main result of this section:

\begin{theorem}[Reverse approximation property]\label{thm:revapp}
    Under Assumption \ref{ass:diffsol}, there exists a sequence 
	$(z_\epsilon,\ell_\epsilon)_{\epsilon>0}\subset H^1(0,T;\ZZ)\times H^1(0,T;\VV^*)$ 
	such that $z_\epsilon$ is a solution of the $\epsilon$-viscous regularization \eqref{eq:RISeps} 
	with external load $\ell_\epsilon$, i.e.,
	\begin{align}\label{eq:riselleps}
		0\in\partial\RR_\epsilon(z_\epsilon'(t))+\mathrm{D}_z\II(\ell_\epsilon(t),z_\epsilon(t)) \text{ f.a.a. } t\in(0,T),
		\quad z_\epsilon(0)=z_0
	\end{align}
	and there holds that
	\begin{align}\label{eq:convzell}
		z_\epsilon\to \tilde z \text{ in } H^1(0,T;\ZZ),\quad \ell_\epsilon\to \ell \text{ in } H^1(0,T;\VV^*), \text{ as } \epsilon\searrow 0,
	\end{align}
	where $\tilde z$ is the reparametrized function from \eqref{eq:tildezdef}.
	Moreover, for all $\epsilon > 0$, we have
	\begin{align}\label{eq:enddist}
	 	\ell_\epsilon(0)=\ell(0)
	 	\quad \text{and} \quad
	 	\dist_{\VV^*}(-\mathrm{D}_z\II(\ell_\epsilon(T),z_\epsilon(T)),\partial \RR(0))\leq C\, \sqrt{\epsilon}
	\end{align}
	with $C>0$ independent of $\epsilon$.
\end{theorem}

\begin{proof}
	If we set $\ell_\epsilon:= \ell-\bar\eta(z_\epsilon-\tilde z)$, then we obtain
	\begin{align*}
		\mathrm{D}_z\II(\ell_\epsilon(t),z_\epsilon(t))
		&=Az_\epsilon(t)+\mathrm{D}_z\FF(z_\epsilon(t)) - \ell(t) + \bar\eta(z_\epsilon(t) - \tilde z(t))\\
		&=\mathrm{D}_z\JJ_{\bar\eta}(t, \ell(t),z_\epsilon(t))
	\end{align*}
	such that \eqref{eq:riselleps} equals \eqref{eq:RISepseta} (with $\eta = \bar\eta$). 
	Thus, Lemma \ref{lem:strongconvH1} yields $z_\epsilon\to \tilde z $ in $H^1(0,T;\ZZ)$ 
	such that $ \ell_\epsilon\to \ell$ in  $H^1(0,T;\VV^*)$. 
	Eventually, exploiting $z_\epsilon(0)=\tilde z(0)=z_0$ and Lemma \ref{lem:dist<eps} gives \eqref{eq:enddist}, 
	which completes the proof.
\end{proof}


\section{Approximation by Regularized Optimal Control Problems}\label{sec:ocpapprox}

Using the reverse approximation result from Theorem~\ref{thm:revapp}, we can now prove the two main 
results of this work, namely the convergence of (subsequences of) optimal solutions to \eqref{eq:optcontroleps} 
to an optimal solution of \eqref{eq:optcontrol} in Theorem~\ref{thm:ocpapproxeps} 
and a similar result for the double viscous regularization in Theorem~\ref{thm:ocpapproxepsdelta}.

\begin{theorem}[Convergence of \eqref{eq:optcontroleps}]\label{thm:ocpapproxeps}
    Suppose that Assumption~\ref{ass:diffsol} is fulfilled and let 
    $(z_\epsilon^*, \ell_\epsilon^*)_{\epsilon>0}\subset H^1(0,T;\ZZ) \times H^1(0,T;\VV^*)$ 
    be a sequence of global minimizers of \eqref{eq:optcontroleps}.
    Then there exists a subsequence (denoted by the same symbol) such that 
    the associated parametrized solution $(S_\epsilon^*,\hat t_\epsilon^*,\hat z_\epsilon^*)$ 
    as defined in \eqref{eq:defs}--\eqref{eq:defhatz} satisfies 
    	\begin{gather}
		S_\epsilon^*\to S^* \text{ in } \R,\quad \hat t_\epsilon^*\rightharpoonup^*\hat t^* \text{ in }W^{1,\infty}(0,S^*;\R), \\
		\ell_\epsilon^*\to \ell^* \text{ in } H^1(0,T;\VV^*), \label{eq:convloadsstrong}\\
		\hat z_\epsilon^*\rightharpoonup^* \hat z^* \text{ in } L^\infty(0,S^*;\ZZ), 
		\quad \hat z_\epsilon^*\to \hat z^* \text{ in } C([0,S^*];\VV).
	\end{gather}
	Moreover, every such accumulation point satisfies 
    \begin{equation*}
        (S^*, \hat t^*, \hat z^*, \ell^*) \in [T,\infty)\times W^{1,\infty}(0,T)\times 
        \mathrm{AC}^\infty([0,S];\RR)\cap L^\infty (0,S;\ZZ) \times H^1(0,T;\VV^*)
    \end{equation*}    	
	and is a global minimizer of \eqref{eq:optcontrol}. 
\end{theorem}

\begin{proof}
    \emph{(i) Existence of a converging subsequence}\\
    Consider the optimal solution $(S, \hat t, \hat z, \ell)$ from Assumption~\ref{ass:diffsol}. 
    By Theorem~\ref{thm:revapp}, we know that there is a sequence $(z_\epsilon, \ell_\epsilon)$ such that 
    $\ell_\epsilon \to \ell$ in $H^1(0,T;\VV^*)$ and 
 	$z_\epsilon(T)\to \tilde z(T) = \hat z(\hat t^{-1}(T))=\hat z(S)$ in $\VV$.
    The continuity of $j$ thus implies    
    	\begin{equation}\label{eq:recovery}
    	\begin{aligned}
		\lim_{\epsilon\searrow 0} J_\epsilon(z_\epsilon, \ell_\epsilon)
		&=\lim_{\epsilon\searrow 0} j(z_\epsilon(T))+\frac{\beta}{2}\,\|\ell_\epsilon\|_{H^1(0,T;\VV^*)}^2\\
		&=j(\hat z(S))+\frac{\beta}{2}\,\|\ell\|_{H^1(0,T;\VV^*)}^2 = J(S,\hat z,\ell).    	
    	\end{aligned}
	\end{equation}
   Moreover, $(z_\epsilon, \ell_\epsilon)$ is feasible for \eqref{eq:optcontroleps}, as we will see next.    
    First of all, by Theorem~\ref{thm:revapp} it solves \eqref{eq:RISeps}. Furthermore,  by \eqref{eq:enddist}, 
    there holds 
    \begin{equation*}
        \dist_{\VV^*}(- \Drm_z\II(\ell_\epsilon(0), z_0), \partial\RR(0)) 
        = \dist_{\VV^*}(- \Drm_z\II(\ell(0), z_0), \partial\RR(0)) = 0 \quad \forall\, \epsilon > 0
    \end{equation*}
    due to the feasibility of $\ell$ for \eqref{eq:optcontrol}. In addition, the second statement in \eqref{eq:enddist}
    implies that 
    \begin{equation}\label{eq:distZstarendtime}
    \begin{aligned}
        & \dist_{\ZZ^*}(- \Drm_z\II(\ell_\epsilon(T), z_\epsilon(T)), \partial\RR(0)) \\
        & \qquad \leq C\, \dist_{\VV^*}(- \Drm_z\II(\ell_\epsilon(T), z_\epsilon(T)), \partial\RR(0)) 
        \leq \epsilon^{\frac{1}{4}}
    \end{aligned}
    \end{equation}
    for all $\epsilon>0$ sufficiently small. 
    Consequently, for $\epsilon>0$ small enough, $(z_\epsilon, \ell_\epsilon)$ is feasible for \eqref{eq:optcontroleps}.

    Therefore, the optimality of $(z_\epsilon^*, \ell_\epsilon^*)$ for the latter implies
    \begin{equation}\label{eq:zepsopt}
        j(z_\epsilon^*(T))+\frac{\beta}{2}\,\|\ell_\epsilon^*\|_{H^1(0,T;\VV^*)}^2
        = J_\epsilon(z_\epsilon^*, \ell_\epsilon^*) \leq J_\epsilon(z_\epsilon, \ell_\epsilon)
    \end{equation}
    for all $\epsilon > 0$ sufficiently small. Since the right hand side converges by \eqref{eq:recovery} and 
    is thus bounded and $j$ is bounded from below
    by assumption, this implies that $(\ell_\epsilon^*)_\epsilon$ is bounded in $H^1(0,T;\VV^*)$. 
	Thus, there exists a weakly converging (not relabeled) subsequence 
	such that $\ell_\epsilon^*\rightharpoonup\ell^*$ in $H^1(0,T;\VV^*)$.
    
    \emph{(ii) Feasibility of the limit}\\
    Thanks to Theorem~\ref{thm:exparasol}, there is yet another subsequence (again denoted by the same symbol)
    such that the reparametrized solution $(S_\epsilon^*, \hat t_\epsilon^*, \hat z_\epsilon^*)$ 
    according to \eqref{eq:defs}--\eqref{eq:defhatz} (with $z_\epsilon = z_\epsilon^*$) 
    converges in the sense of \eqref{eq:convparat}--\eqref{eq:convparazs}
    to a paramterized BV solution associated with $\ell^*$ denoted by $(S^*, \hat t^*, \hat z^*)$, i.e., 
    \begin{equation}\label{eq:limitparasol}
        (S^*, \hat t^*, \hat z^*) \in \LL(z_0, \ell^*)
    \end{equation}
    with $\LL(z_0, \ell^*)$ as defined in Section~\ref{sec:optcontrol}. The weak convergence of $\ell^*_\epsilon$ 
    in $H^1(0,T;\VV^*)$ moreover implies $\ell_\epsilon^*(0)\rightharpoonup\ell^*(0)$ in $\VV^*$ and thus 
    \begin{equation}\label{eq:limitlocstab0}
        -\mathrm{D}_z\II(\ell^*(0),z_0)\in\partial \RR(0)
    \end{equation}        
    by the weak closedness of $\partial\RR(0)$. For the end time condition, 
    we use again the pointwise convergence of the loads, i.e., $\ell_\epsilon^*(T)\rightharpoonup\ell^*(T)$ in $\VV^*$, 
    and $\hat z_\epsilon^*(S_\epsilon^*)\rightharpoonup \hat z^*(S^*) \text{ in } \ZZ$ 
    from \eqref{eq:convparat} and \eqref{eq:convparazs} in 
	order to infer from \eqref{DFweak} that 
    \begin{equation}\label{eq:DzIconvZstar}
        - \Drm_z \II(\ell_\epsilon^*(T), z_\epsilon^*(T))
        =  - \Drm_z\II(\ell_\epsilon^*(T), \hat z_\epsilon^*(S_\epsilon))
        \weak - \Drm_z\II(\ell^*(T), \hat z^*(S)) \quad \text{in } \ZZ^*.
    \end{equation}    	
	The weak lower semicontinuity of the distance then implies
	\begin{align*}
        & \dist_{\ZZ^*}(-\mathrm{D}_z\II(\ell^*(T),\hat z^*(S^*)),\partial \RR(0)) \\
		& \qquad 
		\leq \liminf_{\epsilon\searrow 0}
		\dist_{\ZZ^*}(-\mathrm{D}_z\II(\ell_\epsilon^*(T),z_\epsilon^*(T)),\partial \RR(0))
		\leq\liminf_{\epsilon\searrow 0} \epsilon^{\frac{1}{4}}=0.
	\end{align*}
	Together with \eqref{eq:limitparasol} and \eqref{eq:limitlocstab0}, 
	this shows that the limit $(S^*,\hat t^*,\hat z^*,\ell^*)$ is feasible for \eqref{eq:optcontrol}.     
    
    \emph{(iii) Optimality of the limit}\\
    The optimality of the limit immediately follows by means the recovery sequence from part (i).
    Indeed, as seen above, there holds $z_\epsilon^*(T) \weakly \hat z^*(S^*)$ in $\ZZ$ and by the 
    compact embedding $\ZZ \embed^c \VV$, $z_\epsilon^*(T)$ converges strongly in $\VV$ to $\hat z^*(S^*)$. 
    Together with the continuity of $j$ and the weak lower semicontinuity of the squared norm, this implies
    \begin{equation}\label{eq:convobj}
    \begin{aligned}
        J(S^*, \hat z^*, \ell^*) 
        & \leq \liminf_{\epsilon \searrow 0} J_\epsilon(z_\epsilon^*, \ell_\epsilon^*) \\
        &\leq \limsup_{\epsilon \searrow 0} J_\epsilon(z_\epsilon^*, \ell_\epsilon^*)
        \leq \limsup_{\epsilon \searrow 0}  J_\epsilon(z_\epsilon, \ell_\epsilon)
        = J(S, \hat z, \ell),    
    \end{aligned}
    \end{equation}
    where we used \eqref{eq:zepsopt} and \eqref{eq:recovery} for the last inequality and equality, respectively.
    Since $(S, \hat t, \hat z, \ell)$ is a minimizer of \eqref{eq:optcontrol} by assumption and 
    $(S^*, \hat t^*, \hat z^*, \ell^*)$ is feasible for \eqref{eq:optcontrol} by part (2), this shows its optimality.
    
    It remains to verify the strong convergence of the loads in \eqref{eq:convloadsstrong}.
    To this end, observe that \eqref{eq:convobj} implies 
    \begin{equation*}
    \begin{aligned}
        j(\hat z^*(S^*)) +\frac{\beta}{2}\,\|\ell^*\|_{H^1(0,T;\VV^*)}^2 = J(S^*, \hat z^*, \ell^*) 
        &= \lim_{\epsilon \searrow 0} J_\epsilon(z_\epsilon^*, \ell_\epsilon^*) \\
        &= \lim_{\epsilon \searrow 0} j(z^*_\epsilon(T)) +\frac{\beta}{2}\,\|\ell_\epsilon^*\|_{H^1(0,T;\VV^*)}^2. 
    \end{aligned}
    \end{equation*}
    Since $j(z^*_\epsilon(T))\to j(\hat z^*(S^*))$ as seen above, this implies 
    $\|\ell_\epsilon^*\|_{H^1(0,T;\VV^*)} \to \|\ell^*\|_{H^1(0,T;\VV^*)}$, which, together with the weak convergence, 
    implies strong convergence of the loads.
\end{proof}

\begin{remark}
    Similarly to Remark~\ref{rem:ocpexistVstar}, we observe that the above proof also works, if 
    the end time constraint is considered with the $\VV^*$-distance instead of the $\ZZ^*$-distance.
    Analogously to \eqref{eq:distZstarendtime}, the second statement in \eqref{eq:enddist} shows that 
    the recovery sequence satisfies the end time constraint.
    Moreover, the feasibility of the limit of the sequence of minimizers follows as above by noting that 
    the $\ZZ^*$-distance can be estimated from above by the $\VV^*$-distance.
    The benefit of using the $\ZZ^*$-distance will become clear in the next proof, see \eqref{eq:distZstarineq} below.
\end{remark}

Next, we investigate the convergence properties of \eqref{eq:optcontrolepsdelta} for fixed $\epsilon > 0$ 
and $\delta$ tending to zero.

\begin{proposition}[Convergence of \eqref{eq:optcontrolepsdelta}]\label{prop:ocpconvdelta}
    Let $\epsilon > 0$ be fixed and suppose that \eqref{eq:optcontroleps} 
    admits a solution $(\bar z_\epsilon, \bar\ell_\epsilon) \in H^1(0,T;\ZZ)\times H^1(0,T;\VV^*)$.
    Let $(z^*_{\epsilon, \delta}, \ell^*_{\epsilon, \delta})_{\delta > 0}$ be a sequence of optimal solutions of 
    \eqref{eq:optcontrolepsdelta}. Then there exists a subsequence (denoted by the same symbol) such that 
    \begin{equation}\label{eq:ocpconvdelta}
        z_{\epsilon, \delta}^* \weak z_{\epsilon}^* \quad \text{in } H^1(0,T;\ZZ)
        \quad \text{and} \quad 
        \ell_{\epsilon, \delta}^* \to \ell_{\epsilon}^* \quad \text{in } H^1(0,T;\VV^*)
    \end{equation}
    and the limit $(z_\epsilon^*, \ell_\epsilon^*)$ is a global minimizer of \eqref{eq:optcontroleps}.
\end{proposition}

\begin{proof}
    The proof mainly follows the lines of the one of Theorem~\ref{thm:ocpapproxeps}, 
    with the only difference that the construction of the recovery sequence is now elementary.    
    We simply choose the sequence $(\bar z_{\epsilon, \delta}, \bar \ell_\epsilon)_{\delta > 0}$, 
    where $\bar z_{\epsilon, \delta}$ denotes the solution to \eqref{eq:RISepsdel} with $\ell = \bar\ell_\epsilon$. 
    To check the feasibility of $(\bar z_{\epsilon, \delta}, \bar \ell_\epsilon)$, first note that 
    the initial time condition $- \Drm_z \II(\ell_\epsilon, z_0) \in \partial\RR(0)$ is fulfilled, as $\bar\ell_\epsilon$ 
    is feasible for \eqref{eq:optcontroleps}. 
    Moreover, by construction, $\bar z_{\epsilon, \delta}$ fulfills the ODE in \eqref{eq:optcontrolepsdelta}. 
    It remains to verify the end time constraint. Since $(\bar z_\epsilon, \bar\ell_\epsilon)$ satisfy the 
    end time constraint of \eqref{eq:optcontroleps}, we obtain the existance of a constant $C_\epsilon$, 
    depending on $\epsilon$, but independent of $\delta$ such that 
    \begin{equation}\label{eq:distZstarineq}
    \begin{aligned}
        & \dist_{\ZZ^*}(-\Drm_z\II(\bar\ell_\epsilon(T), \bar z_{\epsilon, \delta}(T)), \partial\RR(0)) \\
        & \quad \leq \dist_{\ZZ^*}(-\Drm_z\II(\bar\ell_\epsilon(T), \bar z_{\epsilon}(T)), \partial\RR(0)) \\
        & \qquad + \|\Drm_z\II(\bar \ell_\epsilon(T), \bar z_\epsilon(T)) 
        - \Drm_z\II(\bar \ell_\epsilon(T), \bar z_{\epsilon, \delta}(T))\|_{\ZZ^*}\\
        & \quad \leq \epsilon^{\frac{1}{4}} + \|A(\bar z_{\epsilon, \delta}(T) - \bar z_\epsilon(T))\|_{\ZZ^*}
        + \|\Drm_z \FF(\bar z_{\epsilon, \delta}(T)) - \Drm_z \FF(\bar z_{\epsilon}(T)) \|_{\ZZ^*}\\
        & \quad \leq \epsilon^{\frac{1}{4}} + C_\epsilon\,\|\bar z_{\epsilon, \delta}(T) - \bar z_\epsilon(T)\|_{\ZZ}
        \leq \epsilon^{\frac{1}{4}} + C_\epsilon \, \sqrt{\delta},
    \end{aligned}
    \end{equation}
    where we used Proposition~\ref{prop:deltaest} for the last estimate.    
    Note that the assumption from \eqref{eq:F2} along with the boundedness of $\bar z_{\epsilon, \delta}(T)$ and 
    $\bar z_{\epsilon}$ in $\ZZ$ by Lemma~\ref{lem:exRIS} and Corollary~\ref{cor:boundszeps} implies that
    \begin{equation*}
    \begin{aligned}
        & \|\Drm_z \FF(\bar z_{\epsilon, \delta}(T)) - \Drm_z \FF(\bar z_{\epsilon}(T)) \|_{\ZZ^*} \\
        & \quad \leq \int_0^1 \|\Drm_z^2\FF(\bar z_{\epsilon}(T) + \theta (\bar z_{\epsilon, \delta}(T) - \bar z_{\epsilon}(T))) \,\d\theta 
        \,\| \bar z_{\epsilon, \delta}(T) - \bar z_{\epsilon}(T)\|_{\ZZ} \\
        & \quad \leq C_\epsilon \, \| \bar z_{\epsilon, \delta}(T) - \bar z_{\epsilon}(T)\|_{\ZZ},
    \end{aligned}
    \end{equation*}
    which is used in the estimate above. Thanks to \eqref{eq:distZstarineq}, $(\bar z_\epsilon, \bar\ell_\epsilon)$
    is feasible for all $\delta > 0$ sufficiently small. 
    
    The rest of the proof is entirely along the lines of the proof of Theorem~\ref{thm:ocpapproxeps}. 
    For convenience of the reader, we sketch the arguments. 
    Analogously to \eqref{eq:zepsopt}, the optimality of $(z_{\epsilon, \delta}^*, \ell_{\epsilon, \delta}^*)$
    and the feasibility of $(\bar z_{\epsilon, \delta}, \bar\ell_\epsilon)$ for small values of $\delta$
    imply
    \begin{equation}\label{eq:zepsdeltaopt}
    \begin{aligned}
        j(z_{\epsilon, \delta}^*(T))+\frac{\beta}{2}\,\|\ell_{\epsilon, \delta}^*\|_{H^1(0,T;\VV^*)}^2
        &= J_\epsilon(z_{\epsilon, \delta}^*, \ell_{\epsilon, \delta}^*) \\
        &\leq J_\epsilon(\bar z_{\epsilon, \delta}, \bar\ell_\epsilon) \to J_\epsilon(\bar z_{\epsilon}, \bar\ell_\epsilon)
        \quad \text{as }\delta \searrow 0,
    \end{aligned}
    \end{equation}
    where the last convergence follows from Proposition~\ref{prop:deltastrong} (likewise from Theorem~\ref{thm:exRISeps} and 
    the compact embedding $H^1(0,T;\ZZ) \embed^c C([0,T];\VV)$).
    Therefore, the sequence $(\ell_{\epsilon, \delta})_{\delta > 0}$ is bounded in $H^1(0,T;\VV^*)$ and thus 
    admits a weakly converging subsequence with weak limit $\ell_\epsilon^*$.
    From Theorem~\ref{thm:exRISeps}, it follows that $z_{\epsilon, \delta}^*$ converges weakly in $H^1(0,T;\ZZ)$ 
    to the solution $z_\epsilon^*$ of \eqref{eq:RISeps} with $\ell = \ell_\epsilon^*$.
    By the same arguments that led to \eqref{eq:limitlocstab0}, we see that the initial time constraint is fulfilled 
    by $\ell_\epsilon^*$. Moreover, the weak convergence of $\ell_{\epsilon, \delta}^*$ and $z_{\epsilon, \delta}^*$ 
    in $H^1(0,T;\VV^*)$ and $H^1(0,T;\ZZ)$, respectively, together with the weak continuity of the point evaluation in time
    and \eqref{DFweak} implies that 
    \begin{equation*}
        \Drm_z \II(\ell_{\epsilon, \delta}^*(T), z_{\epsilon, \delta}^*(T)) 
        \weak \Drm_z \II(\ell_\epsilon^*(T), z_\epsilon^*(T)) 
        \quad \text{in } \ZZ^* \text{ as }\delta\searrow 0.
    \end{equation*}
    The weak lower semicontinuity of the distance in combination with the feasibility of
    $(z_{\epsilon, \delta}^*, \ell_{\epsilon, \delta}^*)$ thus gives 
    \begin{equation*}
    \begin{aligned}
        & \dist_{\ZZ^*}(\- \Drm_z\II(\ell_\epsilon^*(T), z_\epsilon^*(T)), \partial\RR(0)) \\
        & \qquad \leq \liminf_{\delta\searrow 0} \dist_{\ZZ^*}(\- \Drm_z\II(\ell_{\epsilon, \delta}^*(T), z_{\epsilon, \delta}^*(T)), \partial\RR(0)) 
        \leq \lim_{\delta\searrow 0} \epsilon^{\frac{1}{4}} + \delta^{\frac{1}{4}} = \epsilon^{\frac{1}{4}}
    \end{aligned}
    \end{equation*}
    and thus, the weak limit $(z_\epsilon^*, \ell_\epsilon^*)$ is feasible for \eqref{eq:optcontroleps}.
    The optimality of the weak limit now follows analogously to \eqref{eq:convobj} by 
    the weak lower semicontinuity of the objective and \eqref{eq:zepsdeltaopt}.
    Strong convergence of the sequence of loads again follows from weak convergence and norm convergence, 
    which in turn follows from the convergence of the objective.
\end{proof}

The estimate in \eqref{eq:distZstarendtime} illustrates why we are considering the end time constraint with the 
$\ZZ^*$-distance instead of the $\VV^*$-distance. 
By modifying the proofs of Corollary~\ref{cor:DzI} and Lemma~\ref{lem:zdelH1Z}, it should well be possible to 
derive an a priori bound on $\Drm_z\II(\bar\ell_\epsilon(T), \bar z_{\epsilon, \delta}(T))$ in $\VV^*$ 
giving in turn weak convergence of that term to $\Drm_z\II(\bar\ell_\epsilon(T), \bar z_{\epsilon}(T))$ in $\VV^*$. 
But, the $\VV^*$-distance is of course not weakly continuous so that weak convergence of the $\VV^*$-distance 
cannot be shown in this way, not to mention strong convergence or even an order of convergence as provided for the 
$\ZZ^*$-distance by Proposition~\ref{prop:deltaest}. 

Finally, we address the case when both, $\epsilon$ and $\delta$, are driven to zero in parallel. 
In contrast to Theorem~\ref{thm:ocpapproxeps}, we study the convergence in physical time based on the 
results at the end of Section~\ref{sec:defparaBV}, see Corollary~\ref{cor:convinphystime}. 
The reason is again the lack of convergence of the term involving the $\VV^*$-distance. 
This distance also enters the parametrization via the vanishing viscosity contact potential, 
see \eqref{eq:defs} and \eqref{eq:frakp}.
Due to the little information about the convergence of this term, we are not able to discuss 
the convergence of the double viscous regularization in the parametrized picture. 
However, in physical time, one obtains the following result:

\begin{theorem}[Convergence of the double viscous regularization]\label{thm:ocpapproxepsdelta}
    Suppose that Assumption~\ref{ass:diffsol} is fulfilled.
    Let null sequences $(\epsilon_n)_{n\in \N}, (\delta_m)_{m\in \N} \subset \R^+$ be given 
    and consider, for each $(n,m) \in \N^2$, 
    a global minimizer $(z_{n, m}^*, \ell_{n, m}^*)$ of \eqref{eq:optcontrolepsdelta} 
    with $\epsilon = \epsilon_n$ and $\delta = \delta_m$. 
    Then there exists a subsequence $(\epsilon_{n_k}, \delta_{m_k})_{k\in \N}$ such that 
    \begin{alignat}{3}
        \ell^*_{n_k, m_k} & \to \ell^* & \quad & \text{in } H^1(0,T;\VV^*),\\
        z^*_{n_ k, m_k}(t) & \weak z^*(t) & \quad & \text{in } \ZZ \quad \forall\, t \in \frakC(\hat t^*) \label{eq:zmnptwise}
    \end{alignat}
    as $k\to \infty$,
    where $z^* \in \frakP(\hat t^*, \hat z^*)$ and $(S^*, \hat t^*, \hat z^*, \ell^*)$ is a global minimizer of \eqref{eq:optcontrol}.
    Herein, $\frakC(\hat t^*)$ again denote the points of continuity of $z^*$, see Lemma~\ref{lem:continuitypts}.
    In addition, there holds that
    \begin{equation}\label{eq:zmnweak}
        z^*_{n_k, m_k}  \weak^* z^* \quad \text{in } L^\infty(0,T;\ZZ).
    \end{equation}
    Moreover, there exists a further subsequence denoted by $(z^*_{n_{k_j}, m_{k_j}})_{j\in \N}$
    and a representative of $z^*$ (denoted by the same symbol) such that
    $z^*_{n_{k_j}, m_{k_j}}  \weak z^*$ in $\ZZ$ for all $t\in [0,T]$ and the limit $z^*$ is a function in $\BV([0,T];\ZZ)$.
    Furthermore, 
    \begin{equation}
        \dot z^*_{n_{k_j}, m_{k_j}}  \weak \dot z^* \quad \text{in } \frakM([0,T];\ZZ),
    \end{equation}
    where the dot refers to the distributional time derivative, see Appendix~\ref{sec:weakBV}.
\end{theorem}

\begin{proof}
    From Theorem~\ref{thm:ocpapproxeps} and Corollary~\ref{cor:convinphystime}, 
    we deduce the existence of a subsequence $(\epsilon_{n_k})_{k\in \N}$ such that 
    $\ell^*_{n_k} \to \ell^*$ in $H^1(0,T;\VV^*)$
    and $z^*_{n_k}(t) \weak z^*(t)$ in $\ZZ$ for all $t\in \frakC(\hat t^*)$ as $k\to \infty$, 
    where $(z^*_{n_k}, \ell^*_{n_k})$ is a solution of \eqref{eq:optcontroleps} with $\epsilon = \epsilon_k$.    
    Let now $t\in \frakC(\hat t^*)$ and $\varepsilon > 0$ be arbitrary. 
    Then, due to the compact embedding $\ZZ\embed^c \VV$, there exists $K \in \N$ such that 
    \begin{equation*}
        \| z^*_{n_k}(t) - z^*(t)\|_\VV + \|\ell^*_{n_k} - \ell^*\|_{H^1(0,T;\VV^*)} \leq \frac{\varepsilon}{2} 
        \quad \forall\, k \geq K.
    \end{equation*}
    Furthermore, by means of Proposition~\ref{prop:ocpconvdelta}, the compact embedding 
    $H^1(0,T;\ZZ) \embed^c C([0,T];\VV)$, and the diagonal method, we find a subsequence $(\delta_{m_k})_{k\in \N}$ 
    of $(\delta_m)$ such that there exists $K' \in \N$ with
     \begin{equation*}
        \sup_{\tau \in [0,T]} \| z^*_{n_k, m_k}(\tau) - z^*_{n_k}(\tau)\|_\VV 
        + \|\ell^*_{n_k, m_k} - \ell_{n_k}^*\|_{H^1(0,T;\VV^*)}
        \leq \frac{\varepsilon}{2} \quad \forall\, k \geq K'.
    \end{equation*}  
    Consequently, it holds that 
    \begin{equation*}
        \| z^*_{n_k, m_k}(t) - z^*(t)\|_\VV + \|\ell^*_{n_k, m_k} - \ell^*\|_{H^1(0,T;\VV^*)} \leq \varepsilon
        \quad \forall\, k \geq \max\{K, K'\}.
    \end{equation*}
    This gives $\ell^*_{n_k, m_k} \to \ell^*$ in $H^1(0,T;\VV^*)$ and $z^*_{n_k, m_k}(t) \to z^*(t)$ in $\VV$
    for every $t\in \frakC(\hat t^*)$. Note that the subsequence $(\epsilon_{n_k}, \delta_{m_k})_{k\in \N}$ 
    is independent of $t$.
    
    To prove the pointwise weak convergence of $z^*_{n_k, m_k}$, first note that the a priori estimate 
    from Lemma~\ref{lem:exRIS} implies the existence of a constant $C>0$ such that
    \begin{equation}\label{eq:zmnbound}
        \sup_{t\in [0,T]} \|z^*_{n_k, m_k}(t)\|_\ZZ \leq C \quad \forall\, k\in \N.
    \end{equation}        
    Thus, for every $t\in \frakC(\hat t^*)$, there is a weakly converging subsequence in $\ZZ$ and, due to 
    $z^*_{n_k, m_k}(t) \to z^*(t)$ in $\VV$, the weak limit must be $z^*(t)$ and is therefore unique, which gives the 
    weak convergence of the whole sequence as claimed. 
    
    In order to prove \eqref{eq:zmnweak}, let $v\in L^1(0,T;\ZZ^*)$ be arbitrary.
    Then \eqref{eq:zmnptwise} along with Lemma~\ref{lem:contset} implies that 
    \begin{equation*}
        \dual{v(t)}{z_{n_k, m_k}^*(t)}_{\ZZ^*, \ZZ} \to 
        \dual{v(t)}{z^*(t)}_{\ZZ^*, \ZZ}  
        \quad \text{f.a.a.\ } t \in (0,T).
    \end{equation*}
    Thanks to \eqref{eq:zmnbound}, Lebesgue's dominated convergence theorem 
    thus implies 
    \begin{equation*}
        \int_0^T \dual{v(t)}{z_{n_k, m_k}^*(t)}_{\ZZ^*, \ZZ}\, \d t
        \to \int_0^T \dual{v(t)}{z^*(t)}_{\ZZ^*, \ZZ}\, \d t
    \end{equation*}
    and, since $v\in L^1(0,T;\ZZ^*)$ was arbitrary, this gives the claim.
    
    Finally, due to its convergence, $(\ell_{n_k, m_k})$ is bounded in $H^1(0,T;\VV^*)$ by a constant $M$ and 
    thus, Theorem~\ref{thm:W11bound} implies the existence of a constant $C_M> 0$ such that
    $\var_\ZZ(z_{n_k, m_k}) = \| z_{n_k, m_k}'\|_{L^1(0,T;\ZZ)}\leq C_M$ for all $k\in \N$.
    Together with the uniform bound from \eqref{eq:zmnbound}, this implies 
    $\|z_{n_k, m_k}\|_{\BV([0,T];\ZZ)} \leq C$ for all $k\in \N$.   
    Therefore, we can apply Helly's selection principle in Hilbert spaces \cite[Theorem~B.5.10]{MR15}, which 
    implies the existence of a further subsequence converging weakly-$*$ in $\ZZ$ everywhere in $[0,T]$
    with a limit in $\BV([0,T];\ZZ)$. Due to \eqref{eq:zmnptwise}, this limit coincides with $z^*$ in $\frakC(\hat t^*)$
    and, since the jump set $[0,T]\setminus \frakC(\hat t^*)$ has zero Lebesgue measure by Lemma~\ref{lem:contset}, 
    the pointwise limit is a representative of $z^*$.
    The weak-$*$ convergence of the distributional time derivatives eventually follows from Lemma~\ref{lem:weakBV} 
    in the appendix.    
\end{proof}

\begin{remark}
    We expect that it should be possible to show that the pointwise limit in $\BV([0,T];\ZZ)$ 
    is a balanced viscosity solution in the sense of \cite[Definition~3.8.10]{MR15}, but this goes beyond 
    the scope of our work and is postponed to future research.
\end{remark}

\begin{appendix}


\section{A Lower Semicontinuity Result}\label{sec:limAF}

\begin{lemma}\label{lem:limAF}
	Let $(z_n)_{n\in\N}\subset H^1(0,T;\ZZ)$ satisfy $z_n(0)=z_0$ for all $n\in\N$. Moreover, assume that
	\begin{align*}
		z_n\rightharpoonup z \text{ in } H^1(0,T;\ZZ) 
	\end{align*} 
	by passing $n\to\infty$ with a limit $z\in H^1(0,T;\ZZ)$.
	Then for all $t\in[0,T]$ and all $v\in L^2(0,T;\ZZ)$, the following convergence results are valid:
	\begin{align*}
		\liminf_{n\to\infty} \int_0^t \langle Az_n(r),z_n'(r)-v(r)\rangle_{\ZZ^*,\ZZ}\,\d r&\geq \int_0^t \langle Az(r),z'(r)-v(r)\rangle_{\ZZ^*,\ZZ}\,\d r\\
		\lim_{n\to\infty} \int_0^t \langle \mathrm{D}_z\FF(z_n(r)),z_n'(r)-v(r)\rangle_{\ZZ^*,\ZZ}\,\d r&= \int_0^t \langle \mathrm{D}_z\FF(z(r)),z'(r)-v(r)\rangle_{\ZZ^*,\ZZ}\,\d r.
	\end{align*}
\end{lemma}

\begin{proof}
	Let $t\in[0,T]$ and $v\in L^2(0,T;\ZZ)$ be fixed but arbitrary. 
	Exploiting the symmetry of $A$, the weak lower semicontinuity of $|\cdot|_\ZZ$ and the assumed (pointwise) weak convergence results in 
	\begin{align*}
		& \liminf_{n\to \infty}\int_0^t \langle Az_n(r),z_n'(r)\rangle_{\ZZ^*,\ZZ}\,\d r\\
		&\quad = \liminf_{n\to \infty} \int_0^t \frac{1}{2}\frac{\d}{\d t}| z_n(r)|_\ZZ^2\,\d r \\
		&\quad = \liminf_{n\to \infty} \frac{1}{2}| z_{n}(t)|_\ZZ^2-\frac{1}{2}| z_0|_\ZZ^2 \\
		&\quad \geq \frac{1}{2} |z(t)|_\ZZ^2-\frac{1}{2}| z(0)|_\ZZ^2
		= \int_0^t \frac{1}{2}\frac{\d}{\d t}| z(r)|_\ZZ^2\,\d r
		= \int_0^t \langle Az(r),{z}'(r)\rangle_{\ZZ^*,\ZZ}\,\d r.
	\end{align*}
    Along with the weak convergence of $A z_n$ in $L^2(0,T;\ZZ^*)$, this	gives the first assertion. 
    
    Due to $H^1(0,T;\ZZ) \embed C([0,T];\ZZ)$, the point evaluation in time as an operator in $H^1(0,T;\ZZ)$ is weakly continuous
    and thus $z_n(t)\rightharpoonup z(t)$ in $\ZZ$ for all $t\in[0,T]$.
    Therefore, Lemma~\ref{lem:convF} yields $\lim_{n\to \infty}\FF(z_n(t))=\FF(z(t))$. Further, applying the 
	mean value theorem and using \eqref{eq:F2} yields for all $n\in\N$ and $r\in[0,T]$ that
	\begin{equation*}
	\begin{aligned}
		\|\mathrm{D}_z\FF(z_n(r))\|_{\VV^*}&\leq \|\mathrm{D}_z\FF(0)\|_{\VV^*}+\int_0^1\|\mathrm{D}_z^2\FF(\tau z_n(r))z_n(r)\|_{\VV^*}\,\d \tau\\
		&\leq \|\mathrm{D}_z\FF(0)\|_{\VV^*}+\int_0^1 C(1+\tau^q\|z_n(r)\|_\ZZ^q)\|z_n(r)\|_\ZZ \,\d \tau\leq C,
	\end{aligned}
	\end{equation*}
	where we exploited that $\|z_n\|_{C([0,T];\ZZ)} \leq C \, H^1(0,T;\ZZ) \leq C$ 
	by the assumed weak convergence of $(z_n)_n$.
    Hence, by $\VV^*\hookrightarrow \ZZ^*$ continuously, we have that 
	$\|\mathrm{D}_z\FF(z_n(r))\|_{\ZZ^*}\leq C$. In combination with the weak-weak convergence of 
	$\mathrm{D}_z\FF$ by \eqref{DFweak} and the pointwise weak convergence of $(z_n)$, 
	Lebesgue's dominated convergence theorem thus yields
	\begin{align*}
		\lim_{n\to \infty} \int_0^t \langle \mathrm{D}_z\FF(z_n(r)),v(r)\rangle_{\ZZ^*,\ZZ}\, \d r=\int_0^t\langle \mathrm{D}_z\FF(z(r)),v(r)\rangle_{\ZZ^*,\ZZ} \,\d r
	\end{align*}
	Altogether, we have shown that
	\begin{align*}
		& \lim_{n\to \infty}\int_0^t \langle \mathrm{D}_z\FF(z_n(r)),z_n'(r)-v(r)\rangle_{\ZZ^*,\ZZ}\, \d r\\
		&\quad = \lim_{n\to \infty}\FF(z_n(t))-\FF(z_0)-\int_0^t\langle \mathrm{D}_z\FF(z_n(r)),v(r)\rangle_{\ZZ^*,\ZZ} \,\d r\\
		&\quad = \FF(z(t))-\FF(z_0)-\int_0^t\langle \mathrm{D}_z\FF(z(r)),v(r)\rangle_{\ZZ^*,\ZZ} \,\d r\\
		&\quad = \int_0^t \langle \mathrm{D}_z\FF(z(r)),{z}'(r)-v(r)\rangle_{\ZZ^*,\ZZ} \,\d r,
	\end{align*}
	which completes the proof.
\end{proof}

\section{Differentiability of the $\epsilon, \delta$-Norm}\label{sec:normdiff}

\begin{lemma}\label{lem:diffnorm}
    Let $v\in H^{1}(0,T;\ZZ)$ be given. Then the mapping $[0,T]\ni t\mapsto \|v(t)\|_{\epsilon,\delta}\in \R$ 
    is differentiable fa.a.\ $t\in(0,T)$ and almost everywhere in $\omega:=\{ t\in[0,T]: \|v(t)\|_{\epsilon,\delta} > 0\}$ 
    the derivative is 
	\begin{equation}\label{eq:diffnorm}
		\frac{\d}{\d t}\|v(t)\|_{\epsilon,\delta}=\bigg\langle v'(t),\frac{v(t)}{\|v(t)\|_{\epsilon,\delta}}\bigg\rangle_{\epsilon,\delta}
	\end{equation}
	Moreover, it holds
	\begin{equation}\label{eq:HDInorm}
		\int_\omega \frac{\d}{\d t}\|v(t)\|_{\epsilon,\delta} \,\d t= \|v(T)\|_{\epsilon,\delta} - \|v(0)\|_{\epsilon,\delta}.
	\end{equation}
\end{lemma}

\begin{proof}
We abbreviate $f(t):= \|v(t)\|_{\epsilon, \delta}$. For all $s,t\in[0,T]$, we have 
\begin{equation*}
	|f(t)-f(s)| 
	\leq \|v(t)-v(s)\|_{\epsilon,\delta}
	\leq  \max\{1, \delta/\epsilon\}\,\|v(t)-v(s)\|_\ZZ
\end{equation*}
such that the absolute continuity of $v$ transfers to $f$ and hence, $f$ is differentiable almost everywhere as claimed. 
For every point $t\in \omega$ of differentiability of $v$, we obtain
\begin{align*}
	&\bigg|\frac{f(t+h)-f(t)}{h}-\bigg\langle v'(t),\frac{v(t)}{\|v(t)\|_{\epsilon,\delta}}\bigg\rangle_{\epsilon,\delta}\bigg|\\
	& \quad \leq 
    \begin{aligned}[t]
        \bigg|\frac{\|v(t)+h\, v'(t)\|_{\epsilon,\delta} - \|v(t)\|_{\epsilon,\delta}}{h}
        & - \bigg\langle v'(t),\frac{v(t)}{\|v(t)\|_{\epsilon,\delta}}\bigg\rangle_{\epsilon,\delta}\bigg| \\
	    & + \frac{\|v(t+h)-v(t)-h\,v'(t)\|_{\epsilon,\delta}}{h} \, . 
    \end{aligned}
\end{align*}
Now, letting $h\searrow 0$, 
the first term on the right hand side tends to zero because of the differentiability of the norm outside of zero, see
\cite[Proposition~4.7.10]{Sch07}.
The second term vanishes, too, since $v$ is differentiable in $t$. As this is the case in almost every $t\in (0,T)$, 
we have proven \eqref{eq:diffnorm}.

In order to prove \eqref{eq:HDInorm}, we assume that $[0,T]\setminus \omega$ has positive Lebesgue measure. 
Otherwise, the claim follows directly from the fundamental theorem of calculus for absolutely continuous functions. 
The idea is now to split $\omega$ into intervals such that the fundamental theorem of calculus is also applicable. 
Therefore, we define
\begin{equation*}
	t_1:=\inf \{t\in[0,T]:f(t)=0\} \in [0,T],\quad t_2:=\sup \{t\in[0,T]:f(t)=0\} \in [0,T].
\end{equation*}
From the continuity of $f$ we deduce $0=f(t_1)=f(t_2)$. Moreover, by continuity of $v$, the set 
\begin{equation*}
	S:=\omega\cap \Big( [0,T]\setminus \big ([0,t_1]\cup [t_2,T]\big)\Big)
\end{equation*}  
is open. By means of the representation theorem for open sets \cite[Lemma~5.7.1]{Ber99}, 
$S$ can be uniquely written as the union 
of countable many disjoint open intervals, i.e., $S=\bigcup_{k\in\N} (s_{2k-1},s_{2k})$ with $t_1\leq s_1<s_2<...\leq t_2$. 
We directly obtain  
$f(s_k)=0$ for all $k\in \N$ because otherwise $s_k\in S$, which contradicts the disjointness of the intervals.
Eventually, from the fundamental theorem of calculus for absolutely continuous functions, we infer
\begin{align*}
	\int_\omega \frac{\d}{\d t} f(t) \d t
	&=\int_0^{t_1} \frac{\d}{\d t} f(t) \d t+\int_{t_2}^T \frac{\d}{\d t} f(t) \d t+\sum_{k\in\N} \int_{s_{2k-1}}^{s_{2k}} \frac{\d}{\d t} f(t) \d t\\
	&= f(t_1)-f(0)+f(T)-f(t_2)+\sum_{k\in\N}f(s_{2k})-f(s_{2k-1})
	=f(T)-f(0).
\end{align*}
Note that this equation is also valid for $t_1=0$, $t_2=T$ or $S=\emptyset$. All in all, this proves \eqref{eq:HDInorm}.
\end{proof}

\section{Continuous BV Solutions are Differential Solutions}\label{sec:diffsol}

\begin{lemma}
    Assume that $\hat t$ satisfies Assumption~\ref{ass:diffsol}. Then 
    $\tilde z := \hat z \circ \hat t^{-1}$  is an element of $H^1(0,T;\ZZ)$ and 
    a differential solution of \eqref{eq:ris}, i.e., it satisfies 
    \begin{equation}\label{eq:RIStilde2}
	    0\in\partial\RR(\tilde z'(t))+\mathrm{D}_z\II(\ell(t),\tilde z(t)) \; \text{ a.e.\ in }(0,T), \quad \tilde z(0)=z_0.
    \end{equation}
\end{lemma}

\begin{proof}
    First of all, by the chain rule in Sobolev-Bochner spaces, we obtain that $\tilde z\in H^1(0,T;\ZZ)$ with 
    derivative $\tilde z'(t)= \hat z'(\hat t^{-1}(t))\frac{\d}{\d t}\hat t^{-1}(t)$ a.e.\ in $(0,T)$, see also \cite[Lemma~A.1]{KT23}.
    Due to \eqref{eq:t'>del} and the complementarity condition \eqref{eq:compl}, we have that 
    \begin{equation}\label{eq:distloc}
        	\dist_{\VV^*}(-\mathrm{D}_z\II(\hat\ell(s),\hat z(s)),\partial\RR(0))=0 \quad \text{f.a.a.\ } s\in (0,S).
    \end{equation}
    Therefore, the energy identity \eqref{eq:ener} can be written as
    \begin{equation}\label{eq:energeqdiffsol}
    \begin{aligned}
        & \II(\hat \ell(s),\hat z(s))+\int_0^s \RR(\hat z'(r))\, \d r \\
        & \qquad =\II(\hat\ell(0),z_0)-\int_0^s \langle \hat\ell'(r),\hat z(r)\rangle_{\VV^*,\VV}\, \d r    
        \quad \forall\, s\in [0,S],
    \end{aligned}
    \end{equation}
    where we used that the generalized metric derivative equals $\RR(\hat z'(s))$ in all points $s\in (0,S)$ 
    of differentiability of $\hat z$
    (which are almost all points in $(0,S)$), cf.\ also \cite[Lemma~B.3 and B.4]{And25}.
    The assertion now follows from transferring \eqref{eq:distloc} and \eqref{eq:energeqdiffsol} into physical time.
    From \eqref{eq:distloc} we infer
    \begin{equation}\label{eq:DIR(0)}
	    -\mathrm{D}_z\II(\ell(t),\tilde z(t)) \in \partial\RR(0)\quad \text{f.a.a. }t\in(0,T).
    \end{equation}
    In view of $\hat \ell =\ell\circ \hat t$ and \eqref{R2}, applying the substitution $s=\hat t^{-1}(t)$ 
    to \eqref{eq:energeqdiffsol} yields along with the chain rule for the energy according to \cite[Proposition~E.1]{KT23} 
    or \cite[Lemma~C.3]{And25} that
    \begin{align*}
	    \int_0^t \RR(\tilde z'(r))\, \d r&= \II(\ell(0),z_0)-\II(\ell(t),\tilde z(t))-\int_0^t \langle \ell'(r),\tilde z(r)\rangle_{\VV^*,\VV} \,\d r\\
	    &=\int_0^t -\frac{\d}{\d t} \II(\ell(r),\tilde z(r)) -\langle \ell'(r),\tilde z(r)\rangle_{\VV^*,\VV} \,\d r\\
	    &=\int_0^t -\langle \mathrm{D}_z\II(\ell(r),\tilde z(r)),\tilde z'(r)\rangle_{\ZZ^*,\ZZ} \,\d r.
    \end{align*}
    Since this holds for all $t\in [0,T]$, we obtain 
    $\RR(\tilde z'(t))=-\langle \mathrm{D}_z\II(\ell(t),\tilde z(t)),\tilde z'(t)\rangle_{\ZZ^*,\ZZ}$ f.a.a.\ $t\in(0,T)$.
    In combination with \eqref{eq:DIR(0)} and the positive homogeneity of $\RR$, this implies
    that \eqref{eq:RIStilde2} is fulfilled so that $\tilde z$ is indeed a differential solution.
\end{proof}

\section{Uniqueness of Differential Solutions}\label{sec:diffuni}

\begin{lemma}
    Let $z_1, z_2\in H^1(0,T;\ZZ)$ be two differential solutions of \eqref{eq:etaris}
    that, in addition, satisfy the estimate in \eqref{eq:uniboundeta} with $\eta = \bar\eta$. 
    Then $z_1 = z_2$. 
\end{lemma}

\begin{proof}
    By definition of the convex subdifferential, $z_1$ and $z_2$ fulfill
    \begin{equation*}
        \RR(v)\geq \RR(z_i'(t)) + 
        \langle -\mathrm{D}_z\II_{\bar\eta}(\ell(t) + \bar\eta\, \tilde z(t),z_i(t)),v-z_i'(t)\rangle_{\ZZ^*,\ZZ}, \quad i = 1,2,
    \end{equation*}
    for all $v\in \ZZ$ and almost all $t\in(0,T)$. Now inserting $v=z_2'(t)$ in the inequality for $z_1$ and vice versa
    and adding up both inequalities results in 
	\begin{equation}\label{eq:estDE<0}
		0
		\geq \langle \mathrm{D}_z\EE_{\bar\eta}(z_1(t))-\mathrm{D}_z\EE_{\bar\eta}(z_2(t)), z_1'(t)-z_2'(t)\rangle_{\ZZ^*,\ZZ}
	\end{equation}
	for almost all $t\in(0,T)$, where we used that the loads enter the energy just linearly. 
	Define now the function $\mu \in H^1(0,T)$ by
	\begin{equation*}
	    \mu(t):= \langle \mathrm{D}_z\EE_{\bar\eta}(z_1(t))-\mathrm{D}_z\EE_{\bar\eta}(z_2(t)),z_1(t)-z_2(t)\rangle_{\ZZ^*,\ZZ}.
    \end{equation*}	 
    Then the mean value theorem along with the uniform convexity of the energy by Remark~\ref{rem:Eetaconvex} yields
	\begin{equation}
	\begin{aligned}\label{eq:estmu}
        \mu(t) &= \int_0^1 \langle \mathrm{D}_z^2\EE_{\bar\eta}(z_1(t)+\tau (z_2(t)-z_1(t)))(z_1(t)-z_2(t)),z_1(t)-z_2(t)\rangle_{\ZZ^*,\ZZ} \,\d \tau\\
        &\geq \frac{\alpha}{2}\, \|z_1(t)-z_2(t)\|_\ZZ^2.
	\end{aligned}
	\end{equation}
	Note that any convex combination of $z_1$ and $z_2$ of course also satisfies the estimate in \eqref{eq:uniboundeta} with $\eta = \bar\eta$.
	For the derivative of $\mu$, one obtains
    \begin{align*}
        \mu'(t)
        &= \langle \mathrm{D}_z^2\EE_{\bar\eta}(z_1(t))z_1'(t)-\mathrm{D}_z^2\EE_{\bar\eta}(z_2(t))z_2'(t),z_1(t)-z_2(t)\rangle_{\ZZ^*,\ZZ}\\
        &\quad + \langle \mathrm{D}_z\EE_{\bar\eta}(z_1(t))-\mathrm{D}_z\EE_{\bar\eta}(z_2(t)),z_1'(t)-z_2'(t)\rangle_{\ZZ^*,\ZZ}\\
		&= \langle \mathrm{D}_z^2\FF(z_1(t))(z_1(t)-z_2(t))+ \mathrm{D}_z\FF(z_2(t))-\mathrm{D}_z\FF(z_1(t)),z_1'(t)\rangle_{\ZZ^*,\ZZ}\\
		&\quad + \langle \mathrm{D}_z^2\FF(z_2(t))(z_2(t)-z_1(t)) + \mathrm{D}_z\FF(z_1(t))-\mathrm{D}_z\FF(z_2(t)),z_2'(t)\rangle_{\ZZ^*,\ZZ}\\
        &\quad +2 \langle \mathrm{D}_z\EE_{\bar\eta}(z_1(t))-\mathrm{D}_z\EE_{\bar\eta}(z_2(t)),z_1'(t)-z_2'(t)\rangle_{\ZZ^*,\ZZ}.
	\end{align*}
	Using the mean value theorem again, the first term on the right hand side is estimated by
    \begin{align*}
		& \big|\langle \mathrm{D}_z^2\FF(z_1(t))(z_1(t)-z_2(t))+ \mathrm{D}_z\FF(z_2(t))-\mathrm{D}_z\FF(z_1(t)),z_1'(t)\rangle_{\ZZ^*,\ZZ}\big|\\
		& \quad \leq \int_0^1\big\|[\mathrm{D}_z^2\FF(z_1(t))- \mathrm{D}_z^2\FF(z_1(t)+\tau (z_2(t)-z_1(t)))](z_1(t)-z_2(t))\big\|_{\ZZ^*}\,\d \tau\, \|z_1'(t)\|_\ZZ\\
		& \quad \leq \frac{L_r}{2}\, \|z_1(t)-z_2(t)\|_\ZZ^2\, \|z_1'(t)\|_\ZZ,
	\end{align*}
	where we used the assumption from \eqref{eq:D2Flip} with $r:= C(\bar\eta+1)$. 
	The second term on the right hand side can be estimated in exactly the same way and thus, \eqref{eq:estDE<0} implies
    \begin{align*}
		\mu'(t) \leq \frac{L_r}{2}\,\|z_1(t)-z_2(t)\|_\ZZ^2 \big( \|z_1'(t)\|_\ZZ+\|z_2'(t)\|_\ZZ \big)
		\leq \frac{L_r}{\alpha}\, \big(\|z_1'(t)\|_\ZZ+\|z_2'(t)\|_\ZZ\big)\mu(t),
	\end{align*}
	where we used \eqref{eq:estmu} for the last inequality. Since $z_1(0) = z_2(0) = z_0$ and thus $\mu(0)=0$ , integrating gives
	\begin{align*}
		\mu(t)\leq  \int_0^t\frac{L_r}{\alpha}\,\big(\|z_1'(r)\|_\ZZ+\|z_2'(r)\|_\ZZ\big) \mu(r)\,\d r
	\end{align*}
    Eventually, applying Gronwall's lemma yields $\frac{\alpha}{2}\, \|z_1(t)-z_2(t)\|_\ZZ^2\leq  \mu(t)\leq 0$
	such that $z_1=z_2$ as claimed.
\end{proof}
    
\section{Weak-$*$ Convergence in $\BV([0,T];\ZZ)$}\label{sec:weakBV}

Let a function $\zeta \in \BV([0,T];\ZZ)$ be given. Then, for every function $v\in C([0,T];\ZZ)$, 
the Riemann-Stieltjes integral is well defined, which we denote by 
$\int_0^T \dual{v(t)}{\d \zeta(t)}_\ZZ$.
For the precise definition of the Riemann-Stieltjes integral for vector-valued functions as well as 
the properties thereof used in the following, we refer to \cite[Section~V.1]{Kre96}.
The mapping 
\begin{equation*}
    C([0,T];\ZZ) \ni v\mapsto \int_0^T \dual{v(t)}{\d \zeta(t)} \in \R, 
\end{equation*}
is a bounded linear form and thus an element of the dual space $C([0,T];\ZZ)^*$, 
which can be identified with $\frakM([0,T];\ZZ)$, the space of $\ZZ$-valued regular Borel measures, 
by means of the Riesz representation theorem for vector-valued measures, see, e.g., 
\cite{Zin57}, \cite[Section~6.5]{DU77}. We denote this measure by $\dot \zeta \in \frakM([0,T];\ZZ)$, since due to 
the formula of integration by parts for the Riemann-Stieltjes integral, there holds 
for every $\psi\in C^\infty_c(0,T)$ and every $v\in \ZZ$ that 
\begin{equation*}
\begin{aligned}
    \Big\langle \int_{[0,T]} \psi(t) \, \d \dot \zeta(t) , v\Big\rangle_\ZZ
    = \int_0^T \dual{\psi(t) v}{\d \zeta(t)}_\ZZ 
    = - \Big\langle \int_0^T \psi'(t)\, \zeta(t)\,\d t , v\Big\rangle_\ZZ
\end{aligned}
\end{equation*}
and so, $\dot\zeta$ is just the distributional time derivative of $\zeta$.
As usual, we say that a sequence $(\mu_n)_{n\in \N} \subset \frakM([0,T];\ZZ)$ converges weakly-$*$ to a limit 
$\mu \in \frakM([0,T];\ZZ)$, if 
\begin{equation*}
    \int_{[0,T]} v\,\d \mu_n \to \int_{[0,T]} v\,\d\mu \quad \forall\, v\in C([0,T];\ZZ).
\end{equation*} 

\begin{lemma}\label{lem:weakBV}
    Assume that a sequence $\{z_n\}_{n\in \N} \subset \BV([0,T];\ZZ)$ satisfies 
    \begin{equation}\label{eq:varbound}
        \sup_{n\in \N} \| z_n \|_{\BV([0,T];\ZZ)} \leq C 
    \end{equation}
    and 
    \begin{equation}
        z_n(t) \weak z(t) \quad \text{in } \ZZ \quad \forall\, t\in [0,T]
    \end{equation}
    with a limit $z\in \BV([0,T];\ZZ)$. Then 
    \begin{equation*}
        \dot z_n \weak^* \dot z \quad \text{in }\frakM([0,T];\ZZ).
    \end{equation*}
\end{lemma}

\begin{proof}
   Let first $v\in C^1([0,T];\ZZ)$ be arbitrary. Then, due to the pointwise weak convergence and 
   the uniform boundedness of $z_n$ in $\BV([0,T];\ZZ)$, which in particular implies that 
   $\sup_{n\in\N}\sup_{t\in [0,T]} \|z_n(t)\|_\ZZ \leq C$, Lebesgue's dominated convergence theorem implies   
   \begin{equation*}
       \int_0^T \dual{v'(t)}{z_n(t)}_\ZZ\,\d t \to \int_0^T \dual{v'(t)}{z(t)}_\ZZ\,\d t .
   \end{equation*}
   Using the formula of integration for the Riemann-Stieltjes integral and the pointwise weak convergence again, 
   we thus obtain 
    \begin{equation}\label{eq:convstieltjes}
    \begin{aligned}
        & \int_0^T \dual{v(t)}{\d z_n(t)}_{\ZZ} \\[-1ex]
        &\quad = \dual{z_n(T)}{v(T)}_{\ZZ} - \dual{z_n(0)}{v(0)}_{\ZZ} 
        - \int_0^T \dual{v'(t)}{z_n(t)}_\ZZ\,\d t \\
        &\qquad \to \dual{z(T)}{v(T)}_{\ZZ} - \dual{z(0)}{v(0)}_{\ZZ} 
        - \int_0^T \dual{v'(t)}{z(t)}_\ZZ\,\d t \\
        & \qquad\qquad = \int_0^T \dual{v(t)}{\d z(t)}_{\ZZ}.
    \end{aligned}
    \end{equation}
    Next, let $\varphi \in C([0,T];\ZZ)$ and $\varepsilon > 0$ be arbitrary. Since $C^1([0,T];\ZZ)$ is dense in $C([0,T];\ZZ)$
    by the Weierstrass approximation theorem, see, e.g., \cite[Satz~IV.1.3]{GGZ74}, 
    there is a function $v\in C^1([0,T];\ZZ)$ such that 
    \begin{equation*}
         \| \varphi - v\|_{C([0,T];\ZZ)} \leq \frac{\varepsilon}{2 (C + \var_\ZZ(z))},
    \end{equation*}
    where $C$ is the bound from \eqref{eq:varbound}. Moreover, by \eqref{eq:convstieltjes},
    there is an index $N\in \N$ such that $|\int_0^T \dual{v}{\d (z_n -z)} |\leq \varepsilon/2$
    for all $n\geq N$. Consequently,  
    \begin{equation*}
    \begin{aligned}
        & \Big| \int_0^T \dual{\varphi(t)}{\d (z_n - z)(t)}_\ZZ\Big|  \\
        & \quad \leq \Big| \int_0^T \dual{v(t)}{\d(z_n - z)(t)} \Big|
        + \| \varphi - v\|_{C([0,T];\ZZ)} \big(\var_\ZZ(z_n) + \var_\ZZ(z)\big) 
        \leq \varepsilon
    \end{aligned}        
    \end{equation*}
    for all $n \geq N$, which implies 
    \begin{equation}\label{eq:weakconvstieltjes}
        \int_0^T \dual{\varphi(t)}{\d z_n(t)}_\ZZ \to \int_0^T \dual{\varphi(t)}{\d z(t)}_\ZZ
        \quad \, \forall\, \varphi \in C([0,T];\ZZ).
    \end{equation}
    The identification of the Riemann-Stieltjes integral with the distributional time derivative in $\frakM([0,T];\ZZ)$ 
    by means of the Riesz representation theorem mentioned above then implies the result.
\end{proof}

\end{appendix}

\renewcommand*{\bibfont}{\footnotesize}
\printbibliography

\end{document}